\documentclass[a4paper,oneside,final,notitlepage,onecolumn]{article}
\usepackage{makeidx}
\usepackage{amsmath}
\usepackage{latexsym}
\usepackage{eepic}
\usepackage{epsfig}
\usepackage{amsfonts}
\usepackage{amssymb}
\evensidemargin\oddsidemargin
\numberwithin{equation}{section}
\newtheorem{theorem}{Theorem}[section]

\newtheorem{lemma}[theorem]{Lemma}
\newtheorem{proposition}[theorem]{Proposition}
\newtheorem{definition}[theorem]{Definition}
\newtheorem{question}[theorem]{Question}
\newtheorem{example}[theorem]{Example}

\newtheorem{remark}[theorem]{Remark}

\def\deltafett{\mbox{\boldmath{$\delta$}}}

\begin{document}

\title{\textbf{All Toric L.C.I.-Singularities \\ Admit
    Projective Crepant Resolutions}}
\author{\textbf{Dimitrios I. Dais}$^{\ast}$\textbf{, Christian
    Haase}$^{\diamond}%
$\textbf{\ and G\"{u}nter
  M. Ziegler}$^{\diamond}\medskip$\\\noindent$^{\ast}${\scriptsize
  Mathematische Fakult\"{a}t, Universit\"{a}t
T\"{u}bingen, Auf der Morgenstelle 10, D-72076 T\"{u}bingen,
Germany}\\{\scriptsize e-mail:
dais@wolga.mathematik.uni-tuebingen.de\smallskip}\\$^{\diamond}${\scriptsize
Fachbereich Mathematik, TU-Berlin, Sekr. MA 7.1,
Stra\ss e des 17. Juni 136, D-10623 Berlin, Germany}
\and {\scriptsize e-mail: [haase, ziegler]@math.tu-berlin.de}}
\date{}
\maketitle
\begin{abstract}
\noindent It is known that the underlying spaces of all abelian quotient
singularities which are embeddable as complete intersections of hypersurfaces
in an affine space can be overall resolved by means of projective
torus-equivariant crepant birational morphisms in all dimensions. In the
present paper we extend this result to the entire class of toric
l.c.i.-singularities. Our proof makes use of Nakajima's classification theorem
and of some special techniques from toric and discrete geometry.
\end{abstract}

\section{Introduction\label{INTRO}}

\noindent In the past two decades ``crepant'' birational morphisms were mainly
used in algebraic geometry to reduce the canonical singularities of algebraic
(not necessarily proper) $d$-folds, $d\geq3$, to $\mathbb{Q}$-factorial
terminal singularities, and to treat \textit{minimal models} in high
dimensions. From the late eighties onwards, crepant desingularizations
$\widehat{Y}\longrightarrow Y$ of projective varieties $Y$ with trivial
dualizing sheaf and well-controllable singularities play also a crucial role
in producing \textit{Calabi-Yau manifolds}, which serve as internal target
spaces for ``non-linear supersymmetric sigma models'' in the framework of the
physical string-theory. This explains the recent mathematical interest in both
local and global versions of the existence-problem of \textit{smooth}
birational models of such $Y$'s.\medskip\newline $\bullet$\textbf{\
  The local
problem. }This was at first focused on the\textbf{\ }high-dimensional ``McKay
correspondence'' for the underlying spaces $\mathbb{C}^{d}/G$, $G\subset$
SL$\left(  d,\mathbb{C}\right)  $, $d\geq2$, of the Gorenstein quotient
singularities, connecting the irreducible representations of $G$ (or dually,
the conjugacy classes of $G)$, on the one hand, and the cohomology ring of the
overlying spaces $\widehat{X}$ of (preferably projective), crepant, full
desingularizations $\widehat{X}\longrightarrow X=\mathbb{C}^{d}/G$ of $X$, on
the other (cf. \cite{Bat3, Reid3}). The problem setting was partially
extended by proving that a one-to-one correspondence of McKay-type is true,
too, for torus-equivariant, crepant, full desingularizations $\widehat
{X}\longrightarrow X=U_{\sigma}$ of the underlying spaces of all Gorenstein
toric singularities \cite[\S4]{BD}. As it turned out, the non-trivial (even)
cohomology groups of $\widehat{X}$'s have the ``expected'' dimensions,
depending on the ``weights'' (also called ``ages'' in \cite{Ito-Reid}) of the
conjugacy classes of the acting groups, and on the Ehrhart polynomial of the
corresponding lattice polytopes, respectively, and being therefore
\textit{independent} of particular choices of representatives among these
$\widehat{X}$'s. (See \cite[Thm. 8.4]{Bat3} and \cite[Thm. 4.4]{BD}). We
should particularly note that, in both cases, $\widehat{X}$'s of this sort
always exist in dimension $d\leq3$. However, if $d\geq4$, \ this is not always
possible because even within their common class
\begin{equation}
\left\{
\begin{array}
[c]{c}%
\text{Gorenstein abelian}\\
\text{quotient singularities}%
\end{array}
\right\}  =\left\{
\begin{array}
[c]{c}%
\text{Gorenstein quotient}\\
\text{singularities}%
\end{array}
\right\}  \ \bigcap\ \left\{
\begin{array}
[c]{c}%
\text{Gorenstein toric}\\
\text{ singularities}%
\end{array}
\right\}  \label{CLA}%
\end{equation}
there are mostly terminal singularities (see Morrison-Stevens \cite{Mo-Ste}).
Hence, the high-dimensional bijections of McKay-type make sense only in
combination with the following:

\begin{question}
[Existence Problem]\label{QUESTION1}Under which conditions \emph{(}or
restrictions\emph{)} on our starting-point data for these two classes of
Gorenstein singularities do projective, crepant, full
resolutions exist in dimensions $\geq4$ \emph{?}
\end{question}

\noindent First answers via toric
geometry for the case of Gorenstein abelian quotient singularities
(\ref{CLA}) were given in \cite{DHH,DH,DHZ1,DHZ2}. To our surprise,
the number of families of these singularities admitting resolutions of such
special kind is \textit{not} negligible as one would at first sight
expect.\medskip\newline $\bullet$\textbf{\ The global problem. }A wide class of
CY-manifolds of particular interest is that one being constructible by
resolving overall the (necessarily Gorenstein) singularities of the so-called
\textit{CY-varieties} via suitable projective crepant morphisms. In the case
in which the singularities of such a CY-variety $Y$ are of ``mild nature''
(like quotient or toroidal singularities), and as long as an appropriate
stratification of the singular locus Sing$\left(  Y\right)  $ of $Y$ is
available, the existence-problem of crepant full resolutions can be mostly
reduced to the local one by performing standard glueing procedures. (In
contrast to this, the conditions which would guarantee the preservation of the
projectivity of the desingularizing morphisms globally are much more
complicated and require additional information about the global geometry
of $Y$). It is worth mentioning that also the Hodge numbers $h^{p,q}%
(\widehat{Y})$ of the overlying spaces of \textit{all} crepant, full, global
desingularizations $\widehat{Y}\longrightarrow Y$ of $Y$ remain
\textit{invariant} (see Kontsevich \cite{Kontsevich} and Batyrev \cite{Bat2}).
A method of working \textit{formally} with $Y$'s, even without assuming the
existence of such special $\widehat{Y}$'s in dimensions $\geq4$, consists in
introducing the so-called \textit{string-theoretic }Hodge numbers\textit{\ }%
$h_{\text{str}}^{p,q}(Y)$ for $Y$'s (cf. \cite{BD,BB})\textit{. }Since
as yet there is only a conjectural description of a candidate for the
cohomology complex which probably leads to a mathematical definition of the
``string-theoretic cohomology theory'' globally\textit{\ }(cf. \cite[4.4]{BB}
and Borisov's new approach in \cite[\S4]{Borisov1},
\cite[Conj. 9.23]{Borisov2}), it would be important to know at least
\textit{when} the existence of \textit{smooth} birational models for
$Y$'s is feasible or not.\medskip

\noindent$\bullet$ \textbf{The L.C.I.'s.} In the present paper we shall
exclusively deal with one aspect of the local problem. We believe that
a purely algebraic, \textit{sufficient} condition for the existence of the
desired resolutions in \textit{all }dimensions is to require from our
singularities to be, in addition, \textit{l.c.i.'s}. In the toric category,
where the Question \ref{QUESTION1} can be translated into a question
concerning the existence of specific lattice triangulations of lattice
polytopes, this conjecture was verified for abelian quotient singularities in
\cite{DHZ1} via Kei-ichi Watanabe's Theorem \cite{Watanabe}. (For non-abelian
groups acting on $\mathbb{C}^{d}$, it remains open). Furthermore, the authors
of \cite[cf. \S8(iii)]{DHZ1} asked for geometric analogues of the ``joins''
and ``dilations'' occuring in their Reduction Theorem also for toric
non-quotient l.c.i.-singularities. As we shall see below, such a
characterization (in a somewhat different context) \textbf{is} indeed possible
by making use of another beautiful classification theorem due to Haruhisa
Nakajima \cite{Nakajima}, which generalizes Watanabe's results to the entire
class of toric l.c.i.'s. Based on this
classification we prove the following:

\begin{theorem}
[Main Theorem]\label{MAIN}The underlying spaces of all toric
l.c.i.-singularities admit torus-equivariant, projective, crepant, full
resolutions \emph{(}i.e., \emph{``}smooth minimal models\emph{'') }in all dimensions.
\end{theorem}

\noindent The proof of \ref{MAIN} relies on considerably simpler techniques
than those of \cite{DHZ1}, basically because the vertices of the Nakajima's
polytopes reside in the standard \textit{rectangular} lattice within
$\mathbb{R}^{d}$. Nevertheless, Watanabe's forests and \textit{skew} lattices
remain the right language if one wishes to read off the \ weights
of\ \textit{abelian group} actions by predeterminated eigencoordinates and
diagonalizations in a direct manner. On the other hand, the common distinctive
feature in both proofs is an \textit{inductive} argument which makes things
work in \textit{all} dimensions.\medskip

\noindent$\bullet$ This paper is organized as follows: After recalling the
algebraic hierarchy of singularities (see (\ref{HIERARCHY}) below), and some
basic notions and facts from toric geometry in \S\ref{TORIC-G}, we explain in
\S\ref{BC-TR} why the existence of the desired desingularizations is
equivalent to the existence of b.c.-triangulations of the lattice polytopes
supporting the Gorenstein cones. Moreover, we give two first examples of
lattice polytopes (namely the so-called Fano and $\mathbb{H}_{d}$-compatible
polytopes) admitting such triangulations, and describe the corresponding
exceptional prime divisors explicitly. In section \ref{NAK} we provide
convenient reformulations of Nakajima's classification. In \S\ref{M-KOS} we
give the proof of Main Theorem \ref{MAIN} by using certain maximal coherent
triangulations, combined with the ``Key-Lemma'' \ref{KEY} which guarantees
their ``basicness''. An immediate algebraic application of \ref{MAIN} is
contained in the second part of section \ref{M-KOS}, where it is shown that
the monoidal ``coordinate rings'' $\mathbb{C}\left[  \tau_{P}\cap
\mathbb{Z}^{d}\right]  $ of $U_{\tau_{P}^{\vee}}$'s for all Nakajima polytopes
$P$ have the Koszul-property. In \S\ref{COHOMOLOGY} we present a simple method
of computing the non-trivial cohomology group dimensions of the overlying
spaces of all crepant, full resolutions of toric l.c.i.-singularities.
Finally, in section \ref{EXTREME} we apply our results for two ``extreme''
classes of toric g.c.i.-singularities which occur as direct generalizations of
the classical $A_{k-1}$-singularities in arbitrary dimensions.\bigskip
\ \newline $\bullet$ \textbf{General terminology.} \textsf{(a)} First we
recall some fundamental definitions from commutative algebra
(cf. \cite{Kunz, Matsumura}). Let $R$ be a commutative ring with
$1$. The \textit{height}
ht$\left(  \mathfrak{p}\right)  $ of a prime ideal $\mathfrak{p}$ of $R$ is the
supremum of the lengths of all prime ideal chains which are contained in
$\mathfrak{p}$, and the \textit{dimension} of $R$ is defined to be dim$\left(
R\right)  :=$ sup$\left\{  \text{ht}\left(  \mathfrak{p}\right)  \left|
\mathfrak{p}\text{ prime ideal of }R\right.  \right\}  $. $R$ is
\textit{Noetherian} if any ideal of it has a finite system of generators. $R$
is a \textit{local ring} if it is endowed with a \textit{unique} maximal ideal
$\mathfrak{m}$. A local ring $R$ is \textit{regular }(resp. \textit{normal}) if
dim$\left(  R\right)  =$ dim$\left(  \mathfrak{m}/\mathfrak{m}^{2}\right)  $ (resp. if
it is an integral domain and is integrally closed in its field of fractions).
A finite sequence $a_{1},\ldots,a_{\nu}$ of elements of a ring $R$ is defined
to be a \textit{regular sequence} if $a_{1}$ is not a zero-divisor in $R$ and
for all $i$, $i=2,\ldots,\nu$, $a_{i}$ is not a zero-divisor of
$R/\left\langle a_{1},\ldots,a_{i-1}\right\rangle $. A Noetherian local ring
$R$ (with maximal ideal $\mathfrak{m}$) is \textit{Cohen-Macaulay} if
depth$\left(  R\right)  =$ dim$\left(  R\right)  $, where the \textit{depth}
of $R$ is defined to be the maximum of the lengths of all regular sequences
whose members belong to $\mathfrak{m}$. A Cohen-Macaulay local ring $R$ is
\textit{Gorenstein }if Ext$_{R}^{\text{dim}\left(  R\right)  }\left(
R/\mathfrak{m},R\right)  \cong R/\mathfrak{m}$. A Noetherian local ring $R$ is said to
be a \textit{complete intersection} if there exists a regular local ring
$R^{\prime}$, such that $R\cong R^{\prime}/\left\langle f_{1},\ldots
,f_{q}\right\rangle $ for a finite set of elements $\left\{  f_{1}%
,\ldots,f_{q}\right\}  \subset R^{\prime}$ whose cardinality equals $q=$
dim$\left(  R^{\prime}\right)  -$ dim$\left(  R\right)  $. The hierarchy by
inclusion of the above types of Noetherian local rings is known to be
described by the following diagram:
\begin{equation}%
\begin{array}
[c]{ccc}%
\left\{  \text{Noetherian local rings}\right\}  & \supset & \left\{
\text{normal local rings}\right\} \\
\cup &  & \cup\\
\left\{  \text{Cohen-Macaulay local rings}\right\}  &  & \left\{
\text{regular local rings}\right\} \\
\cup &  & \cap\\
\left\{  \text{Gorenstein local rings}\right\}  & \supset & \left\{
\text{complete intersections (``c.i.'s'')}\right\}
\end{array}
\label{HIERARCHY}%
\end{equation}
\textsf{(b)} An arbitrary Noetherian ring $R$ and its associated affine scheme
Spec$\left(  R\right)  $ are called Cohen-Macaulay, Gorenstein, normal or
regular, respectively, iff \ all the localizations $R_{\mathfrak{m}}$ with respect
to all the members $\mathfrak{m}\in$ Max-Spec$\left(  R\right)  $ of the maximal
spectrum of $R$ are of this type. In particular, if the $R_{\mathfrak{m}} $'s for
all maximal ideals $\mathfrak{m}$ of $R$ are c.i.'s, then one often says that $R$
is a \textit{local complete intersection }(``l.c.i.'') to distinguish it from
the ``global'' ones. (A \textit{global complete intersection} (``g.c.i.'') is
defined to be a ring $R$ of finite type over a field $\mathbf{k}$ (i.e., an
affine $\mathbf{k}$-algebra), such that $R\cong\mathbf{k}\left[
\mathsf{T}_{1}..,\mathsf{T}_{d}\right]  \,/\,\left\langle \varphi_{1}\left(
\mathsf{T}_{1},..,\mathsf{T}_{d}\right)  ,..,\varphi_{q}\left(  \mathsf{T}%
_{1},..,\mathsf{T}_{d}\right)  \right\rangle $ \ for $q$ polynomials
$\varphi_{1},\ldots,\varphi_{q}$ from $\mathbf{k}\left[  \mathsf{T}%
_{1},..,\mathsf{T}_{d}\right]  $ with $q=d-$ dim$\left(  R\right)  $, cf.
\cite{Ishida, Nakajima}). Hence, the above inclusion hierarchy can be
generalized for all Noetherian rings, just by omitting in (\ref{HIERARCHY})
the word ``local'' and by substituting l.c.i.'s for c.i.'s.\medskip
\newline \textsf{(c) }Throughout the paper we consider only \textit{complex
varieties }$\left(  X,\mathcal{O}_{X}\right)  $, i.e., integral separated
schemes of finite type over $\mathbf{k}=\mathbb{C}$; thus, the punctual
algebraic behaviour of $X$ is determined by the stalks $\mathcal{O}_{X,x}$ of
its structure sheaf $\mathcal{O}_{X}$, and $X$ itself is said to have a given
algebraic property (as in \textsf{(b)}) whenever all $\mathcal{O}_{X,x}$'s
have the analogous property from (\ref{HIERARCHY}) for all $x\in X$.
Furthermore, via the \textsc{gaga}-correspondence (\cite{Serre},
\cite[\S2]{Gro2}) which preserves the above quoted algebraic
properties, we shall always 
work within the \textit{analytic category} by using the so-called
antiequivalence principle \cite{Gro1}, i.e., the usual contravariant functor
$\left(  X,x\right)  \leadsto\mathcal{O}_{X,x}^{\text{hol}}$ between the
category of isomorphy classes of germs of $X$ and the corresponding category
of isomorphy classes of analytic local rings at the marked points
$x$).\medskip\ \newline \textsf{(d) }For\textsf{\ }a complex variety $X$, we
denote by Sing$\left(  X\right)  =\left\{  x\in X\ \left|  \ \right.
\mathcal{O}_{X,x}^{\text{hol}}\text{ is a non-regular local ring}\right\}  $
its singular locus. By a \textit{desingularization} (or \textit{resolution of
singularities}) $f:\widehat{X}\rightarrow X$ of a non-smooth $X$, we mean a
``full'' or ``overall'' desingularization (if not mentioned), i.e.,
Sing$(\widehat{X})=\varnothing$. When we deal with \textit{partial}
desingularizations, we mention it explicitly. A partial desingularization
$f:X^{\prime}\rightarrow X$ of a normal, Gorenstein complex variety $X$ is
called\ \textit{non-discrepant} or simply \textit{crepant}, if the (up to
rational equivalence uniquely determined) difference $K_{X^{\prime}}-f^{\ast
}\left(  K_{X}\right)  $ vanishes. ($K_{X}$ and $K_{X^{\prime}}$ denote here
canonical divisors of $X$ and $X^{\prime}$, respectively). Furthermore,
$f:X^{\prime}\rightarrow X$ is \textit{projective} if $X^{\prime}$ admits an
$f$-ample Cartier divisor.

\section{Some basic facts from toric geometry\label{TORIC-G}}

\noindent In this section we introduce the brief toric glossary \textsf{(a)}%
-\textsf{(k) }and the notation which will be used in the subsequent sections.
For further details the reader is referred to the textbooks of Oda \cite{Oda},
Fulton \cite{Fulton} and Ewald \cite{Ewald}, and to the lecture notes
\cite{KKMS}. \medskip

\noindent\textsf{(a)} The \textit{linear hull, }the\textit{\ affine hull}, the
\textit{positive hull} and \textit{the convex hull} of a set $B$ of vectors of
$\mathbb{R}^{r}$, $r\geq1,$ will be denoted by lin$\left(  B\right)  $,
aff$\left(  B\right)  $, pos$\left(  B\right)  $ (or $\mathbb{R}_{\geq0}\,B$)
and conv$\left(  B\right)  $, respectively. The \textit{dimension} dim$\left(
B\right)  $ of a $B\subset\mathbb{R}^{r}$ is defined to be the dimension of
its affine hull. \newline \newline \textsf{(b) }Let $N$ be a free $\mathbb{Z}%
$-module of rank $r\geq1$. $N$ can be regarded as a \textit{lattice }in
$N_{\mathbb{R}}:=N\otimes_{\mathbb{Z}}\mathbb{R}\cong\mathbb{R}^{r}$.
An $n\in N$ is called \textit{primitive} if conv$\left(  \left\{
    \mathbf{0},n\right\} \right)  \cap N$ contains no other points
except $\mathbf{0}$ and $n$.\smallskip

Let $N$ be as above, $M:=$ Hom$_{\mathbb{Z}}\left(  N,\mathbb{Z}\right)  $ its
dual lattice, $N_{\mathbb{R}},M_{\mathbb{R}}$ their real scalar extensions,
and $\left\langle .,.\right\rangle :M_{\mathbb{R}}\times N_{\mathbb{R}%
}\rightarrow\mathbb{R}$ the natural $\mathbb{R}$-bilinear pairing. A subset
$\sigma$ of $N_{\mathbb{R}}$ is called \textit{convex polyhedral cone}
(\textit{c.p.c.}, for short) if there exist $n_{1},\ldots,n_{k}\in
N_{\mathbb{R}}$, such that $\sigma=$ pos$\left(  \left\{  n_{1},\ldots
,n_{k}\right\}  \right)  $. Its \textit{relative interior }int$\left(
\sigma\right)  $ is the usual topological interior of it, considered as subset
of lin$\left(  \sigma\right)  =\sigma+\left(  -\sigma\right)  $. The
\textit{dual cone} $\sigma^{\vee}$ of a c.p.c. $\sigma$ is a c.p. cone defined
by
\[
\sigma^{\vee}:=\left\{  \mathbf{y}\in M_{\mathbb{R}}\ \left|  \ \left\langle
\mathbf{y},\mathbf{x}\right\rangle \geq0,\ \forall\mathbf{x},\ \mathbf{x}%
\in\sigma\right.  \right\}  \;.\;
\]
Note that $\left(  \sigma^{\vee}\right)  ^{\vee}=\sigma$ and dim$\left(
\sigma\cap\left(  -\sigma\right)  \right)  +$ dim$\left(  \sigma^{\vee
}\right)  =$ dim$\left(  \sigma^{\vee}\cap\left(  -\sigma^{\vee}\right)
\right)  +$ dim$\left(  \sigma\right)  =r.$ A subset $\tau$ of a c.p.c.
$\sigma$ is called a \textit{face} of $\sigma$ (notation: $\tau\prec\sigma$),
if $\tau=\left\{  \mathbf{x}\in\sigma\ \left|  \ \left\langle m_{0}%
,\mathbf{x}\right\rangle =0\right.  \right\}  $, for some $m_{0}\in
\sigma^{\vee}$. A c.p.c. $\sigma=$ pos$\left(  \left\{  n_{1},\ldots
,n_{k}\right\}  \right)  $ is called \textit{simplicial} (resp.
\textit{rational}) if $n_{1},\ldots,n_{k}$ are $\mathbb{R}$-linearly
independent (resp. if $n_{1},\ldots,n_{k}\in N_{\mathbb{Q}}$, where
$N_{\mathbb{Q}}:=N\otimes_{\mathbb{Z}}\mathbb{Q}$). A \textit{strongly convex
polyhedral cone }(\textit{s.c.p.c.}, for short) is a c.p.c. $\sigma$ for which
$\sigma\cap\left(  -\sigma\right)  =\left\{  \mathbf{0}\right\}  $, i.e., for
which dim$\left(  \sigma^{\vee}\right)  =r$. The s.c.p. cones are
alternatively called \textit{pointed cones} (having $\mathbf{0}$ as their
apex).\newline \newline \textsf{(c) }If $\sigma\subset N_{\mathbb{R}}$ is a
rational c.p. cone, then the subsemigroup $\sigma\cap N$ of $N$ is a monoid.
The following proposition is due to Gordan, Hilbert and van der Corput and
describes its fundamental properties.

\begin{proposition}
[Minimal generating system]\label{MINGS}$\sigma\cap N$ is finitely generated
as additive semigroup. Moreover, if $\sigma$ is strongly convex, then among
all the systems of generators of $\sigma\cap N$, there is a system
$\mathbf{Hilb}_{N}\left(  \sigma\right)  $ of \emph{minimal cardinality}, which
is uniquely determined \emph{(}up to the ordering of its elements\emph{)} by
the following characterization\emph{:\smallskip}
\begin{equation}
\mathbf{Hilb}_{N}\left(  \sigma\right)  =\left\{  n\in\sigma\cap\left(
N\smallsetminus\left\{  \mathbf{0}\right\}  \right)  \ \left|  \
\begin{array}
[c]{l}%
n\ \text{\emph{cannot be expressed as the sum of two }}\\
\text{\emph{other vectors belonging to\ } }\sigma\cap\left(  N\smallsetminus
\left\{  \mathbf{0}\right\}  \right)
\end{array}
\right.  \right\}  \label{Hilbbasis}%
\end{equation}
$\mathbf{Hilb}_{N}\left(  \sigma\right)  $ \emph{is called }\textit{the
Hilbert basis of }$\sigma$ w.r.t. $N.$
\end{proposition}

\noindent

\noindent\textsf{(d)} For a lattice $N$ of rank $r$ having $M$ as its dual, we
define an $r$-dimensional \textit{algebraic torus }$T_{N}\cong\left(
\mathbb{C}^{\ast}\right)  ^{r}$ by setting $T_{N}:=
\thinspace$Hom$_{\mathbb{Z}}\left(
M,\mathbb{C}^{\ast}\right)  =N\otimes_{\mathbb{Z}}\mathbb{C}^{\ast}$. Every
$m\in M$ assigns a character $\mathbf{e}\left(  m\right)  :T_{N}%
\rightarrow\mathbb{C}^{\ast}$. Moreover, each $n\in N$ determines an
$1$-parameter subgroup
\[
\vartheta_{n}:\mathbb{C}^{\ast}\rightarrow T_{N}\ \ \ \text{with\ \ \ }%
\vartheta_{n}\left(  \lambda\right)  \left(  m\right)  :=\lambda^{\left\langle
m,n\right\rangle }\text{, \ \ for\ \ \ }\lambda\in\mathbb{C}^{\ast},\ m\in
M\ .\
\]
We can therefore identify $M$ with the character group of $T_{N}$ and $N$ with
the group of $1$-parameter subgroups of $T_{N}$. On the other hand, for a
rational s.c.p.c. $\sigma$ with $M\cap\sigma^{\vee}=\mathbb{Z}_{\geq0}%
\ m_{1}+\cdots+\mathbb{Z}_{\geq0}\ m_{\nu}$, we associate to the finitely
generated monoidal subalgebra $\mathbb{C}\left[  M\cap\sigma^{\vee}\right]
=\oplus_{m\in M\cap\sigma^{\vee}}\mathbf{e}\left(  m\right)  $ of the
$\mathbb{C}$-algebra $\mathbb{C}\left[  M\right]  =\oplus_{m\in M}%
\mathbf{e}\left(  m\right)  $ an affine complex variety
\[
U_{\sigma}:=\text{Max-Spec}\left(  \mathbb{C}\left[  M\cap\sigma^{\vee
}\right]  \right)  ,
\]
which can be identified with the set of semigroup homomorphisms :
\[
U_{\sigma}=\left\{  u:M\cap\sigma^{\vee}\ \rightarrow\mathbb{C\ }\left|
\begin{array}
[c]{c}%
\ u\left(  \mathbf{0}\right)  =1,\ u\left(  m+m^{\prime}\right)  =u\left(
m\right)  \cdot u\left(  m^{\prime}\right)  ,\smallskip\ \\
\text{for all \ \ }m,m^{\prime}\in M\cap\sigma^{\vee}%
\end{array}
\right.  \right\}  \ ,
\]
where $\mathbf{e}\left(  m\right)  \left(  u\right)  :=u\left(  m\right)
,\ \forall m,\ m\in M\cap\sigma^{\vee}\ $ and\ $\forall u,\ u\in U_{\sigma}$.

\begin{proposition}
[Embedding by binomials]\label{EMB}In the analytic category, $U_{\sigma}$,
identified with its image under the injective map $\left(  \mathbf{e}\left(
m_{1}\right)  ,\ldots,\mathbf{e}\left(  m_{\nu}\right)  \right)  :U_{\sigma
}\hookrightarrow\mathbb{C}^{\nu}$, can be regarded as an analytic set
determined by a system of equations of the form\emph{:} \emph{(monomial) =
(monomial).} This analytic structure induced on $U_{\sigma}$ is independent of
the semigroup generators $\left\{  m_{1},\ldots,m_{\nu}\right\}  $ and each
map $\mathbf{e}\left(  m\right)  $ on $U_{\sigma}$ is holomorphic w.r.t. it.
In particular, for $\tau\prec\sigma$, $U_{\tau}$ is an open subset of
$U_{\sigma}$. 
Moreover, if $\sigma$ is $r$-dimensional and 
$\#\left(  \mathbf{Hilb}_{M}\left(  \sigma^{\vee
}\right)  \right)  =k\ \left(  \leq\nu\right)  $, then $k$ is
nothing but the {\em embedding dimension} 
\textit{of} $U_{\sigma}$\textit{, i.e.
the }minimal \textit{number of generators of the maximal ideal of the local}
$\mathbb{C}$\textit{-algebra} $\mathcal{O}_{U_{\sigma},\ \mathbf{0}%
}^{\text{hol}}$.
\end{proposition}

\noindent\textit{Proof. }See Oda \cite[Prop. 1.2 and 1.3., pp. 4-7]{Oda}.
$_{\Box}\bigskip$

\noindent\textsf{(e)} A \textit{fan }w.r.t. a free $\mathbb{Z}$%
-module\textit{\ }$N$ is a finite collection $\Delta$ of rational s.c.p. cones
in $N_{\mathbb{R}}$, such that :\smallskip\newline (i) any face $\tau$ of
$\sigma\in\Delta$ belongs to $\Delta$, and\smallskip\newline (ii) for
$\sigma_{1},\sigma_{2}\in\Delta$, the intersection $\sigma_{1}\cap\sigma_{2}$
is a face of both $\sigma_{1}$ and $\sigma_{2}.\smallskip$\newline By $\left|
\Delta\right|  :=\cup\left\{  \sigma\ \left|  \ \sigma\in\Delta\right.
\right\}  $ one denotes the \textit{support} and by $\Delta\left(  i\right)  $
the set of all $i$-dimensional cones of a fan $\Delta$ for $0\leq i\leq r$. If
$\varrho\in\Delta\left(  1\right)  $ is a ray, then there exists a unique
primitive vector $n\left(  \varrho\right)  \in N\cap\varrho$ with
$\varrho=\mathbb{R}_{\geq0}\ n\left(  \varrho\right)  $ and each cone
$\sigma\in\Delta$ can be therefore written as
\[
\sigma=\sum_{\varrho\in\Delta\left(  1\right)  ,\ \varrho\prec\sigma
}\ \mathbb{R}_{\geq0}\ n\left(  \varrho\right)  \ \ .
\]
The set Gen$\left(  \sigma\right)  :=\left\{  n\left(  \varrho\right)
\ \left|  \ \varrho\in\Delta\left(  1\right)  ,\varrho\prec\sigma\right.
\right\}  $ is called the\textit{\ set of minimal generators }(within the pure
first skeleton) of $\sigma$. For $\Delta$ itself one defines analogously
Gen$\left(  \Delta\right)  :=\bigcup_{\sigma\in\Delta}$ Gen$\left(
\sigma\right)  .\bigskip$\newline \textsf{(f) }The \textit{toric variety
X}$\left(  N,\Delta\right)  $ associated to a fan\textit{\ }$\Delta$ w.r.t.
the lattice\textit{\ }$N$ is by definition the identification space
\begin{equation}
X\left(  N,\Delta\right)  :=((%
{\textstyle\coprod\limits_{\sigma\in\Delta}}
\ U_{\sigma})\ /\ \sim) \label{TORVAR}%
\end{equation}
with $U_{\sigma_{1}}\ni u_{1}\sim u_{2}\in U_{\sigma_{2}}$ if and only if
there is a $\tau\in\Delta,$ such that $\tau\prec\sigma_{1}\cap\sigma_{2}$ and
$u_{1}=u_{2}$ within $U_{\tau}$. $X\left(  N,\Delta\right)  $ is called
\textit{simplicial} if all the cones of $\Delta$ are simplicial. $X\left(
N,\Delta\right)  $ is compact iff $\left|  \Delta\right|  =N_{\mathbb{R}}$
\cite[Thm. 1.11, p. 16]{Oda}. Moreover, $X\left(  N,\Delta\right)  $ admits
a canonical $T_{N}$-action which extends the group multiplication of
$T_{N}=U_{\left\{  \mathbf{0}\right\}  }$:
\begin{equation}
T_{N}\times X\left(  N,\Delta\right)  \ni\left(  t,u\right)  \longmapsto
t\cdot u\in X\left(  N,\Delta\right)  \label{torus action}%
\end{equation}
where, for $u\in U_{\sigma}\subset X\left(  N,\Delta\right)  $, $\left(
t\cdot u\right)  \left(  m\right)  :=t\left(  m\right)  \cdot u\left(
m\right)  ,\ \forall m,\ m\in M\cap\sigma^{\vee}$ . The orbits w.r.t. the
action (\ref{torus action}) are parametrized by the set of all the cones
belonging to $\Delta$. For a $\tau\in\Delta$, we denote by orb$\left(
\tau\right)  $ (resp. by $V\left(  \tau\right)  $) the orbit (resp. the
closure of the orbit) which is associated to $\tau$. If $\tau\in\Delta$, then
$V\left(  \tau\right)  :=V\left(  \tau;\Delta\right)  :=X\left(  N\left(
\tau\right)  ,\text{ Star}\left(  \tau;\Delta\right)  \right)  $ is itself a
toric variety w.r.t.
\[
N\left(  \tau\right)  :=N\;/\;N_{\tau}\ ,\ \ \ \ \ \text{Star}\left(
\tau;\Delta\right)  :=\left\{  \overline{\sigma}\ \left|  \ \sigma\in
\Delta,\ \tau\prec\sigma\right.  \right\}  \ ,
\]
where $N_{\tau}$ is the sublattice $N\cap$ lin$\left(
  \mathbb{\tau}\right)$ of $N$ and $\overline{\sigma}\ =\left(
\sigma+\left(  N_{\tau}\right)  _{\mathbb{R}}\right)  /\left(  N_{\tau
}\right)  _{\mathbb{R}}$ denotes the image of $\sigma$ in $N\left(
\tau\right)  _{\mathbb{R}}=N_{\mathbb{R}}/\left(  N_{\tau}\right)
_{\mathbb{R}}$.\medskip\medskip

\noindent\textsf{(g)}\textit{\ }The behaviour of toric varieties with regard
to the algebraic properties (\ref{HIERARCHY}) has as follows.\textit{\ }

\begin{theorem}
[Normality and CM-property]All toric varieties are normal and Cohen-Macaulay.
\end{theorem}

\noindent\textit{Proof. }For a proof of the normality property see
\cite[Thm. 1.4, p. 7]{Oda}. The CM-property for toric varieties was
first shown by 
Hochster in \cite{Hochster}. See also Kempf \cite[Thm. 14, p. 52]{KKMS}, and
Oda \cite[3.9, p. 125]{Oda}. $_{\Box}\medskip$\newline In fact, by the
definition (\ref{TORVAR}) of $X\left(  N,\Delta\right)  $, all the algebraic
properties of this kind are \textit{local }with respect to its affine
covering, i.e., it is enough to be checked for the affine toric varieties
$U_{\sigma}$ for all (maximal) cones $\sigma$ of the fan $\Delta$.

\begin{definition}
[Multiplicities and basic cones]\emph{Let }$N$\emph{\ be a free }%
$Z$\emph{-module of rank }$r$\emph{\ and} $\sigma\subset N_{\mathbb{R}}$
\emph{a simplicial,} \emph{rational s.c.p.c. of dimension} $d\leq r$. $\sigma$
\emph{can be obviously written as} $\sigma=\varrho_{1}+\cdots+\varrho_{d}%
$\emph{, for distinct rays} $\varrho_{1},\ldots,\varrho_{d}$. The
\textit{multiplicity} \emph{mult}$\left(  \sigma;N\right)  $ \emph{of}
$\sigma$ \emph{with respect to} $N$ \emph{is defined as the index}
\[
\text{\emph{mult}}\left(  \sigma;N\right)  := 
| \medspace N_\sigma \medspace : \medspace \mathbb{Z} n(\varrho_1) +
\cdots + \mathbb{Z} n(\varrho_d) \medspace |
\]
\emph{If mult}$\left(  \sigma;N\right)  =1$\emph{, then }$\sigma$\emph{\ is
called a} \textit{basic cone} \emph{w.r.t.} $N$.
\end{definition}

\begin{theorem}
[Smoothness criterion]\label{SMCR}The affine toric variety $U_{\sigma}$ is
smooth iff $\sigma$ is basic \textit{w.r.t.} $N$. \emph{(}Correspondingly, an
arbitrary toric variety $X\left(  N,\Delta\right)  $ is smooth if and only if
it is simplicial and each s.c.p. cone $\sigma\in\Delta$ is basic
\textit{w.r.t.} $N$.\emph{)}
\end{theorem}

\noindent\textit{Proof. }See \cite[ch. I, Thm. 4, p. 14]{KKMS}, and
\cite[Thm. 1.10, p. 15]{Oda}. $_{\Box}\medskip$\newline Next Theorem is due to
Stanley \cite[\S6]{Stanley1}, who worked directly with the monoidal
$\mathbb{C}$-algebra $\mathbb{C}\left[  M\cap\sigma^{\vee}\right]  $, as well
as to Ishida \cite[\S7]{Ishida}, Danilov and Reid \cite[p. 294]{Reid1},
who provided a purely algebraic-geometric characterization of the
Gorensteinness property.

\begin{theorem}
[Gorenstein property]\label{GOR-PR}The following conditions are
equivalent\emph{:\medskip}\newline \emph{(i)} $\ \ U_{\sigma}$ is
Gorenstein.\medskip\ \newline \emph{(ii)} \ There exists an element
$m_{\sigma}$ of $M$, such that $M\cap\left(  \text{\emph{int}}\left(
\sigma^{\vee}\right)  \right)  =m_{\sigma}+M\cap\sigma^{\vee}$.\medskip
\newline $\smallskip$\emph{(iii)} \emph{Gen}$\left(  \sigma\right)
\subset\mathbf{H} $, where $\mathbf{H}$ denotes an affine hyperplane of
\ $\left(  N_{\sigma}\right)  _{\mathbb{R}}$ that contains a lattice
basis of $N_\sigma$.\medskip\newline Moreover, if
\emph{dim}$\left(  \sigma\right)  =r$, then $m_{\sigma}$ in \emph{(ii) }is a
uniquely determined primitive element of $M\cap\left(  \text{\emph{int}%
}\left(  \sigma^{\vee}\right)  \right)  $ and $\mathbf{H}$ in \emph{(iii)
}equals $\mathbf{H}=\left\{  \mathbf{x}\in N_{\mathbb{R}}\ \left|
\ \left\langle m_{\sigma},\mathbf{x}\right\rangle =1\right.  \right\}  $.
\end{theorem}

\noindent A geometric interpretation of the remaining ``finer'' algebraic
property, namely whether $U_{\sigma}$ is a l.c.i. or not, in terms of the
defining fan, is due to Nakajima and will be presented separately in
\S\ref{NAK}, Thm. \ref{Nak-thm}.\bigskip

\noindent\textsf{(h)}\textit{\ }A \textit{map of fans\ }$\varpi:\left(
N^{\prime},\Delta^{\prime}\right)  \rightarrow\left(  N,\Delta\right)  $ is a
$\mathbb{Z}$-linear homomorphism $\varpi:N^{\prime}\rightarrow N$ whose scalar
extension $\varpi\otimes_{\mathbb{Z}}$id$_{\mathbb{R}}:N_{\mathbb{R}}^{\prime
}\rightarrow N_{\mathbb{R}}$ satisfies the property:
\[
\forall\sigma^{\prime},\ \sigma^{\prime}\in\Delta^{\prime}\ \ \text{ }%
\exists\ \sigma,\ \sigma\in\Delta\ \ \text{ with\ \
  }\varpi\otimes_{\mathbb{Z}}\text{id}_{\mathbb{R}} \left(
\sigma^{\prime}\right)  \subset\sigma\,.
\]
$\varpi\otimes_{\mathbb{Z}}$id$_{\mathbb{C}^{\ast}}:T_{N^{\prime}}=N^{\prime
}\otimes_{\mathbb{Z}}\mathbb{C}^{\ast}\rightarrow T_{N}=N\otimes_{\mathbb{Z}%
}\mathbb{C}^{\ast}$ is a homomorphism from $T_{N^{\prime}}$ to $T_{N}$ and the
scalar extension $\varpi^{\vee}\otimes_{\mathbb{Z}}$id$_{\mathbb{R}%
}:M_{\mathbb{R}}\rightarrow M_{\mathbb{R}}^{\prime}$ of the dual $\mathbb{Z}%
$-linear map $\varpi^{\vee}:M\rightarrow M^{\prime}$ induces canonically an
\textit{equivariant holomorphic map }$\varpi_{\ast}:X\left(  N^{\prime}%
,\Delta^{\prime}\right)  \rightarrow X\left(  N,\Delta\right)  $. This map
is\textit{\ proper} if and only if $\varpi^{-1}\left(  \left|  \Delta\right|
\right)  =\left|  \Delta^{\prime}\right|  .$ In particular, if $N=N^{\prime}$
and $\Delta^{\prime}$ is a refinement of $\Delta$, then id$_{\ast}:X\left(
N,\Delta^{\prime}\right)  \rightarrow X\left(  N,\Delta\right)  $ is
\textit{proper}\emph{\ }and \textit{birational }cf. \cite[Thm. 1.15 and
Cor. 1.18]{Oda}.\bigskip

\noindent\textsf{(i)}\textit{\ }By Carath\'{e}odory's Theorem concerning
convex polyhedral cones (cf. \cite[III 2.6 \& V 4.2]{Ewald}) one can choose a
refinement $\Delta^{\prime}$ of any given fan $\Delta$, so that $\Delta
^{\prime}$ becomes simplicial. Since further subdivisions of $\Delta^{\prime}$
reduce the multiplicities of its cones, we may arrive (after finitely many
subdivisions) at a fan $\widetilde{\Delta}$ having only basic cones. Hence,
for every toric variety $X\left(  N,\Delta\right)  $ there exists a refinement
$\widetilde{\Delta}$ of $\Delta$ consisting of exclusively basic cones w.r.t.
$N$, i.e., such that $f=$ id$_{\ast}:X(N,\widetilde{\Delta})\longrightarrow
X\left(  N,\Delta\right)  $ is a $T_{N}$-equivariant (full)
desingularization.\bigskip\newline \textsf{(j) }The group of $T_{N}$-invariant
Weil divisors of a toric variety $X\left(  N,\Delta\right)  $ has the set
$\{V(\varrho)\,\left|  \,\varrho\in\Delta(1)\right.  \}$ as $\mathbb{Z}%
$-basis. In fact, such a divisor $D$ is of the form $D=D_{\psi}$, where
$D_{\psi}:=-\sum_{\varrho\in\Delta\left(  1\right)  }\psi(n(\varrho
))V(\varrho)$ and $\psi:\left|  \Delta\right|  \rightarrow\mathbb{R}$ a
\textit{PL-}$\Delta$\textit{-support function}, i.e., an $\mathbb{R}$-valued,
positively homogeneous function on $\left|  \Delta\right|  $ with $\psi
(N\cap\left|  \Delta\right|  )\subset\mathbb{Z}$ which is piecewise linear and
upper convex on each $\sigma\in\Delta$. (\textit{Upper convex} on a $\sigma
\in\Delta$ means that $\psi\left|  _{\sigma}(\mathbf{x}+\mathbf{x}^{\prime
})\right.  \geq\psi\left|  _{\sigma}(\mathbf{x})+\right.  \psi\left|
_{\sigma}(\mathbf{x}^{\prime})\right.  $, for all $\mathbf{x},\mathbf{x}%
^{\prime}\in\sigma$). For example, the canonical divisor $K_{X\left(
N,\Delta\right)  }$ of $X\left(  N,\Delta\right)  $ equals $D_{\psi}$ for
$\psi$ a PL-$\Delta$-support function\textit{\ }with $\psi(n(\varrho))=1$, for
all rays $\varrho\in\Delta\left(  1\right)  $. A divisor $D=D_{\psi}$ is
Cartier iff $\psi$ is a \textit{linear} $\Delta$-support function (i.e.,
$\psi\left|  _{\sigma}\right.  $ is overall linear on each $\sigma\in\Delta$).
Obviously, $D_{\psi}$ is $\mathbb{Q}$-Cartier iff $k\cdot\psi$ is a linear
$\Delta$-support function for some $k\in\mathbb{N}$.

\begin{theorem}
[Ampleness criterion]\label{AMPLE}A $T_{N}$-invariant \emph{(}$\mathbb{Q}%
$-\emph{) }Cartier divisor $D=D_{\psi}$ of a toric variety $X\left(
N,\Delta\right)  $ of dimension $r$ is ample if and only if there exists a
$\kappa\in\mathbb{N}$, such that $\kappa\cdot\psi$ is a \emph{strictly upper
convex} linear $\Delta$-support function, i.e., iff \ for every $\sigma
\in\Delta(r)$ there is a unique $m_{\sigma}\in M=$ \emph{Hom}$_{\mathbb{Z}%
}(N,\mathbb{Z})$, such that $\kappa\cdot\psi(\mathbf{x})\leq\left\langle
m_{\sigma},\mathbf{x}\right\rangle $, for all $\mathbf{x}\in\left|
\Delta\right|  $ , with equality being valid iff $\mathbf{x}\in\sigma$.
\end{theorem}

\noindent\textit{Proof. }It follows from \cite[Thm. 13, p. 48]{KKMS}. $_{\Box
}\bigskip$\textit{\newline }\textsf{(k) }Throughout the paper, by a
\textit{polytope} in an euclidean space, is meant the convex hull of finitely
many points or, equivalently, a bounded polyhedron. A\emph{\ }\textit{lattice
polytope}\emph{\ }$P$\emph{\ }embedded in a given euclidean space is a
polytope whose set vert$\left(  P\right)  $ of vertices belongs to a reference
lattice within this space. If $M$ is a free $\mathbb{Z}$-module of rank $r$,
$N=$ Hom$_{\mathbb{Z}}(M,\mathbb{Z})$ its dual, and $P\subset M_{\mathbb{R}%
}\cong\mathbb{R}^{r}$ an $r$-dimensional lattice polytope w.r.t. $M$, then
there is a unique fan $\Delta^{\left(  P\right)  }$ in $N_{\mathbb{R}}$, the
so-called \textit{normal fan }of $P$, so that the corresponding $r$%
-dimensional toric variety $X(N,\Delta^{\left(  P\right)  })$ is projective
and endowed with a distinguished $T_{N}$-invariant ample Cartier divisor
$D_{P}:=D_{\psi}$ which is induced by the strictly upper convex support
function $\psi:N_{\mathbb{R}}\rightarrow\mathbb{R}$, with $\psi\left(
\mathbf{x}\right)  :=$ min$\{\left\langle \mathbf{y},\mathbf{x}\right\rangle
\,\left|  \,\mathbf{y}\in P\right.  \}$; and conversely, regarding a
projective toric variety $X\left(  N,\Delta\right)  $ and a $T_{N}$-invariant
ample Cartier divisor $D=D_{\psi}$ on it as our starting-point data, we win a
characteristic $r$-dimensional lattice polytope $P=P_{D}$ assigned to $D$,
with $P_{D}=\{\mathbf{y}\in M_{\mathbb{R}}\,\left|  \,\right.  \left\langle
\mathbf{y},\mathbf{x}\right\rangle \geq\psi\left(  \mathbf{x}\right)
,\forall\mathbf{x},\,\mathbf{x}\in N_{\mathbb{R}}\}$ (cf. Oda
\cite[\S2.4]{Oda}).

\section{Torus-equivariant crepant projective resolutions \newline of
Gorenstein toric singularities via b.c.-triangulations\label{BC-TR}}

\noindent We shall henceforth focus our attention to Gorenstein toric
singularities and to their desired resolutions.\bigskip

\noindent\textsf{(a) }Let $N$ be a free $\mathbb{Z}$-module of rank
$r\geq2$ and 
$\sigma\subset N_{\mathbb{R}}$ a rational s.c.p.c. of dimension $d\leq r $. We
identify $U_{\sigma}$ with $X\left(  N,\Delta\right)  $, where $\Delta$
denotes the fan consisting of $\sigma$ together with all of its faces. Since
$N\left(  \sigma\right)  =N/N_{\sigma}$ is torsion free, there exists a
lattice decomposition $N=N_{\sigma}\oplus\breve{N}$, inducing a decomposition
of its dual $M=M_{\sigma}\oplus\breve{M}$, where $M_{\sigma}=$
Hom$_{\mathbb{Z}}\left(  N_{\sigma},\mathbb{Z}\right)  $ and $\breve{M}=$
Hom$_{\mathbb{Z}}(\breve{N},\mathbb{Z)}$. Writing $\sigma$ as $\sigma
=\sigma^{\prime}\oplus\left\{  \mathbf{0}\right\}  $ with $\sigma^{\prime} $ a
$d$-dimensional cone in $\left(  N_{\sigma}\right)  _{\mathbb{R}}$, we obtain
decompositions
\[
T_{N}\ \cong T_{N_{\sigma}}\ \times T_{\breve{N}}\ \ \ \ \ \text{and\ \ \ \ }%
M\cap\sigma^{\vee}=\left(  M\cap\left(  \sigma^{\prime}\right)  ^{\vee
}\right)  \oplus\breve{M}\ ,
\]
which give rise to the analytic isomorphisms:
\[%
\begin{array}
[b]{ccc}%
U_{\sigma} & \cong\ U_{\sigma^{\prime}}\times T_{\breve{N}}\ \cong
\ U_{\sigma^{\prime}}\times T_{N\left(  \sigma\right)  }\ \cong &
U_{\sigma^{\prime}}\times\left(  \mathbb{C}^{\ast}\right)  ^{r-d}\smallskip\\
\parallel &  & \parallel\\
X\left(  N,\Delta\right)  &  & X\left(  N_{\sigma},\Delta^{\prime}\right)
\times\left(  \mathbb{C}^{\ast}\right)  ^{r-d}%
\end{array}
\]
with $\Delta^{\prime}$ the fan consisting of $\sigma^{\prime}$ together with
all of its faces (cf. \cite[p. 29]{Fulton}, and \cite[Thm.VI.2.12, p.
223]{Ewald}). $U_{\sigma}$ can be therefore viewed as as a fiber bundle over
$U_{\sigma^{\prime}}$ having an $\left(  r-d\right)  $-dimensional algebraic
torus as its typical fibre. Obviously, the study of the algebraic properties
(mentioned in \S\ref{INTRO}) for $U_{\sigma}$ can be reduced to that of the
corresponding properties of $\ U_{\sigma^{\prime}}$. (For instance, the
singular locus of $U_{\sigma}$ equals Sing$\left(  U_{\sigma}\right)  =$
Sing$\left(  U_{\sigma^{\prime}}\right)  \times\left(  \mathbb{C}^{\ast
}\right)  ^{r-d}$). In fact, the main reason for preferring to work with
$U_{\sigma^{\prime}}$ (or with the germ $\left(  U_{\sigma^{\prime}%
},\text{orb}\left(  \sigma^{\prime}\right)  \right)  )$ instead of $U_{\sigma
}$, is that since lin$\left(  \sigma^{\prime}\right)  =\left(  N_{\sigma
}\right)  _{\mathbb{R}}$, the orbit orb$\left(  \sigma^{\prime}\right)  \in
U_{\sigma^{\prime}}$ is the unique fixed closed point under the action of
$T_{N_{\sigma}}$ on $U_{\sigma^{\prime}}$.

\begin{definition}
[Singular representatives]\label{SINR}\emph{If }$\sigma$\emph{\ is non-basic
w.r.t. }$N$\emph{, then }$U_{\sigma^{\prime}}$ \emph{will be called
}\textit{the singular representative}\emph{\ of }$U_{\sigma}$ \emph{and
orb}$\left(  \sigma^{\prime}\right)  \in U_{\sigma^{\prime}}$ \emph{the
associated }\textit{distinguished}\emph{\ singular point within the singular
locus Sing}$\left(  U_{\sigma^{\prime}}\right)  $ \emph{of }$U_{\sigma
^{\prime}}=X\left(  N_{\sigma},\Delta^{\prime}\right)  $.\emph{\ }
\end{definition}

\begin{definition}
[Splitting codimension]\label{SPL}\emph{If }$\sigma$\emph{\ is non-basic
w.r.t. }$N$\emph{, then} \emph{it is also useful to introduce the notion of
the ``}\textit{splitting codimension}\emph{'' of} \emph{orb}$\left(
\sigma^{\prime}\right)  \in U_{\sigma^{\prime}}$ \emph{as the number}
\[
\text{\emph{min\ }}\left\{  \varkappa\in\left\{  2,\ldots,d\right\}  \ \left|
\
\begin{array}
[c]{c}%
U_{\sigma^{\prime}}\cong U_{\sigma^{\prime\prime}}\times\mathbb{C}%
^{d-\varkappa}\emph{,}\text{\ }\emph{\ }\text{\emph{for\ some\ } }%
\sigma^{\prime\prime}\prec\sigma^{\prime}\text{ }\\
\text{\emph{with \ dim}}\left(  \sigma^{\prime\prime}\right)  =\varkappa\text{
\ \ \emph{and \ Sing}}\left(  U_{\sigma^{\prime\prime}}\right)  \neq
\varnothing
\end{array}
\right.  \right\}  \ .
\]
\emph{(In \cite[p. 231]{DHH} and \cite[p. 202]{DHZ1} there is a misprint 
  in this definition: one must replace therein max by min.)}
\emph{If this number equals} $d$\emph{, then }$\left(  U_{\sigma^{\prime}%
},\emph{orb}\left(  \sigma^{\prime}\right)  \right)  $\emph{\ will be called
an} \textit{msc-singularity}\emph{, i.e., a singularity having the maximum}
\emph{splitting codimension.\smallskip}
\end{definition}

\noindent\textsf{(b) }Gorenstein toric affine varieties are completely
determined by suitable lattice polytopes.

\begin{definition}
[Lattice equivalence]\label{LP}\emph{\ If }$\ N_{1}$\emph{\ and }$N_{2}%
$\emph{\ are two free }$\mathbb{Z}$\emph{-modules (not necessarily of the same
rank) and }$P_{1}\subset\left(  N_{1}\right)  _{\mathbb{R}}$\emph{, }%
$P_{2}\subset\left(  N_{2}\right)  _{\mathbb{R}}$\emph{\ two lattice polytopes
w.r.t. them, we shall say that }$P_{1}$\emph{\ and }$P_{2}$\emph{\ are}
\textit{lattice equivalent }\emph{to each other, and denote this by }%
$P_{1}\sim P_{2}$\emph{, if }$P_{1}$\emph{\ is affinely equivalent to }$P_{2}%
$\emph{\ via an affine map }$\varpi:\left(  N_{1}\right)  _{\mathbb{R}%
}\rightarrow\left(  N_{2}\right)  _{\mathbb{R}}$\emph{, such that the}
\emph{restiction} $\varpi\left|  _{\text{\emph{aff}}\left(  P\right)
}\right.  :$ \emph{aff}$\left(  P\right)  \rightarrow$ \emph{aff}$\left(
P^{\prime}\newline \right)  $ \emph{is a bijection} \emph{mapping }$P_{1}%
$\emph{\ onto the (necessarily equidimensional) polytope }$P_{2}$
\emph{,} \emph{and, in
addition, }$N_{P_{1}}$\emph{\ is mapped bijectively onto the lattice
}$N_{P_{2}}$\emph{, where
}$N_{P_{j}}$\emph{\ is the affine sublattice aff}$\left(  P_{j}\right)
\cap N_{j}$\emph{\ of }$N_j$\emph{, }$j=1,2.$\emph{\ If
}$N_{1}=N_{2}=:N$\emph{\ and rk}$\left(  N\right)  =$\emph{\ dim}$\left(
P_{1}\right)  =$ \emph{dim}$\left(  P_{2}\right)  $\emph{, then these }%
$\varpi$\emph{'s are exactly the} \textit{affine integral transformations}
\emph{which are composed of unimodular }$N$-\emph{transformations and
  }$N$-\emph{translations.}
\end{definition}

\noindent Let now $U_{\sigma}=X\left(  N,\Delta\right)  $ be a $d$-dimensional
affine toric variety as in \textsf{(a)} and $U_{\sigma^{\prime}}=X\left(
N_{\sigma},\Delta^{\prime}\right)  $. Assuming that $U_{\sigma} $ is
Gorenstein, we may pass to another analytically isomorphic ``standard''
representative as follows: Denote by $\mathbb{Z}^{d}$ the standard rectangular
lattice in $\mathbb{R}^{d}$ and by $(\mathbb{Z}^{d})^{\vee}$ its dual lattice
within $(\mathbb{R}^{d})^{\vee}=$ Hom$_{\mathbb{R}}(\mathbb{R}^{d}%
,\mathbb{R})$. Since dim$\left(  \sigma^{\prime}\right)  =$ rk$\left(
N_{\sigma}\right)  =d$, or equivalently, since $\left(  \sigma^{\prime
}\right)  ^{\vee}$ is strongly convex in $\left(  M_{\sigma}\right)
_{\mathbb{R}}$, Thm. \ref{GOR-PR} (iii) implies
\[
\text{Gen}\left(  \sigma^{\prime}\right)  \subset\mathbf{H}^{\left(  d\right)
}\text{ \ \ with\ \ \ \ }\mathbf{H}^{\left(  d\right)  }:=\left\{
\mathbf{x}\in\left(  N_{\sigma}\right)  _{\mathbb{R}}\ \left|  \ \left\langle
m_{\sigma^{\prime}},\mathbf{x}\right\rangle =1\right.  \right\}  ,
\]
for a unique primitive $m_{\sigma^{\prime}}\in M_{\sigma}.$ Clearly,
$\sigma^{\prime}\cap\mathbf{H}^{\left(  d\right)  }$ is a $\left(  d-1\right)
$-dimensional lattice polytope (w.r.t. $N_{\sigma}$). We choose \ a specific
$\mathbb{Z}$-module isomorphism $\Upsilon:N_{\sigma}\overset{\cong
}{\longrightarrow}\mathbb{Z}^{d}$ inducing an $\mathbb{R}$-vector space
isomorphism $\Phi=\Upsilon\otimes_{\mathbb{Z}}$id$_{\mathbb{R}}:\left(
N_{\sigma}\right)  _{\mathbb{R}}\overset{\cong}{\longrightarrow}\mathbb{R}%
^{d}$, such that
\[
\Phi\left(  m_{\sigma^{\prime}}\right)  =(1,\underset{\left(  d-1\right)
\text{-times}}{\underbrace{0,0,...,0,0}})\Longrightarrow\Phi\left(
\mathbf{H}^{\left(  d\right)  }\right)  =\left\{  \mathbf{x}=\left(
x_{1},\ldots,x_{d}\right)  \in\mathbb{R}^{d}\ \left|  \ x_{1}=1\right.
\right\}  \mathbf{=:\bar{H}}^{\left(  d\right)  }\text{.}%
\]
Obviously, $P:=\Phi\left(  \sigma^{\prime}\cap\mathbf{H}^{\left(  d\right)
}\right)  \subset\mathbf{\bar{H}}^{\left(  d\right)  }$ is a lattice $\left(
d-1\right)  $-polytope (w.r.t. $\mathbb{Z}^{d}$). Defining
\[
\tau_{P}:=\text{pos}(P)=\left\{  \kappa\ \mathbf{x}\in\mathbb{R}%
^{d}\mathbf{\ }\left|  \ \kappa\in\mathbb{R}_{\geq0},\ \mathbf{x}\in P\right.
\right\}  ,\ \ \ \Delta_{P}:=\left\{  \text{\emph{\ }}\tau_{P}\text{ together
with all of its faces}\right\}  ,
\]
(cf. Figure \textbf{1}) we obtain easily the following Lemma:

\begin{lemma}
\label{ISOM}\emph{(i) }There exists a torus-equivariant analytic isomorphism
\[
U_{\sigma^{\prime}}=X\left(  N_{\sigma},\Delta^{\prime}\right)  \cong
U_{\tau_{P}}=X(\mathbb{Z}^{d},\Delta_{P})\ \ \ (=\text{\emph{Max-Spec}%
}(\mathbb{C}\left[  (\mathbb{Z}^{d})^{\vee}\cap\tau_{P}^{\vee}\right]  ))
\]
mapping \emph{orb}$\left(  \sigma^{\prime}\right)  $ onto \emph{orb}$\left(
\tau_{P}\right)  .\medskip$\newline \emph{(ii)} If \ $Q\subset\mathbf{\bar{H}%
}^{\left(  d\right)  }$ is a lattice $\left(  d-1\right)  $-polytope
\emph{(}w.r.t. $\mathbb{Z}^{d}$\emph{)}, then $P\sim Q$ $\ $iff there exists a
torus-equivariant analytic isomorphism $U_{\tau_{P}}\cong U_{\tau_{Q}}$
mapping $\emph{orb}\left(  \tau_{P}\right)  $ onto $\emph{orb}\left(  \tau
_{Q}\right)  $.
\end{lemma}

\begin{center}
\begin{picture}(0,0)%
\epsfig{file=tauP.pstex}%
\end{picture}%
\setlength{\unitlength}{987sp}%
\begingroup\makeatletter\ifx\SetFigFont\undefined%
\gdef\SetFigFont#1#2#3#4#5{%
  \reset@font\fontsize{#1}{#2pt}%
  \fontfamily{#3}\fontseries{#4}\fontshape{#5}%
  \selectfont}%
\fi\endgroup%
\begin{picture}(15624,9924)(1189,-11773)
\put(7951,-3811){\makebox(0,0)[lb]{\smash{\SetFigFont{11}{13.2}{\rmdefault}{\mddefault}{\updefault}$\mathbf{\bar{H}^{(d)}}$}}}
\put(6751,-6136){\makebox(0,0)[lb]{\smash{\SetFigFont{11}{13.2}{\rmdefault}{\mddefault}{\updefault}$P$}}}
\put(8626,-8836){\makebox(0,0)[lb]{\smash{\SetFigFont{11}{13.2}{\rmdefault}{\mddefault}{\updefault}$\tau_P$}}}
\end{picture}

\\[.5\baselineskip]
\textbf{Figure 1}
\end{center}

\begin{definition}
[Standard representatives]\label{STAN}\emph{\ Any member of the isomorphy
class of the underlying space} $U_{\tau_{P}}=X(\mathbb{Z}^{d},\Delta_{P}%
)$\emph{\ of the distinguished Gorenstein point orb}$\left(  \tau_{P}\right)
$\emph{\ (as in \ref{ISOM}(ii)) is said to be} \textit{a standard
representative }\emph{of }$U_{\sigma}$ \emph{associated to the lattice
polytope }$P$\emph{, and, in particular, }\textit{a singular standard
representative }\emph{of }$U_{\sigma}$\emph{, whenever }$\sigma$\emph{\ is
non-basic w.r.t. }$N$. \emph{(In this case, the splitting codimension of }
\emph{orb}$\left(  \tau_{P}\right)  $\emph{\ is defined to be the splitting
codimension of orb}$\left(  \sigma^{\prime}\right)  $\emph{.)}
\end{definition}

\noindent\textsf{(c) }Suppose that $\sigma$ is a non-basic c.p. cone w.r.t.
$N$. From the above discussion it is now clear that for desingularizing
$U_{\sigma}$, it suffices to resolve a singular representative $U_{\sigma
^{\prime}}$, and for $U_{\sigma}$ Gorenstein, a standard singular
representative $U_{\tau_{P}}$ of it. In the latter case, for any
torus-equivariant partial desingularization \ $f=$ id$_{\ast}:X(\mathbb{Z}%
^{d},\widehat{\Delta}_{P})\longrightarrow X(\mathbb{Z}^{d},\Delta_{P}%
)=U_{\tau_{P}}$ coming from a refinement $\widehat{\Delta}_{P}$ of $\Delta
_{P}$ (cf. \S\ref{TORIC-G}, \textsf{(h)}-\textsf{(i)}) there are one-to-one
correspondences:
\begin{equation}%
\begin{array}
[c]{ccc}%
\varrho & \in & \widehat{\Delta}_{P}\left(  1\right)  \smallsetminus\Delta
_{P}\left(  1\right) \\
\updownarrow &  & \updownarrow\\
n\left(  \varrho\right)  & \in & \text{Gen}(\widehat{\Delta}_{P}%
)\,\smallsetminus\text{Gen}\left(  \tau_{P}\right) \\
\updownarrow &  & \updownarrow\\
D_{n\left(  \varrho\right)  }:=V\left(  \varrho\right)  =V(\varrho
;\widehat{\Delta}_{P}) & \in & \left\{
\begin{array}
[c]{c}%
\text{exceptional prime divisors}\\
\text{with respect to }f
\end{array}
\right\}
\end{array}
\label{exc-div}%
\end{equation}
Moreover, as we shall see below in proposition \ref{CREPANT}, it is possible
to describe certain intrinsic algebraic-geometric properties of those $f$'s
which are crepant or\thinspace/\thinspace and projective exclusively in terms
of lattice triangulations of the polytope $P$ defining $U_{\tau_{P}}$. For
this reason, before proceeding to this description, we recall some central
notions from the theory of polytopal subdivisions which will be crucially
utilized in the rest of the paper.

\begin{definition}
[Polytopal subdivisions and refinements]\label{PSUB}\emph{(i)} \emph{A
}\textit{polytopal complex }\emph{is a finite family} $\mathcal{S}$ \ \emph{of
polytopes in an euclidean space }$\mathbb{R}^{\ell}$%
\emph{, so that the intersection of any two of its
polytopes constitutes always a common face of each of them. The}
\textit{dimension} \emph{dim}$\left(  \mathcal{S}\right)  $ \emph{of such an
}$\mathcal{S}$ \emph{is defined to be the largest possible dimension of a
polytope belonging to it. }$\mathcal{S}$ \emph{is called a} \textit{pure}
\emph{polytopal complex if every polytope in $\mathcal{S}$ is contained in one
of dimension dim}$\left(  \mathcal{S}\right)  $.\medskip\newline \emph{(ii)
Let }$\mathcal{V}$ \emph{denote a finite set of points in an euclidean space,
such that }$P=$ \emph{conv}$\left(  \mathcal{V}\right)  $ \emph{is a }%
$k$\emph{-dimensional polytope. A }\textit{polytopal subdivision}
$\mathcal{S}$ \emph{of }$P$ \emph{is a finite family} $\mathcal{S}=\left\{
P_{1},P_{2},\ldots,P_{\nu}\right\}  $ \emph{of} $k$\emph{-dimensional
polytopes, such that:\smallskip}\newline \textsc{a.}\emph{\ $\mathcal{S}$ is a
pure }$k$\emph{-dimensional polytopal complex.\smallskip\newline }\textsc{b.}
\emph{The space supporting }$P$\emph{\ is the union of spaces supporting
}$P_{1},P_{2},\ldots,P_{\nu}.\smallskip$\emph{\newline }\textsc{c.}
\emph{vert}$\left(  P_{i}\right)  \subseteq\mathcal{V}$ \emph{, for all}
$i\in\left\{  1,2,\ldots,\nu\right\}  .\medskip$\newline \emph{(iii) A
polytopal subdivision} $\mathcal{S}$ \emph{of }$P$ \emph{as in (ii) is called
a }\textit{triangulation}\emph{\ of }$P$\emph{\ if each }$P_{i}$\emph{,}
$1\leq i\leq\nu$\emph{, is a }$k$\emph{-dimensional simplex.\medskip}%
\newline \emph{(iv) Suppose that} $\mathcal{S}=\left\{  P_{1},P_{2}%
,\ldots,P_{\nu}\right\}  $\emph{,} $\mathcal{S}^{\prime}=\left\{
P_{1}^{\prime},P_{2}^{\prime},\ldots,P_{\mu}^{\prime}\right\}  $ \emph{are two
polytopal subdivisions of} $P.$ \emph{Then} $\mathcal{S}^{\prime}$\emph{\ is
a} \textit{refinement}\emph{\ of }$\mathcal{S}$ \emph{if for each} $j$, $1\leq
j\leq\mu$\emph{, there exists an} $i$, $1\leq i\leq\nu$\emph{, such}
\emph{that }$P_{j}^{\prime}\subseteq P_{i}$.
\end{definition}

\begin{definition}
[Coherent subdivisions]\label{COHSUB}\emph{A polytopal subdivision
$\mathcal{S}$ of} \emph{a polytope }$P\subset\mathbb{R}^{k}$\emph{\ is called}
\textit{coherent} \emph{(or, alternatively, }\textit{regular}\emph{, cf.
\cite[ 5.3]{Ziegler}) if \ }$P$ \emph{is the image }$\pi\left(  Q\right)  =P$
\emph{of a polytope} $Q$\emph{\ }$\subset\mathbb{R}^{k+1}$ \emph{under the
projection map}
\begin{equation}
\mathbb{R}^{k+1}\ni\left(  x_{1},\ldots,x_{k},x_{k+1}\right)  \overset{\pi
}{\longmapsto}\left(  x_{1},\ldots,x_{k}\right)  \in\mathbb{R}^{k}%
\ \label{projection}%
\end{equation}
\emph{so that }$\mathcal{S}=\left\{  \pi\left(  F\right)  :F\ \text{\emph{is a
lower face of }}Q\right\}  $\emph{, where }\textit{the lower faces of }$Q$
\emph{are the faces for which some outward normal vector has negative }$\left(
k+1\right)  $\emph{-st coordinate. (The set of all lower faces of }$Q$
\emph{is sometimes called the }\textit{lower envelope} \emph{of }$Q$\emph{).}
\end{definition}

\noindent The next two Lemmas describe further useful conditions which are
equivalent to the coherency of \emph{$\mathcal{S}$.}

\begin{lemma}
[Coherency and strictly upper convex functions]\label{COHUCF}A polytopal
subdivision\emph{\ $\mathcal{S}$} $\ $of $P$ is coherent iff there exists
\emph{a strictly upper convex} $\mathcal{S}$\emph{-support function}
$\psi:\left|  \mathcal{S}\right|  \rightarrow\mathbb{R}$, i.e. a
piecewise-linear real function defined on the underlying\emph{\ }%
space\emph{\ }$\left|  \mathcal{S}\right|  $, for which
\[
\psi(t\ \mathbf{x}+\left(  1-t\right)  \mathbf{\ y)}\geq t\ \psi\left(
\mathbf{x}\right)  +\left(  1-t\right)  \ \psi\left(  \mathbf{y}\right)
\emph{,}\text{\emph{\ for all} \ }\mathbf{x},\mathbf{y}\in\left|
\mathcal{S}\right|  ,\emph{\ }\text{\emph{and}\ }\emph{\ }t\in\left[
0,1\right]  \ ,
\]
so that its domains of linearity are exactly the polytopes of $\mathcal{S}$
\emph{\ }having maximal dimension.
\end{lemma}

\noindent\textit{Proof}. If $\mathcal{S}$ is coherent, then $\mathcal{S}%
=\left\{  \pi\left(  F\right)  :F\ \text{is a lower face of\emph{\ }%
}Q\right\}  $, with $\pi:\mathbb{R}^{k+1}\rightarrow\mathbb{R}^{k}\ $\ the
projection (\ref{projection}) and $Q$ a polytope in $\mathbb{R}^{k+1}$. The
function $\psi:\left|  \mathcal{S}\right|  \rightarrow\mathbb{R}$ defined by
setting
\[
\psi(\mathbf{x}):=\text{ max }\left\{  t\in\mathbb{R\ }\left|  \ \right.
(\mathbf{x},-t)\in Q\right\}  ,\text{ \ for all }\mathbf{x}=(x_{1}%
,\ldots,x_{k})\in\left|  \mathcal{S}\right|  =P,
\]
is strictly upper convex. (Some authors prefer to work with convex support
functions instead of upper convex ones and use min and $(\mathbf{x},t)$
instead of max and $(\mathbf{x},-t)$. But this is just a sign
convention).\medskip\newline Conversely, if $\psi:\left|  \mathcal{S}\right|
\rightarrow\mathbb{R}$ is assumed to be a strictly upper convex support
function, then $\mathcal{S}$ is coherent in the sense of \ref{COHSUB} by
defining $Q$ to be the polytope conv$\left\{  (\mathbf{x,}-\psi(\mathbf{x}%
))\in\mathbb{R}^{k+1}\mathbb{\ }\left|  \ \right.  \mathbf{x}\in P\right\}  $.
$_{\Box}$

\begin{lemma}
[Coherency and ``heights'']\label{HEIGHTS}Let $\mathcal{V}$ be a finite set of
points in $\mathbb{R}^{k} $ and $P=$ \emph{conv}$\left(  \mathcal{V}\right)  $
\emph{.} A \emph{height function} on $\mathcal{V}$ \ is defined to be a
function $\omega:\mathcal{V}\rightarrow\mathbb{R}$. \emph{(}The values
$\omega\left(  \mathbf{v}\right)  $, $\mathbf{v}\in\mathcal{V}$, are called
\emph{``heights'').} Every height function $\omega$ on $\mathcal{V}$ induces a
coherent polytopal subdivision $\mathcal{S}_{\omega}$ of the polytope $P=$
\emph{conv}$\left(  \mathcal{V}\right)  $ with \emph{vert}$\left(
\mathcal{S}_{\omega}\right)  \subseteq\mathcal{V}$\emph{; } and
conversely, each 
coherent polytopal subdivision $\mathcal{S}$ of $P=$ \emph{conv}$\left(
\mathcal{V}\right)  $ \emph{\ }with \emph{vert}$\left(  \mathcal{S}\right)
\subseteq\mathcal{V}$\emph{\ }is of the form
$\mathcal{S}=\mathcal{S}_{\omega}$, for some height function
$\omega$. \emph{\ }
\end{lemma}

\noindent\textit{Proof}. Let $\omega$ be a height function on $\mathcal{V}$.
The heights can be used to ``lift'' the point configuration $\mathcal{V}$ into
the next dimension and to define $Q_{\omega}:=$ conv$(\{(\mathbf{v}%
,\omega(\mathbf{v}))\in\mathbb{R}^{k+1}\ \left|  \ \mathbf{v}\in
\mathcal{V}\right.  \})$. The lower envelope of the polytope $Q_{\omega}$ is a
pure polytopal complex having dimension equal to dim$(P)$. Its image under the
projection (\ref{projection}) determines a (necessarily coherent) polytopal
subdivision $\mathcal{S}_{\omega}$ of $P$ with vert$\left(  \mathcal{S}%
_{\omega}\right)  \subseteq\mathcal{V}$. In fact, if $\{\mathbf{v}_{i_{1}%
},..,\mathbf{v}_{i_{\mu}}\}$ are the vertices a polytope belonging to
$\mathcal{S}_{\omega}$, then $\{(\mathbf{v}_{i_{1}},\omega(\mathbf{v}_{i_{1}%
})),..,(\mathbf{v}_{i_{\mu}},\omega(\mathbf{v}_{i_{\mu}}))\}$ is the vertex
set of a face of the lower envelope of $Q_{\omega}$.\medskip\newline Let now
$\mathcal{S}$ denote an arbitrary coherent polytopal subdivision $\mathcal{S}$
of $P$ \emph{\ }with vert$\left(  \mathcal{S}\right)
\subseteq\mathcal{V} $\emph{\ }.
By Lemma \ref{COHUCF} there exists a strictly upper convex support function
$\psi:\left|  \mathcal{S}\right|  \rightarrow\mathbb{R}$. Using the height
function $\omega:=(-\psi)\left|  _{\mathcal{V}}\right.  $ we obtain
$\mathcal{S}=\mathcal{S}_{\omega}$. $_{\Box}$

\begin{remark}
\emph{For ``generic'' choices of }$\omega$\emph{'s the coherent polytopal
subdivisions} $\mathcal{S}_{\omega}$ \emph{are }\textit{triangulations}%
\emph{\ of} $P$ \emph{(cf. \cite[p. 215 and p. 228]{GKZ}, and
  \cite[p. 64]{Sturmfels}).}
\end{remark}

\begin{definition}
[Lattice subdivisions]\emph{A} \textit{lattice subdivision} $\mathcal{S}$
$\ $\emph{of a lattice polytope} $P$ \emph{is a} \emph{polytopal subdivision
of }$P$\emph{, such that the set vert}$\left(  \mathcal{S}\right)  $\emph{\ of
the vertices of }$\mathcal{S}$\emph{\ \ belongs to the reference lattice (and
vert}$\left(  P\right)  \subseteq$\emph{\ vert}$\left(  \mathcal{S}\right)
$\emph{). A }\textit{lattice triangulation} \emph{of a lattice polytope} $P$
\emph{is a} \emph{lattice subdivision of }$P$\emph{\ which, in addition, is a
triangulation\ (in the sense of \ref{PSUB}).}
\end{definition}

\begin{definition}
[Maximal and basic triangulations]\emph{(i) A lattice polytope }$P$\emph{\ is
called} elementary \emph{if the lattice points belonging to it are exactly its
vertices. A lattice simplex is said to be} \textit{basic} \emph{or
}\textit{unimodular }\emph{if its vertices constitute a part of an affine }%
$\mathbb{Z}$\emph{-basis of the reference lattice (or equivalently, if its
relative, normalized volume} \emph{equals} $1$\emph{).\medskip}\newline
\emph{(ii) A lattice triangulation }$\mathcal{T}$ $\ $\emph{of a lattice
polytope} $P$ \emph{is defined to be} \textit{maximal} \emph{(resp.}
\textit{basic}\emph{), if it consists only of elementary (resp. basic) simplices.}
\end{definition}

\begin{definition}
[``b.c.''-triangulations]\emph{A }\textit{b.c.-triangulation}\emph{\ will be
used as abbreviation for a basic, coherent triangulation of a lattice polytope.}
\end{definition}

\noindent Reverting to Gorenstein affine toric varieties, we explain how
torus-equivariant crepant or\thinspace/\thinspace and projective
desingularizations can be constructed by means of lattice triangulations.

\begin{proposition}
[Crepant desingularizations and triangulations]\label{CREPANT}Every
torus-equivariant partial crepant desingularization of a standard
representative $U_{\tau_{P}}$ of a Gorenstein affine toric variety
\ $U_{\sigma}$ \emph{(}as in \emph{\ref{STAN}}, with $P\subset\mathbf{\bar{H}%
}^{\left(  d\right)  }$ a lattice polytope w.r.t. $\mathbb{Z}^{d}$\emph{)},
induced by a subdivision of $\Delta_{P}$ into simplicial s.c.p. cones, is of
the form
\begin{equation}
f=f_{\mathcal{T}}:X(\mathbb{Z}^{d},\widehat{\Delta}_{P}\left(  \mathcal{T}%
\right)  )\longrightarrow X(\mathbb{Z}^{d},\Delta_{P})=U_{\tau_{P}%
}\ \label{TAU}%
\end{equation}
where $\widehat{\Delta}_{P}=$ $\widehat{\Delta}_{P}\left(  \mathcal{T}\right)
:=\left\{  \sigma_{\mathbf{s}},\ \!\mathbf{s}\in\mathcal{T}\right\}  $ is
determined by a lattice triangulation $\mathcal{T}$ of $P$ with
\[
\sigma_{\mathbf{s}}:=\left\{  \kappa\ \mathbf{x}\in\mathbb{R}^{d}%
\mathbf{\ }\left|  \ \kappa\in\mathbb{R}_{\geq0},\ \mathbf{x}\in
\mathbf{s}\right.  \right\}  .
\]
By \emph{(\ref{exc-div}) }the set of exceptional prime divisors equals
$\left\{  D_{n}=V\left(  \mathbb{R}_{\geq0}\,n\right)  \ \left|  \ n\in\left(
P\mathbb{r}\text{\emph{vert}}\left(  P\right)  \right)  \cap\right.
\mathbb{Z}^{d}\right\}  $. Moreover, such an $f_{\mathcal{T}}$ has the
following properties\emph{:\smallskip}\newline \emph{(i)\ \ \ }$f_{\mathcal{T}%
}$ \ \emph{is maximal w.r.t. discrepancy }$\Longleftrightarrow\mathcal{T}$ is
maximal.\medskip\newline \emph{(ii)\ }$f_{\mathcal{T}}$ \ \emph{is full
(}i.e., $X(\mathbb{Z}^{d},\widehat{\Delta}_{P}\left(  \mathcal{T}\right)  )$
is overall smooth\emph{) }$\Longleftrightarrow\mathcal{T}$ is basic.\medskip
\newline \emph{(iii) }$f_{\mathcal{T}}$ \ \emph{is projective (}%
i.e.,\emph{\ }$X(\mathbb{Z}^{d},\widehat{\Delta}_{P}\left(  \mathcal{T}%
\right)  )$ is quasiprojective\emph{) }$\Longleftrightarrow\mathcal{T}$ is
coherent.\smallskip\newline 
\end{proposition}

\noindent\textit{Proof. }$\ $ Let \ $f:X(\mathbb{Z}^{d},\widehat{\Delta}%
_{P})\longrightarrow X(\mathbb{Z}^{d},\Delta_{P})=U_{\tau_{P}}$ denote an
arbitrary torus-equivariant partial desingularization of $U_{\tau_{P}}$
induced by a subdivision $\widehat{\Delta}_{P}$ of $\Delta_{P}$ into
simplicial s.c.p. cones. The discrepancy of $f$ equals
\[
K_{X\left(  \mathbb{Z}^{d},\text{ }\widehat{\Delta}_{P}\right)  }-f^{\ast
}\left(  K_{U_{\tau_{P}}}\right)  =\left[  -\sum_{\varrho^{\prime}\in\left(
\widehat{\Delta}_{P}\left(  1\right)  \mathbb{r\Delta}_{P}\left(  1\right)
\right)  }\ \widehat{D}_{n\left(  \varrho^{\prime}\right)  }-\sum_{\varrho
\in\Delta_{P}\left(  1\right)  }\ \widehat{D}_{n\left(  \varrho\right)
}\right]  -f^{\ast}\left(  -\sum_{\varrho\in\Delta_{P}\left(  1\right)
}\ D_{n\left(  \varrho\right)  }\right)
\]
where $D_{n\left(  \varrho\right)  }:=V\left(  \varrho,\Delta_{P}\right)  $,
$\widehat{D}_{n\left(  \varrho\right)  }:=V(\varrho,\widehat{\Delta}_{P})$,
for all rays $\varrho$ of $\Delta_{P}$ and $\widehat{\Delta}_{P}$,
respectively, and
\[
f^{\ast}\left(  -\sum_{\varrho\in\Delta_{P}\left(  1\right)  }\ D_{n\left(
\varrho\right)  }\right)  =-\sum_{\varrho\in\Delta_{P}\left(  1\right)
}\ \widehat{D}_{n\left(  \varrho\right)  }-\sum_{\varrho^{\prime}\in\left(
\widehat{\Delta}_{P}\left(  1\right)  \mathbb{r\Delta}_{P}\left(  1\right)
\right)  }\mu_{\varrho^{\prime}}\ \widehat{D}_{n\left(  \varrho^{\prime
}\right)  }%
\]
with $\mu_{\varrho^{\prime}}$'s $\in\mathbb{Q}_{\geq0}$. If $\phi
_{T_{\mathbb{Z}^{d}}}$ is the rational differential form generating the
dualizing sheaf of the torus $T_{\mathbb{Z}^{d}}$, then the dualizing sheaf of
$U_{\tau_{P}}$ is isomorphic to $\mathbb{C}\left[  (\mathbb{Z}^{d})^{\vee}%
\cap\text{ int}(\tau^{\vee})\right]  \cdot\phi_{T_{\mathbb{Z}^{d}}}$. Since
$U_{\tau_{P}}$ is Gorenstein, $\mathbb{C}\left[  (\mathbb{Z}^{d})^{\vee}%
\cap\text{ int}(\tau^{\vee})\right]  \cdot\phi_{T_{\mathbb{Z}^{d}}}$ is
generated by $\mathbf{e}((1,0,...,0,0))\cdot\phi_{T_{\mathbb{Z}^{d}}}$ (cf.
Thm. \ref{GOR-PR} and subsection \textsf{(b)}), and $K_{U_{\tau_{P}}}$ is
trivial. The preservation of Gorensteinness for $X(\mathbb{Z}^{d}%
,\widehat{\Delta}_{P})$ is equivalent to say that, for each member of its
affine cover $\{U_{\widehat{\sigma}}\,\left|  \,\widehat{\sigma}\in
\widehat{\Delta}_{P}\right.  (d)\}$, the sheaf of sections of the canonical
divisor $K_{X\left(  \mathbb{Z}^{d},\text{ }\widehat{\Delta}_{P}\right)  }$
over $U_{\widehat{\sigma}}$ is isomorphic to $\mathbb{C}\left[  (\mathbb{Z}%
^{d})^{\vee}\cap\text{ int}((\widehat{\sigma})^{\vee})\right]  \cdot
\phi_{T_{\mathbb{Z}^{d}}}$ and is therefore generated by $\mathbf{e}%
((1,0,...,0,0))\cdot\phi_{T_{\mathbb{Z}^{d}}}$. The order of vanishing for the
divisor div$(\mathbf{e}((1,0,...,0,0))\cdot\phi_{T_{\mathbb{Z}^{d}}})$ which
is associated to this common single generator along the $\widehat{D}_{n\left(
\varrho^{\prime}\right)  }$'s, $\varrho^{\prime}\in(\widehat{\Delta}%
_{P}\left(  1\right)  \mathbb{r\Delta}_{P}\left(  1\right)  )$, equals
\[
\mu_{\varrho^{\prime}}=\text{ord}_{\widehat{D}_{n\left(  \varrho^{\prime
}\right)  }}\left(  \text{div}(\mathbf{e}((1,0,...,0,0))\cdot\phi
_{T_{\mathbb{Z}^{d}}})\right)  =\,\langle(1,\underset{\left(  d-1\right)
\text{-times}}{\underbrace{0,0,...,0,0}}),n\left(
\varrho^{\prime}\right) \rangle
\]
(cf. Fulton \cite[Lemma of p. 61]{Fulton}). From the above equations we
deduce
\[
K_{X\left(  \mathbb{Z}^{d},\text{ }\widehat{\Delta}_{P}\right)  }-f^{\ast
}\left(  K_{U_{\tau_{P}}}\right)  =\sum_{\varrho^{\prime}\in\left(
\widehat{\Delta}_{P}\left(  1\right)  \mathbb{r\Delta}_{P}\left(  1\right)
\right)  }\ (\langle(1,\underset{\left(  d-1\right)  \text{-times}}{\underbrace
{0,0,...,0,0}}),n\left(  \varrho^{\prime}\right)  \rangle-1)\ \widehat{D}_{n\left(
\varrho^{\prime}\right)  }\ .
\]
Thus, $f$ is crepant iff
\begin{equation}
\text{Gen}(\widehat{\Delta}_{P})\subset\mathbf{\bar{H}}^{\left(  d\right)  }
\label{GEN-HYP}%
\end{equation}
i.e., iff $f$ is of the form (\ref{TAU}). Now property (i) is obvious. For
(ii) observe that (\ref{GEN-HYP}) implies for all $\ \!\mathbf{s}%
\in\mathcal{T}$:$\;\sigma_{\mathbf{s}}$ is a basic cone $\Longleftrightarrow
\mathbf{s}$ is a basic simplex. Concerning (iii), note that all
torus-invariant Weil divisors of $X(\mathbb{Z}^{d},\widehat{\Delta}_{P}\left(
\mathcal{T}\right)  )$ are $\mathbb{Q}$-Cartier because this toric variety is
$\mathbb{Q}$-factorial. Clearly, for every strictly upper convex linear
$\widehat{\Delta}_{P}\left(  \mathcal{T}\right)  $-support function $\psi$ (in
the sense of \S\ref{TORIC-G} \textsf{(j)}), the restriction $\psi\left|
_{\mathcal{T}}\right.  $ is a strictly upper convex $\mathcal{T}$-support
function (as in \ref{COHUCF}); and conversely, as it was explained in
\cite[\S4]{DHZ1}, to any $\mathcal{T}$-support function $\psi$, one may
canonically assign (eventually after suitable perturbation of the defining
inequalities and\thinspace/\thinspace or multiplication by a scalar) a
strictly upper convex linear $\widehat{\Delta}_{P}\left(  \mathcal{T}\right)
$-support function $\psi^{\prime}$ (with $\psi^{\prime}(\left|  \widehat
{\Delta}_{P}\left(  \mathcal{T}\right)  \right|  \cap\mathbb{Z}^{d}%
)\subset\mathbb{Q}$ or even in $\mathbb{Z}$). To finish the proof we apply
Lemma \ref{COHUCF} for $\psi^{\prime}$ and Theorem \ref{AMPLE} for the divisor
$D_{\psi^{\prime}}$. $_{\Box}$

\begin{remark}
\emph{The birational morphisms }$f_{\mathcal{T}}$\emph{, for} $\mathcal{T}%
$\emph{'s} \emph{maximal and coherent, can be decomposed into more elementary
toric contractions (see Reid \cite[($0.2$)-($0.3$)]{Reid2}). In
several cases, these contractions are directly expressible as chains of
blow-downs (cf. \cite[7.2 and \S9]{DH}).\smallskip}
\end{remark}

\noindent$\bullet$ Every lattice polytope $P$ can be clearly embedded, up to
an affine transformation, into $\mathbf{\bar{H}}^{\left(  d\right)  }$, with
$d=$ dim$\left(  P\right)  +1$, and its supporting cone $\tau_{P}%
\subset\mathbb{R}^{d}$ gives rise to the construction of an affine Gorenstein
variety $U_{\tau_{P}}$. Consequently, if we restrict the initial Question
\ref{QUESTION1} of the introduction to the category of Gorenstein toric
singularities (and their \textit{torus-equivariant} resolutions), our previous
discussion in subsection \textsf{(b)}, together with the ``bridge'' which is
built by proposition \ref{CREPANT} and connects algebraic with discrete
geometric statements, enable us to reformulate it as follows:

\begin{question}
\label{QUESTION2}Under which conditions does a given lattice polytope $P$ of
dimension $\geq3$ admit of b.c.-triangulations\thinspace\emph{?}
\end{question}

\begin{remark}
\label{Un-Sec}\emph{(i) All elementary triangles are basic, but already in
dimension }$3$\emph{\ there exist counterexamples of elementary simplices
which are non-basic. Moreover, already in dimension }$2$ \emph{(i.e., for
certain lattice polygons) there is a plethora of non-coherent (but necessarily
basic) maximal triangulations. Hence, the problem of the existence of
b.c.-triangulations turns out to be very subtle in general. The required extra
``conditions'' in the formulation of \ Question \ref{QUESTION2} depend
essentially on the representatives of the coordinates of vertices of the given
lattice polytope }$P$ \emph{within its lattice equivalence class.
Unfortunately, regarding these integer coordinates as freely moving
``parameters'', we see that in high dimensions} \emph{they are ``too many'' to
handle (even for simplices and even if we reduce them by suitable unimodular
transformations like Hermite normal form transformations). This is why a first
realistic attempt to answer \ref{QUESTION2} partially (or at least to find
}\textit{sufficient conditions} \emph{for the above existence problem) seems
to be feasible only by the consideration of some }\textit{special}%
\emph{\ }\textit{families}\emph{\ of }$P$\emph{'s. In the present paper we
deal with altogether three families of lattice polytopes and prove that they
admit b.c.-triangulations (see below \ref{FANO}, \ref{FANOS}, \ref{COMP},
\ref{ZONOS}, \ref{NAKPOL}, and \ref{MAIN2}). The third one is exactly that
corresponding to the toric l.c.i.-singularities and has some
interesting members in common (and also not in common) with the first
two (see \S\ref{EXTREME}).\medskip
\newline (ii) For any finite set of points $\mathcal{V}$ in an }%
$\mathbb{R}^{d}$\emph{,} \textit{all triangulations} $\mathcal{T}$ \emph{of
the polytope} $P=$ \emph{conv}$\left(  \mathcal{V}\right)  $ \emph{with}
\emph{vert}$\left(  \mathcal{T}\right)  \subseteq\mathcal{V}$ $\ $\emph{are
parametrized by the vertices of a ``gigantic'' polytope }$\mathbf{Un}\left(
\mathcal{V}\right)  $\emph{, the so-called} \textit{universal polytope}
\emph{of }$P$\emph{\ (see Billera, Filliman \& Sturmfels \cite[\S3]{BFS}, and
de Loera, Ho\c{s}ten, Santos \& Sturmfels \cite[\S1-\S4]{DL-H-S-S}).
}$\mathbf{Un}\left(  \mathcal{V}\right)  $ \emph{contains a subpolytope}
$\mathbf{Sec}\left(  \mathcal{V}\right)  $ \emph{whose vertices parametrize
only the coherent }$\mathcal{T}$\emph{'s}. $\mathbf{Sec}\left(  \mathcal{V}%
\right)  $ \emph{is} \emph{in most of the cases considerably ``big'' too, and
is called} \emph{the }\textit{secondary polytope}\emph{\ of }$P$\emph{. (For
the main concepts of the theory of secondary polytopes the reader is referred
to \cite{BFS}, Oda \& Park \cite{OP}, Ziegler \cite[Lecture $9$]{Ziegler}}%
\emph{, as well as to the treatment of Gelfand, Kapranov \& Zelevinsky
\cite[Ch. 7]{GKZ}. In practice, working with examples for which the
cardinality of the given }$\mathcal{V}$\emph{\ 's}$\ $\emph{is relatively
small, an enumeration of the vertices of }$\mathbf{Sec}\left(  \mathcal{V}%
\right)  $ \emph{can be easily achieved by making use of the \textsc{maple}%
}-\emph{package \textsc{puntos} \cite{DELOERA} of de Loera).\medskip}%
\newline \emph{(iii) In the particular case in which }$P\subset\mathbf{\bar
{H}}^{\left(  d\right)  }\subset\mathbb{R}^{d}$\emph{\ is a lattice polytope
(w.r.t. }$\mathbb{Z}^{d}$\emph{) and }$\mathcal{V}=P\cap\mathbb{Z}^{d}$\emph{,
the b.c.-triangulations of }$P$\emph{\ correspond to a }\textit{very
special}\emph{\ (not necessarily non-empty) ``mysterious'' subset
}$\mathbf{BC}(\mathcal{V})$ \emph{of vert}$\left(  \mathbf{Sec}\left(
\mathcal{V}\right)  \right)  $\emph{. Thus, since \ref{QUESTION2} asks for
conditions under which }$\mathbf{BC}(\mathcal{V})\neq\varnothing$\emph{, the
expected theoretical answer(s) would surely require a much more extensive
study for }$\mathbf{Sec}\left(  \mathcal{V}\right)  $ \emph{itself. At this
point, we should also stress that in high dimensions ``exotic pathological
counterexamples'' exist! For instance, Hibi and Ohsugi \cite{Hi-O}
discovered recently a }$9$\emph{-dimensional }$0/1$\emph{-polytope (with }$15
$ \emph{vertices) having basic triangulations, but none of whose coherent
triangulations is basic.\medskip}\newline \emph{(iv) Passing by a connected
vertex path from one vertex of }$\mathbf{Sec}\left(  \mathcal{V}\right)  $
\emph{to another, we perform a finite series of ``bistellar operations'' which
are nothing but ``flops'' in the algebraic-geometric terminology
\cite[\S3]{OP}.\smallskip}
\end{remark}

\noindent\textsf{(d)} \ We next present two characteristic families of lattice
polytopes which admit specific b.c.-triangulations $\mathcal{T}$ leading to
projective, crepant, full desingularization morphisms $f_{\mathcal{T}}$ with
explicitly describable exceptional prime divisors.

\begin{definition}
[Fano polytopes]\label{FANO}\emph{A lattice polytope }$Q$ \emph{is called a
}\textit{Fano polytope}\emph{\ if }$Q\sim P$\emph{, where }$P\subset
\mathbb{R}^{d}$ \emph{denotes a lattice polytope (w.r.t.} $\mathbb{Z}^{d}%
$\emph{)} \emph{containing exactly one lattice point in its relative
  interior, which, together with the vertices of each facet, forms an
  affine lattice basis of }$(\mathbb{Z}^d)_P$ .
\end{definition}

\begin{proposition}
\label{FANOS}Let $P\subset\mathbf{\bar{H}}^{\left(  d\right)  }\subset
\mathbb{R}^{d}$ be a Fano polytope \emph{(}w.r.t. $\mathbb{Z}^{d}$\emph{)
}with \emph{int}$\left(  P\right)  \cap\mathbb{Z}^{d}=\left\{  n_{0}\right\}
$.\smallskip\ \newline \emph{(i) }The canonical lattice triangulation
$\mathcal{T}^{\text{\emph{can}}}:=\left\{  \left\{  n_{0}\right\}  \star
F\ \left|  \ F\ \text{\emph{face of }}P\right.  \right\}  $ constructed by
``joins'' \emph{(}i.e., by considering the pyramids over the faces of $P$ with
$n_{0}$ as apex\emph{) }is a b.c.-triangulation of $P$.\smallskip
\newline \emph{(ii) }The induced torus-equivariant projective, crepant, full
desingularization
\[
f_{\mathcal{T}^{\text{\emph{can}}}}:X(\mathbb{Z}^{d},\widehat{\Delta}%
_{P}\left(  \mathcal{T}^{\text{\emph{can}}}\right)  )\longrightarrow
X(\mathbb{Z}^{d},\Delta_{P})=U_{\tau_{P}}\
\]
possesses exactly one exceptional prime divisor
\[
D_{n_{0}}=V\left(  \mathbb{R}_{\geq0}n_{0}\right)  =X((\mathbb{Z}%
^{d})({\mathbb{R}_{\geq0}n_{0}}),\emph{Star}(\mathbb{R}_{\geq0}n_{0}%
;\widehat{\Delta}_{P}\left(  \mathcal{T}^{\text{\emph{can}}}\right)  ))
\]
which is a projective, toric Fano manifold.
\end{proposition}

\noindent\textit{Proof. }(i) follows directly from
\cite[Thm. 3.5]{DHZ1} and Lemma \ref{LA} below. For
(ii) note that $P$ is, in particular, a \textit{reflexive} polytope (cf.
\cite[4.1.5]{Bat1}). This is equivalent to say that its polar polytope
$P^{\ast}\subset(\mathbb{R}^{d})^{\vee}$ with respect to aff$\left(  P\right)
$ (having $n_{0}$ as its ``origin'') is again a lattice polytope (w.r.t.
$(\mathbb{Z}^{d})^{\vee}$). Since the rays of the fan Star$(\mathbb{R}_{\geq
0}n_{0};\widehat{\Delta}_{P}\left(  \mathcal{T}^{\text{can}}\right))  $ are
exactly the $1$-dimensional cones determined by joining $n_{0}$ with the
vertices of $P$, $D_{n_{0}}$ is the $(d-1)$-dimensional projective toric
variety associated to the normal fan of $P^{\ast}$ (see \S\ref{TORIC-G}
\textsf{(k)}). Thus, the fan Star$(\mathbb{R}_{\geq0}n_{0};\widehat{\Delta
}_{P}\left(  \mathcal{T}^{\text{can}}\right) )$ is \textit{strongly polytopal}
(see \cite[V.4.3 and V.4.4, p. 159]{Ewald}), and is composed of exclusively
basic cones. Consequently, the exceptional prime divisor $D_{n_{0}}$ has to
carry the analytic structure of a smooth toric variety which is \textit{Fano},
i.e., whose antidualizing sheaf is ample (cf. \cite[2.1.6 \& 2.2.23]{Bat1}).
$_{\Box}$

\begin{definition}
[$\mathbb{H}_{d}$-compatible polytopes]\label{COMP}\emph{Let }$\mathbb{H}_{d}$
\emph{denote the affine hyperplane arrangement (of} \emph{type} $\widetilde
{\mathcal{A}}_{d}$\emph{) in} $\mathbb{R}^{d}$ \emph{consisting of the union
of hyperplanes}
\[
\left\{  \left\{  \mathbf{x\in\,}\mathbb{R}^{d}\,\left|  \,x_{i}%
=\kappa\right.  \right\}  ,1\leq i\leq d,\ \kappa\in\mathbb{Z}\right\}
\cup\left\{  \left\{  \mathbf{x\in\,}\mathbb{R}^{d}\,\left|  \,x_{i}%
-x_{j}=\lambda\right.  \right\}  ,1\leq i<j\leq d,\ \lambda\in\mathbb{Z}%
\right\}  \ .
\]
\emph{A lattice polytope }$Q$ \emph{will be called a }$\mathbb{H}_{d}%
$\textit{-compatible polytope}\emph{\ if }$Q\sim P$\emph{, where }%
$P\subset\mathbb{R}^{d}$ \emph{denotes a lattice polytope (w.r.t.}
$\mathbb{Z}^{d}$\emph{), such that the affine hulls aff}$\left(  F\right)  $
\emph{for all facets }$F$\emph{\ of }$P$\emph{\ belong to} $\mathbb{H}_{d}$.
\emph{The affine hyperplane arrangement }$\mathbb{H}_{d}$ \emph{induces a
basic triangulation }$\mathcal{T}_{\mathbb{H}_{d}}$\emph{\ of the
}whole\emph{\ space }$\mathbb{R}^{d}$. \emph{In fact, }$\mathcal{T}%
_{\mathbb{H}_{d}}$\emph{\ is also coherent because there exists an overall
well-defined strictly upper convex function on }$\left|  \mathcal{T}%
_{\mathbb{H}_{d}}\right|  $ \emph{being constructible by means of appropriate
sums of Heaviside functions (see \cite[Ch. 3]{KKMS}, and
\cite[Prop. 6.1]{DHZ1}).}
\end{definition}

\begin{theorem}
\label{ZONOS}For $d\geq2$, let 
$P\subset\mathbf{\bar{H}}^{\left(  d\right)}\subset\mathbb{R}^{d}$
be a $\left(  d-1\right)  $-dimensional $\mathbb{H}_{d}$-compatible
polytope w.r.t.\ $\mathbb{Z}^{d}$, and let 
$\mathcal{T}_{\mathbb{H}_{d}}\left|  _{P}\right.  $
denote the triangulation
of $P$ determined by the restriction of $\mathcal{T}_{\mathbb{H}_{d}}$ to
$\left|  P\right|  $. Then $\mathcal{T}_{\mathbb{H}_{d}}\left|  _{P}\right.  $
is a b.c-triangulation, too, and the corresponding torus-equivariant
projective, crepant, full desingularization
\[
f_{\mathcal{T}_{\mathbb{H}_{d}}\left|  _{P}\right.  }:X(\mathbb{Z}%
^{d},\widehat{\Delta}_{P}\left(  \mathcal{T}_{\mathbb{H}_{d}}\left|
_{P}\right.  \right)  )\longrightarrow X(\mathbb{Z}^{d},\Delta_{P}%
)=U_{\tau_{P}}\
\]
possesses exceptional prime divisors $D_{n}=V\left(  \mathbb{R}_{\geq
0}n\right)  =X((\mathbb{Z}^{d})({\mathbb{R}_{\geq0}n}),$ \emph{Star}%
$(\mathbb{R}_{\geq0}n;\widehat{\Delta}_{P}(\mathcal{T}_{\mathbb{H}_{d}}\left|
_{P}\right.  )))$ for which
\[
D_{n}\cong\widehat{W}^{\left(  d\right)  },\ \ \text{\emph{if \ }}n\in\left(
\text{\emph{int}}\left(  P\right)  \right)  \cap\mathbb{Z}^{d},
\]
and
\[
D_{n}\cong\left\{  \text{\emph{a quasiprojective }}\left(  d-1\right)
\text{\emph{-dimensional subvariety of \ }}\widehat{W}^{\left(  d\right)
}\right\}  \ ,
\]
if $n\in
(\partial P \setminus $\emph{vert}$(P) )\cap \mathbb{Z}^{d}$, where
\emph{\ }$\widehat{W}^{\left(  d\right)  }$ denotes the projective toric Fano
manifold obtained by a torus-equivariant, crepant, projective, full resolution
$\widehat{W}^{\left(  d\right)  }\longrightarrow W^{\left(  d\right)  }$ of a
$\left(  d-1\right)  $-dimensional projective, toric Fano variety $W^{\left(
d\right)  }$ \emph{(}with at most Gorenstein singularities\emph{)}, induced by
the triangulation $\mathcal{T}_{\mathbb{H}_{d}}$\emph{.} In particular, as
projective variety, $W^{\left(  d\right)  }$ admits an embedding $W^{\left(
d\right)  }\hookrightarrow\mathbb{P}_{\mathbb{C}}^{\,d\left(  d-1\right)  }$
and has the degree\emph{\ }$\tbinom{2\,\left(  d-1\right)  }{d-1}$
w.r.t. this embedding.
\end{theorem}

\noindent\textit{Proof.} The first assertion is obvious. Let now $n\in
($int$(P))\cap\mathbb{Z}^{d}$. The \textit{star} of $n$
with respect to $\mathcal{T}_{\mathbb{H}_{d}}$ (in the sense of the theory of
simplicial complexes) is lattice equivalent to a pure simplicial complex
consisting of the triangulation $\mathcal{T}_{\mathbb{H}_{d}}\left|
_{\mathcal{Z}^{\left(  d\right)  }}\right.  $ of a $(d-1)$-dimensional
\textit{lattice zonotope} $\mathcal{Z}^{\left(  d\right)  }\subset
\mathbb{R}^{d-1}$ with $d$ ``zones'' into basic simplices. (Our reference
lattice here is that one being generated by aff$\left(  P\right)
\cap\mathbb{Z}^{d}$). The zonotope $\mathcal{Z}^{\left(  d\right)  }$ can be
also viewed as the convex hull of the union of the $\left[  -1,0\right]
$-cube with the $\left[  0,1\right]  $-cube, or, alternatively, as the
Minkowski sum of $d$ segments:
\begin{align*}
\mathcal{Z}^{\left(  d\right)  }  &  =\medskip\left\{  \mathbf{x}=\left(
x_{1},..,x_{d-1}\right)  \in\mathbb{R}^{d-1}\ \left|  \
\begin{array}
[c]{l}%
\left|  x_{i}\right|  \leq1,\medskip\ \forall i,\ \ 1\leq i\leq d-1,\ \text{
and}\\
\left|  x_{i}-x_{j}\right|  \leq1\text{, \ for all \ }i,j\text{, s.t. \ }1\leq
i,j\leq d-1
\end{array}
\right.  \text{ }\right\}  =\\
&  =\medskip\text{ conv}\left(  (\left[  -1,0\right]  ^{d-1})\cup(\left[
0,1\right]  ^{d-1})\right)  =\medskip\\
&  = \frac{1}{2} \left(\left[  -e_{1},e_{1}\right]  + \cdots + \left[
    -e_{d-1} ,e_{d-1}\right]  + \left[  -\left(  e_{1}+e_{2}+\cdots
+e_{d-1}\right)  ,e_{1}+e_{2}+\cdots+e_{d-1}\right] \right)
\end{align*}
with $\left\{  e_{1},e_{2},\ldots,e_{d-1}\right\}  $ denoting the standard
basis of unit vectors of $\mathbb{R}^{d-1}$. Obviously,
\[
\text{vert}\left(  \mathcal{Z}^{\left(  d\right)  }\right)  =\left\{
\pm\left(  e_{i_{1}}+e_{i_{2}}+\cdots+e_{i_{k}}\right)  \ \left|
\begin{array}
[c]{l}%
\text{for all subsets of indices }\\
1\leq i_{1}<i_{2}<\cdots<i_{k}\leq d-1\\
\text{and all \ }k\text{, \ }1\leq k\leq d-1
\end{array}
\right.  \right\}
\]
and
\[
\#(\text{vert}(\mathcal{Z}^{\left(  d\right)  }))=\#(\{\text{facets of
}\mathcal{Z}^{\left(  d\right)  \,\ast}\})=2\left[  \tbinom{d-1}{1}%
+\tbinom{d-1}{2}+\cdots+\tbinom{d-1}{d-1}\right]  =2\left(  2^{d-1}-1\right),
\]
where $\mathcal{Z}^{\left(  d\right)  \,\ast}$ denotes the polar of
$\mathcal{Z}^{\left(  d\right)  }$,
\[
\mathcal{Z}^{\left(  d\right)  \,\ast}=\left\{  \mathbf{y}=\left(
y_{1},\ldots,y_{d-1}\right)  \in(\mathbb{R}^{d-1})^{\vee}\ \left|  \
\begin{array}
[c]{l}%
\left|  y_{i_{1}}+y_{i_{2}}+\cdots+y_{i_{k}}\right|  \leq1,\ \text{for all }\\
1\leq i_{1}<i_{2}<\cdots<i_{k}\leq d-1\\
\text{and all \ }k\text{, \ }1\leq k\leq d-1
\end{array}
\right.  \text{ }\right\}  \ ,
\]
having the following $2\left(  d-1\right)  +2\tbinom{d-1}{2}=d\left(
d-1\right)  $ vertices:
\[
\text{vert}\left(  \mathcal{Z}^{\left(  d\right)  \,\ast}\right)  =\left\{
\pm e_{1}^{\vee},\pm e_{2}^{\vee},\ldots,\pm e_{d-1}^{\vee}\right\}
\cup\left\{  \pm\left(  e_{i}^{\vee}-e_{j}^{\vee}\right)  \ \left|  \ 1\leq
i<j\leq d-1\right.  \right\}  .
\]
($\left\{  e_{1}^{\vee},e_{2}^{\vee},\ldots,e_{d-1}^{\vee}\right\}  $ denotes
here the $\mathbb{R}$-basis of $(\mathbb{R}^{d-1})^{\vee}$ which is dual to
$\left\{  e_{1},e_{2},\ldots,e_{d-1}\right\}  $).\medskip\ \newline Note that
$\mathcal{Z}^{\left(  d\right)  \,}$, $\mathcal{Z}^{\left(  d\right)  \,\ast}$
are reflexive polytopes and can be inscribed in the cube $\left[  -1,1\right]
^{d-1}.$ (Figures $\mathbf{2}$ and $\mathbf{3}$ illustrate them for $d=3$ and
$d=4$). Let now $W^{\left(  d\right)  }$ be the $(d-1)$-dimensional projective
toric variety being associated to the normal fan of $\mathcal{Z}^{\left(
d\right)  \,\ast}$ (as in \S\ref{TORIC-G} \textsf{(k)}). This fan is strongly
polytopal because its rays are exactly the $1$-dimensional cones in
$\mathbb{R}^{d-1}$ determined by joining the origin with the vertices of
$\mathcal{Z}^{\left(  d\right)  \,}$. On the other hand, by Theorem
\ref{GOR-PR} and \cite[2.2.23]{Bat1}, we see that $W^{\left(  d\right)  }$ is
a Gorenstein toric Fano variety (which is \textit{singular} for $d\geq4$).
Moreover, the exceptional prime divisor $D_{n}=V\left(  \mathbb{R}_{\geq
0}n\right)
=X((\mathbb{Z}^{d})({\mathbb{R}_{\geq0}n}),\thinspace$Star$(\mathbb{R}%
_{\geq0}n;\widehat{\Delta}_{P}(\mathcal{T}_{\mathbb{H}_{d}}\left|
_{P}\right.  )))$ is analytically isomorphic to $\widehat{W}^{\left(
d\right)  }$, where $\widehat{W}^{\left(  d\right)  }$ is that projective,
toric Fano manifold which occurs as the overlying space of the
torus-equivariant, projective, crepant desingularization of $W^{\left(
d\right)  }$ induced by restricting the b.c.-triangulation $\mathcal{T}%
_{\mathbb{H}_{d}}$ onto $\mathcal{Z}^{\left(  d\right)  \,}$. Of course,
$W^{\left(  d\right)  }$ can be embedded into $\mathbb{P}_{\mathbb{C}%
}^{\,d\left(  d-1\right)  }$ via the map
\[
W^{\left(  d\right)  }\ni w\longmapsto\left[  w:\mathbf{e}\left(  e_{1}^{\vee
}\right)  \left(  w\right)  :\mathbf{e}\left(  -e_{1}^{\vee}\right)  \left(
w\right)  :\cdots:\mathbf{e}\left(  -(e_{d-2}^{\vee}-e_{d-1}^{\vee})\right)
\left(  w\right)  \right]  \in\mathbb{P}_{\mathbb{C}}^{\,d\left(  d-1\right)
},
\]
defined by evaluating the torus characters at the points of $\mathcal{Z}%
^{\left(  d\right)  \,\ast}\cap\mathbb{Z}^{d}=$ $\{\mathbf{0}\}\cup$
vert$(\mathcal{Z}^{\left(  d\right)  \,\ast})$ (cf. Fulton \cite[p.
69]{Fulton}), because $\mathcal{Z}^{\left(  d\right)  \,\ast}$ is the
lattice polytope
which is determined by the anticanonical divisor of $W^{\left(  d\right)  }$.
By \cite[Thm. 4.16, p. 36]{Sturmfels}, the degree of $W^{\left(  d\right)  }$
with respect to this embedding is equal to the normalized volume
Vol$_{\text{norm}}\left(  \mathcal{Z}^{\left(  d\right)  \,\ast}\right)  $ of
$\mathcal{Z}^{\left(  d\right)  \,\ast}$. It is worth mentioning that the
facets of $\mathcal{Z}^{\left(  d\right)  \,\ast}$ are exactly the
subpolytopes of the form
\begin{equation}
\pm\text{ conv}(\{e_{i_{1}}^{\vee},e_{i_{2}}^{\vee},\ldots,e_{i_{k}}^{\vee
}\}\cup\{e_{i_{1}}^{\vee}-e_{j_{1}}^{\vee},e_{i_{2}}^{\vee}-e_{j_{2}}^{\vee
},\ldots,e_{i_{k}}^{\vee}-e_{j_{k}}^{\vee}\}) \label{FACETZ}%
\end{equation}
for all subsets of indices $1\leq i_{1}<i_{2}<\cdots<i_{k}\leq d-1$, with
$1\leq k\leq d-1$, and all possible indices
\[
j_{1}\in\left\{  1,\ldots,d-1\right\}  \mathbb{r}\left\{  i_{1},i_{2}%
,\ldots,i_{k}\right\}  ,\ldots\ldots,j_{k}\in\left\{  1,\ldots,d-1\right\}
\mathbb{r}\left\{  i_{1},i_{2},\ldots,i_{k}\right\}  \ .
\]
Each facet (\ref{FACETZ}) is nothing but the direct product of a
$(d-(k+1))$-dimensional basic simplex with a $(k-1)$-dimensional basic
simplex, for any $k$, $1\leq k\leq d-1$, and consequently its relative,
normalized volume equals $\binom{d-(k+1)+k-1}{k-1}=\binom{d-2}{k-1}$. Hence,
since the normalized volume of a reflexive polytope is equal to the sum of the
relative, normalized volumes of its facets, we get
\[
\text{Vol}_{\text{norm}}\left(  \mathcal{Z}^{\left(  d\right)  \,\ast}\right)
=2\,\left[  \,\sum_{k=1}^{d-1}\ \tbinom{d-1}{k}\tbinom{d-2}{k-1}\,\right]
=\tfrac{2}{d-1}\,\left[  \,\sum_{k=1}^{d-1}\ \tbinom{d-1}{k}^{2}\cdot
k\,\right]  =\tbinom{2\,\left(  d-1\right)  }{d-1}\ .
\]
Finally, let us point out that if $n\in (\partial P \setminus
$vert$(P) )\cap \mathbb{Z}^{d}$, then by
construction Star$(\mathbb{R}_{\geq0}n;\widehat{\Delta}_{P}(\mathcal{T}%
_{\mathbb{H}_{d}}\left|  _{P}\right.  ))$ is a subfan of the fan induced by
the star of $n$ with respect to the \textit{entire} $\mathcal{T}%
_{\mathbb{H}_{d}} $. $D_{n}$ can be therefore viewed as a torus-invariant
non-compact subvariety of a projective toric variety which is analytically
isomorphic to the above defined $W^{\left(  d\right)  }$. $_{\Box}$

\vspace{1\baselineskip}
\begin{center}
\begin{minipage}{50mm}
{\centering $\mathcal{Z}^{(3)}$
\\[.5\baselineskip]
\epsfig{file=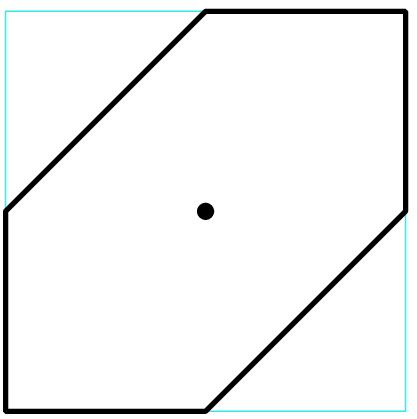}}
\end{minipage}
\hfill
\raisebox{0mm}{$\underleftrightarrow{
    \raisebox{2mm}{\text{  polarity  }}}$}
\hfill
\begin{minipage}{50mm}
{\centering $\mathcal{Z}^{(3) \thinspace *}$
\\[.5\baselineskip]
\epsfig{file=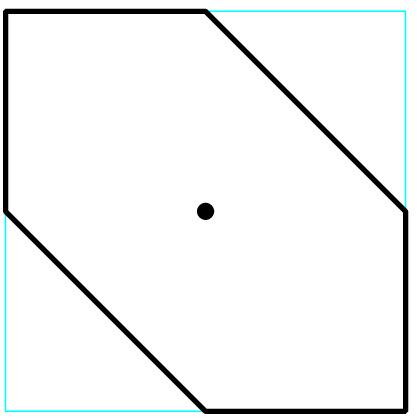}}
\end{minipage}
\\[\baselineskip]
after triangulating
\\[.5\baselineskip]
\epsfig{file=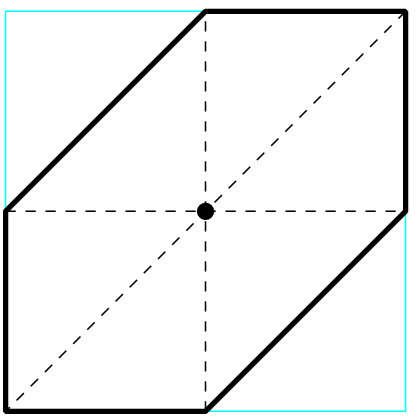}
\\[\baselineskip]
\textbf{Figure 2}
\end{center}

\begin{center}
\enlargethispage{100mm}
\begin{minipage}{70mm}
{\centering $\mathcal{Z}^{(4)}$
\\[.5\baselineskip]
\epsfig{file=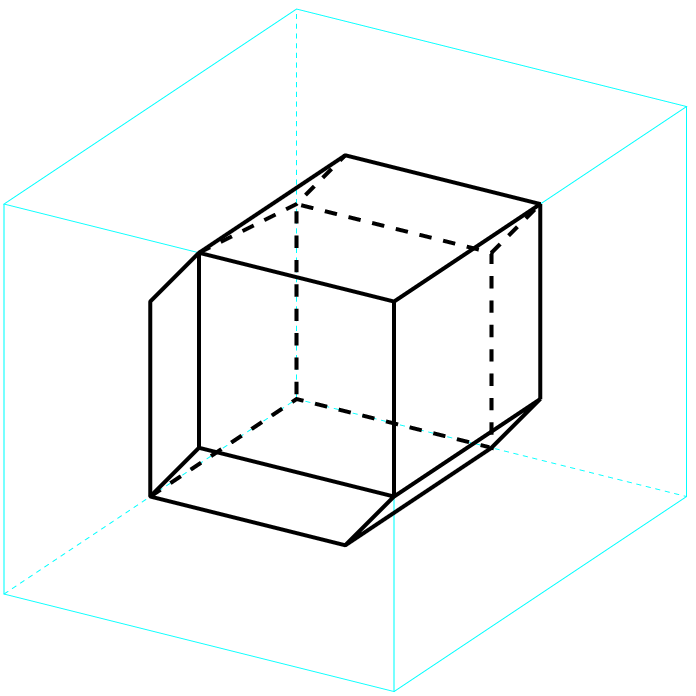}}
\end{minipage}
\hfill
\raisebox{0mm}{$\underleftrightarrow{
    \raisebox{2mm}{\text{  polarity  }}}$}
\hfill
\begin{minipage}{70mm}
{\centering $\mathcal{Z}^{(4) \thinspace *}$
\\[.5\baselineskip]
\epsfig{file=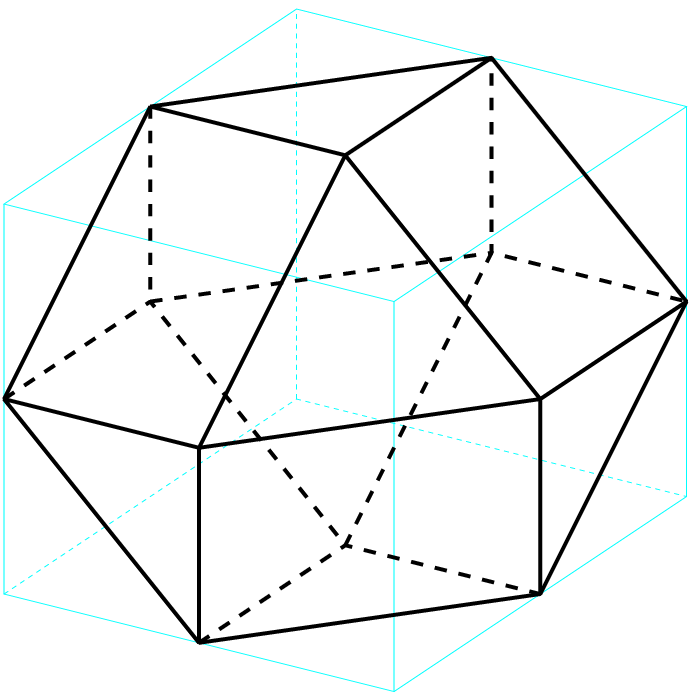}}
\end{minipage}
\\[\baselineskip]
after triangulating
\\[.5\baselineskip]
\epsfig{file=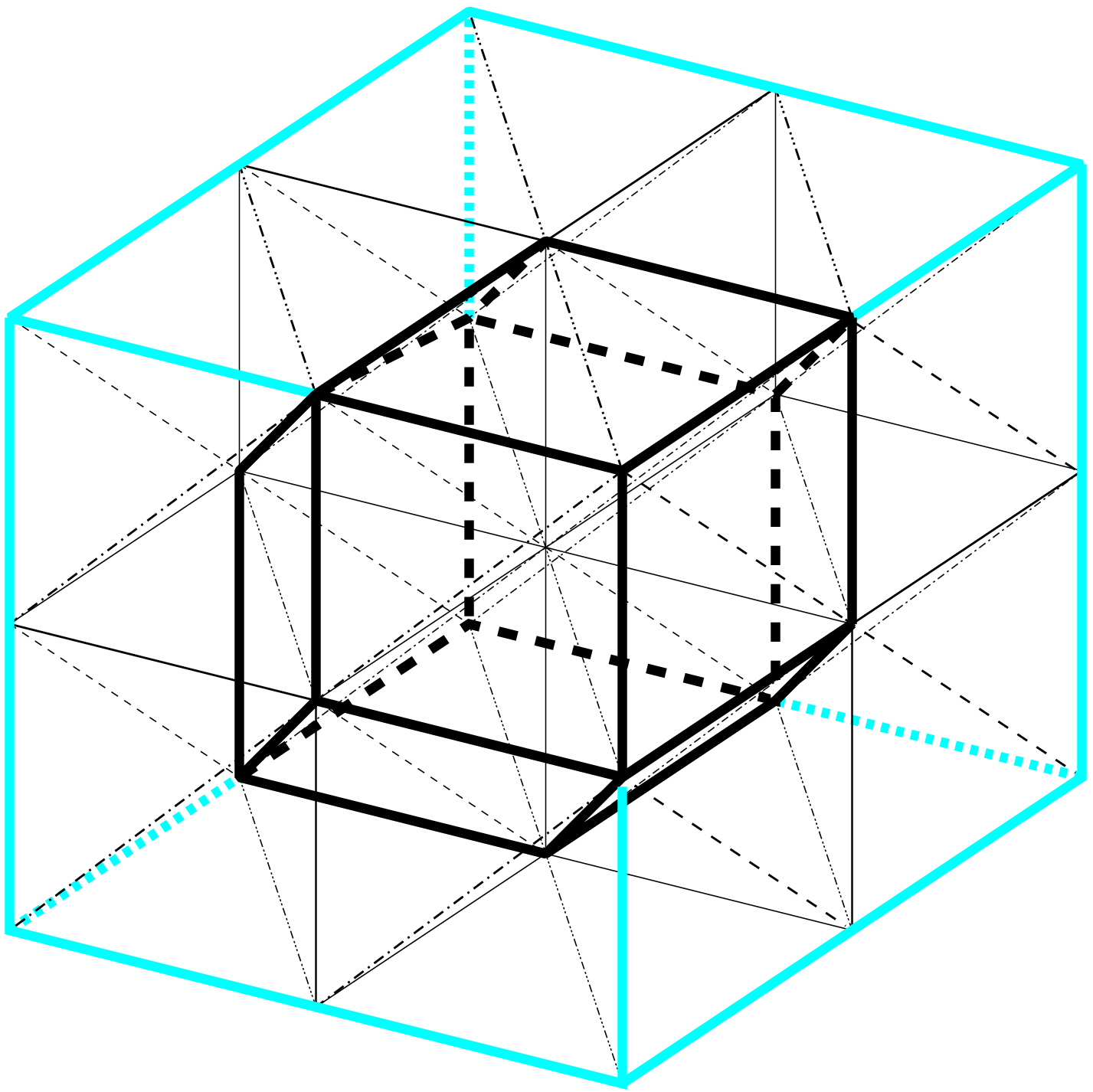}
\\[\baselineskip]
\textbf{Figure 3}
\end{center}

\pagebreak

\section{Nakajima's polytopes and Classification Theorem\label{NAK}}

\noindent Let $\mathbb{R}^{d}$ be the usual $d$-dimensional euclidean space,
$\mathbb{Z}^{d}$ the usual rectangular lattice in $\mathbb{R}^{d}$ and
$(\mathbb{Z}^{d})^{\vee}$ its dual lattice in $(\mathbb{R}^{d})^{\vee}=$
Hom$_{\mathbb{R}}(\mathbb{R}^{d},\mathbb{R)}$. From now on, we shall represent
the points of $\mathbb{R}^{d}$ by column vectors and the points of its dual
$(\mathbb{R}^{d})^{\vee}$ by row vectors.

\begin{definition}
\emph{A}\textit{\ sequence\emph{\ }of free parameters of length} $\ell$
\emph{(w.r.t. }$\mathbb{Z}^{d}$\emph{) is defined to be a finite sequence}
\[
\mathbf{m}:=\left(  m_{1},m_{2},\ldots,m_{\ell}\right)  ,\ \ \ 1\leq\ell\leq
d-1,
\]
\emph{consisting of vectors} $m_{i}:=\left(  m_{i,1},m_{i,2},\ldots
,m_{i,d}\right)  ,1\leq i\leq\ell,$ \emph{of} $\left(  \mathbb{Z}^{d}\right)
^{\vee}\mathbb{r}\left\{  \left(  0,\ldots,0\right)  \right\}  \ $ \emph{for
which }$m_{i,j}=0$ \emph{for all }$i$\emph{, }$1\leq i\leq\ell$\emph{, and all
}$j$\emph{, }$1\leq j\leq d$\emph{, with }$i<j$. \emph{As }$\left(  \ell\times
d\right)  $\emph{-matrix such an}\textit{\ }$\mathbf{m}$ \emph{has the
form:\smallskip}
\begin{equation}
\mathbf{m}=\left(
\begin{array}
[c]{cc}%
\begin{array}
[c]{cccccc}%
m_{1,1} & 0 & 0 & \cdots & \cdots & 0\\
m_{2,1} & m_{2,2} & 0 & \cdots & \cdots & 0\\
m_{3,1} & m_{3,2} & m_{3,3} & \cdots & \cdots & 0\\
\vdots & \vdots & \vdots & \ddots & \cdots & \vdots\\
m_{\ell-1,1} & m_{\ell-1,2} & m_{\ell-1,3} & \cdots & \ddots & 0\\
m_{\ell,1} & m_{\ell,2} & m_{\ell,3} & \cdots & \cdots &  m_{\ell,\ell}%
\end{array}
\!\!\smallskip & \underset{d-\ell\emph{\ zero-columns}}{\underbrace{%
\begin{array}
[c]{ccc}%
0 & \cdots & 0\smallskip\\
0 & \cdots & 0\smallskip\\
0 & \cdots & 0\smallskip\\
\vdots & \cdots & \vdots\\
0 & \cdots & 0\\
0 & \cdots & 0
\end{array}
}}%
\end{array}
\right)  \label{MATRIX-F}%
\end{equation}
\end{definition}

\begin{definition}
[Nakajima's polytopes]\label{NAKPOL}\emph{Fixing the dimension }$d$\emph{\ of
our reference space, we define the} \emph{polytopes }$\left\{  P_{\mathbf{m}%
}^{\left(  i\right)  }\subset\mathbf{\bar{H}}^{\left(  d\right)  }%
\mathbf{\ }\left|  \ i\in\mathbb{N}\text{\emph{, }}1\leq i\leq d\right.
\right\}  $ \emph{lying on }$\mathbf{\bar{H}}^{\left(  d\right)  }%
\mathbf{=}\left\{  \mathbf{x}=(x_{1},..,x_{d})^{\intercal}\in\mathbb{R}%
^{d}\ \left|  \ x_{1}=1\right.  \right\}  $ \emph{and being associated to an
``admissible'' free-parameter-sequence} $\mathbf{m}$ \emph{as in
(\ref{MATRIX-F}) }$\emph{w.r.t.\ }\mathbb{Z}^{d}$\emph{\ (with length }%
$\ell=i-1$\emph{, } \emph{for }$2\leq i\leq d$\emph{) by using induction on}
$i$\emph{; namely we define} $P_{\mathbf{m}}^{\left(  1\right)  }%
:=\{(1,\underset{\left(  d-1\right)  \text{\emph{-times}}}{\underbrace
{0,0,\ldots,0,0}})\}$\emph{, } \emph{and for }$2\leq i\leq d$\emph{,}
\[
P_{\mathbf{m}}^{\left(  i\right)  }:=\text{\emph{conv}}\left(  P_{\mathbf{m}%
}^{\left(  i-1\right)  }\cup\left\{  \left.  (\mathbf{x}^{\prime},\left\langle
m_{i-1},\mathbf{x}\right\rangle ,\underset{\left(  d-i\right)
\text{\emph{-times}}}{\underbrace{0,..,0}})^{\intercal}\ \right|
\ \mathbf{x}=(\underset{\underset{\mathbf{x}^{\prime}}{\mathbf{\parallel}}%
}{\underbrace{x_{1},x_{2},..,x_{i-1}}},\underset{\left(  d-i+1\right)
\text{\emph{-times}}}{\underbrace{0,..,0})^{\intercal}}\in P_{\mathbf{m}%
}^{\left(  i-1\right)  }\right\}  \right)  .\smallskip
\]
\emph{(}$P_{\mathbf{m}}^{\left(  i\right)  }$ \emph{is obviously }$\left(
i-1\right)  $\emph{-dimensional). For} $\mathbf{m}$ \emph{to be ``}%
\textit{admissible}\emph{'' means that}
\[
\left\langle m_{i-1},\mathbf{x}\right\rangle \geq0,\ \ \forall
\mathbf{x,\ \ \ x}=(x_{1},x_{2},\ldots,x_{i-1},\underset{\left(  d-i+1\right)
\text{\emph{-times}}}{\underbrace{0,\ldots,0})^{\intercal}}\in P_{\mathbf{m}%
}^{\left(  i-1\right)  }.
\]
\emph{Any lattice }$\left(  i-1\right)  $-\emph{polytope} $P$ \emph{which is
lattice equivalent to a }$P_{\mathbf{m}}^{\left(  i\right)  }$ \emph{(as
defined} \emph{above) will be called} \emph{a} \textit{Nakajima polytope}
$\emph{(w.r.t.\ }\mathbb{R}^{d}\emph{).}$
\end{definition}

\begin{example}
\label{Ex}\emph{(i) For }$i=d=1,$\emph{\ we have trivially }$P_{\mathbf{m}%
}^{\left(  1\right)  }=\left\{  1\right\}  .\smallskip$\emph{\newline (ii) For
}$d=2,\ \mathbf{m}=\left(  m_{1,1},0\right)  $\emph{\ we have }$P_{\mathbf{m}%
}^{\left(  1\right)  }=\left\{  \left(  1,0\right)  ^{\intercal}\right\}
$\emph{\ and}
\[
P_{\mathbf{m}}^{\left(  2\right)  }=\text{\emph{\ conv}}\left(  \left\{
\left(  1,0\right)  ^{\intercal}\right\}  \cup\left\{  \left(  1,\left\langle
m_{1},(1,0)\right\rangle \right)  ^{\intercal}\right\}  \right)
=\text{\emph{conv}}\left(  \left\{  \left(  1,0\right)  ^{\intercal}\right\}
\cup\left\{  \left(  1,m_{1,1}\right)  ^{\intercal}\right\}  \right)
\text{\emph{, \ }}m_{1,1}>0.
\]
\emph{(iii)} \ \emph{For }$i=d=3,$ \emph{and}
\[
\mathbf{m=}\left(
\begin{array}
[c]{ccc}%
m_{1,1} & 0 & 0\\
m_{2,1} & m_{2,2} & 0
\end{array}
\right)
\]
\emph{we obtain}
\[
P_{\mathbf{m}}^{\left(  3\right)  }=\text{\emph{\ conv}}\left(  \left\{
\left(  1,0,0\right)  ^{\intercal},\left(  1,m_{1,1},0\right)  ^{\intercal
},\left(  1,0,m_{2,1}\right)  ^{\intercal},\left(  1,m_{1,1},m_{2,1}%
+m_{1,1}m_{2,2}\right)  ^{\intercal}\right\}  \right)
\]
\emph{with}
\[
m_{1,1}>0,\ \ m_{2,1}\geq0,\ \ m_{2,1}+m_{1,1}m_{2,2}\geq0,\ \ \left(
m_{2,1},m_{2,2}\right)  \neq\left(  0,0\right)  .
\]
\newline \emph{(iv) Finally, for }$i=d=4,$ \emph{and}
\[
\mathbf{m=}\left(
\begin{array}
[c]{cccc}%
m_{1,1} & 0 & 0 & 0\\
m_{2,1} & m_{2,2} & 0 & 0\\
m_{3,1} & m_{3,2} & m_{3,3} & 0
\end{array}
\right)
\]
\emph{we get}
\[
P_{\mathbf{m}}^{\left(  4\right)  }=\text{\emph{conv}}\left(  \left\{
\begin{array}
[c]{c}%
\left(  1,0,0,0\right)  ^{\intercal},\left(  1,m_{1,1},0,0\right)
^{\intercal},\left(  1,0,m_{2,1},0\right)  ^{\intercal},\left(  1,m_{1,1}%
,m_{2,1}+m_{1,1}m_{2,2},0\right)  ^{\intercal},\smallskip\\
\left(  1,0,0,m_{3,1}\right)  ^{\intercal},\left(  1,m_{1,1},0,m_{3,1}%
+m_{1,1}m_{3,2}\right)  ^{\intercal},\left(  1,0,m_{2,1},m_{3,1}%
+m_{3,3}m_{2,1}\right)  ^{\intercal},\smallskip\\
\left(  1,m_{1,1},m_{2,1}+m_{1,1}m_{2,2},m_{3,1}+m_{3,2}m_{1,1}+m_{2,1}%
m_{3,3}+m_{1,1}m_{2,2}m_{3,3}\right)  ^{\intercal}%
\end{array}
\right\}  \right)
\]
\emph{with}
\[
\left\{
\begin{array}
[c]{l}%
m_{1,1}>0,\ m_{2,1}\geq0,\ m_{2,1}+m_{1,1}m_{2,2}\geq0,\ m_{3,1}%
\geq0,\ m_{3,1}+m_{1,1}m_{3,2}\geq0,\smallskip\\
m_{3,1}+m_{3,3}m_{2,1}\geq0,\ m_{3,1}+m_{3,2}m_{1,1}+m_{2,1}m_{3,3}%
+m_{1,1}m_{2,2}m_{3,3}\geq0,\smallskip\\
\ \left(  m_{2,1},m_{2,2}\right)  \neq\left(  0,0\right)  ,\ \ \left(
m_{3,1},m_{3,2},m_{3,3}\right)  \neq\left(  0,0,0\right)
\end{array}
\right.
\]
\emph{In the Figures \textbf{4} and \textbf{5} \ we illustrate the
  lattice polytopes }$P_{\mathbf{m}}^{\left(  3\right)  } , P_{\mathbf{m}}^{\left(  4\right)  }$
\emph{, respectively, for}
\[
\mathbf{m=}\left(
\begin{array}
[c]{ccc}%
2 & 0 & 0\\
2 & 1 & 0
\end{array}
\right)
\text{ \ \ \ \ \emph{and} \ \ \ \ }
\mathbf{m=}\left(
\begin{array}
[c]{rrrr}%
1 & 0 & 0 & 0\\
1 & 0 & 0 & 0\\
2 & -1 & -1 & 0
\end{array}
\right) \qquad .
\]
\end{example}

\begin{center}
\begin{minipage}{25mm}
{\centering
\epsfig{file=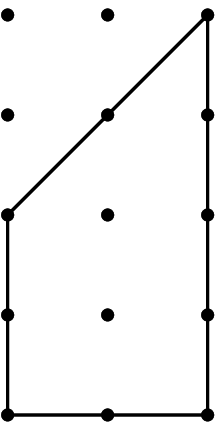}
\\[.5\baselineskip]
{\centering \textbf{Figure 4}}}
\end{minipage}
\hspace{40mm}
\begin{minipage}{28mm}
{\centering
\epsfig{file=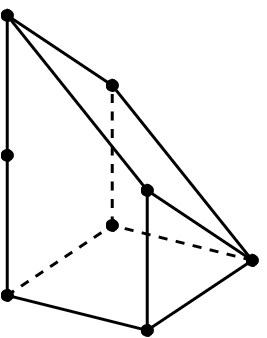}
\\[.5\baselineskip]
{\centering\textbf{Figure 5}}}
\end{minipage}
\end{center}

\begin{lemma}
[Description by inequalities]\label{INEQ}The $\left(  i-1\right)
$-dimensional polytope $P_{\mathbf{m}}^{\left(  i\right)  }$ \emph{(}%
$w.r.t.\emph{\ }\mathbb{R}^{d}$\emph{) }associated to an admissible
free-parameter-sequence $\mathbf{m}$\textbf{\ }can be written as a bounded
solution set of a finite system of linear inequalities as follows
\[
P_{\mathbf{m}}^{\left(  i\right)  }=\left\{  \mathbf{x}=(x_{1},..,x_{d}%
)^{\intercal}\in\mathbb{R}^{d}\ \left|
\begin{array}
[c]{l}%
x_{1}=1\text{\emph{,\ }}0\leq x_{j+1}\leq\,\left\langle m_{j},\mathbf{x}%
\right\rangle ,\text{ }\forall\text{\emph{\ }}j,\,1\leq j\leq i-1\medskip\\
\ \text{\emph{and\ \ \ }}x_{\mu}=0,\text{ }\forall\text{\emph{\ }}%
\mu,\ \ \ \,i+1\leq\mu\leq d
\end{array}
\right.  \right\}  .
\]
\end{lemma}

\noindent Another useful geometric description of Nakajima polytopes can be
provided by means of suitably cutted half-line prisms.

\begin{definition}
[Half-line prisms]\emph{Let} $\left(  i,d\right)
\in\mathbb{N}^{\thinspace 2}$ \emph{,
}$2\leq i\leq d$\emph{, and }$Q\subset \mathbf{\bar{H}}^{\left(
  d-1\right)}$
\emph{be a} $\left(  i-2\right)
$\emph{-dimensional polytope. We define the} 
\textit{half-line prism} $\mathbf{Pr}_{\text{\emph{hl}}}^{\left(  i\right)
}\left(  Q\right)  $ \emph{over} $Q$ \emph{as the }$\left(  i-1\right)
$-\emph{dimensional polyhedron}
\[
\mathbf{Pr}_{\text{\emph{hl}}}^{\left(  i\right)  }\left(  Q\right)
:=Q\times\mathbb{R}_{\geq0} = \left\{ (\mathbf{x},t) \in
  \mathbf{\bar{H}}^{\left( d \right)} \medspace | \medspace \mathbf{x}
  \in Q , t \in \mathbb{R}_{\ge 0} \right\}  ,
\]
\emph{and identify }$Q\subset\mathbf{\bar{H}}^{\left(  d-1\right)  }$
\emph{with }$Q\times\left\{  0\right\}  \subset\mathbf{\bar{H}}^{\left(
d-1\right)  }\times\left\{  0\right\}  \hookrightarrow\mathbf{\bar{H}%
}^{\left(  d\right)  } $\emph{.} \emph{(The only difference between}
$\mathbf{Pr}_{\text{\emph{hl}}}^{\left(  i\right)  }\left(  Q\right)  $
\emph{and a usual prism, is that the first one is} \emph{``open from above'').}
\end{definition}

\begin{lemma}
[Reduction Lemma]\label{REDUCTION}A lattice polytope $P\subset\mathbf{\bar{H}%
}^{\left(  d\right)  }\hookrightarrow\mathbb{R}^{d}$ is a Nakajima polytope of
dimension $i-1$\textbf{\ }$\emph{(w.r.t.\ }\mathbb{R}^{d}\emph{)}$ \textbf{iff
} $P\sim\left\{  \text{\emph{a lattice point in }}\mathbf{\bar{H}}^{\left(
d\right)  }\mathbf{\cap\,}\mathbb{Z}^{d}\right\}  $, for $i=1$, while for
$2\leq i\leq d$,
\[
P\sim(\mathbf{Pr}_{\text{\emph{hl}}}^{\left(  i\right)  }\left(  Q\right)
)\cap\left\{  \mathbf{x}=(1,x_{2},..,x_{d})^{\intercal}\in\mathbf{\bar{H}%
}^{\left(  d\right)  }\left|
\begin{array}
[c]{l}%
x_{i}\medskip\leq\sum_{j=1}^{i-1}\lambda_{j}x_{j}\text{ \ \emph{and\medskip}%
}\\
x_{\mu}=0,\text{ }\forall\text{\emph{\ }}\mu,\ \,i+1\leq\mu\leq d
\end{array}
\right.  \ \right\}  \ ,
\]
where the facet\emph{\ }$Q$ of the right-hand side is a Nakajima polytope of
dimension $i-2$ $\emph{(}$w.r.t.$\emph{\ }\mathbb{R}^{d-1}$, identified with
$\mathbb{R}^{d-1}\times\left\{  0\right\}  \hookrightarrow\mathbb{R}%
^{d}\emph{)},$ and $(\lambda_{1},\ldots,\lambda_{i-1},-1,\underset{\left(
d-i\right)  \text{\emph{-times}}}{\underbrace{0,0,\ldots,0,0})}\in
(\mathbb{Z}^{d})^{\vee}$ expresses a \emph{(}not identically zero\emph{) }
functional with non-negative values on $Q\times\left\{  0\right\}
\hookrightarrow\mathbf{\bar{H}}^{\left(  d\right)  }$.
\end{lemma}

\noindent\textit{Proof. }\ For $i=1$ there is nothing to be shown. Let
$i\in\left\{  2,\ldots,d\right\}  $. If $P\sim P_{\mathbf{m}}^{\left(
i\right)  }$ for $\mathbf{m}$ an admissible sequence of free parameters of
length $i-1$ (w.r.t. $\mathbb{Z}^{d}$), then
\[
P_{\mathbf{m}}^{\left(  i\right)  }=(\mathbf{Pr}_{\text{hl}}^{\left(
i\right)  }(P_{\overline{\mathbf{m}}}^{\left(  i-1\right)  }))\cap\left\{
\mathbf{x}=(1,x_{2},..,x_{d})^{\intercal}\in\mathbf{\bar{H}}^{\left(
d\right)  }\left|
\begin{array}
[c]{l}%
x_{i}\medskip\leq\sum_{j=1}^{i-1}\lambda_{j}x_{j}\text{ \ and\medskip}\\
x_{\mu}=0,\text{ }\forall\text{\emph{\ }}\mu,\ \,i+1\leq\mu\leq d
\end{array}
\right.  \ \right\}  \ ,
\]
where $P_{\overline{\mathbf{m}}}^{\left(  i-1\right)  }$ is a Nakajima
polytope of dimension $i-2$, with $P_{\overline{\mathbf{m}}}^{\left(
i-1\right)  }=\left\{  (1,0,...,0)\right\}  $ for $i=2$, and determined by the
admissible sequence of free parameters $\overline{\mathbf{m}}=\left(
m_{1},m_{2},\ldots,m_{i-2}\right)  $ of length $i-2$, for $i\geq3$, so that
$\ \mathbf{m}=\left(
\begin{array}
[c]{cc}%
\overline{\mathbf{m}}\ \ 0 & 0\\
m_{i-1} & 0
\end{array}
\right)  $ under the identification of $P_{\overline{\mathbf{m}}}^{\left(
i-1\right)  }$ with $P_{\overline{\mathbf{m}}}^{\left(  i-1\right)  }%
\times\left\{  0\right\}  \hookrightarrow P_{\mathbf{m}}^{\left(  i\right)  }%
$. Hence, it is enough to take $Q=P_{\overline{\mathbf{m}}}^{\left(
i-1\right)  }$ and $\lambda_{1}=m_{i-1,1}$, $\lambda_{2}=m_{i-1,2}$,
$\ldots$, $\lambda_{i-1}=m_{i-1,i-1}$. And conversely, having the
intersection of the
half-line prism $\mathbf{Pr}_{\text{hl}}^{\left(  i\right)  }\left(  Q\right)
$ with the half-space determined by the non-trivial non-negatively-valued
functional $(\lambda_{1},\ldots,\lambda_{i-1},-1,0,\ldots,0)\in(\mathbb{Z}%
^{d})^{\vee}$ on $Q\times\left\{  0\right\}  $ as our starting-point, we may
construct (by the backtracking method, i.e., by passing from the last to the
last but one row etc.) an admissible sequence $\mathbf{m}$ of free parameters
of length $i-1$ (w.r.t. $\mathbb{Z}^{d}$), such that $P\sim P_{\mathbf{m}%
}^{\left(  i\right)  }$. $_{\Box}\bigskip$

\noindent Now Nakajima's Classification Theorem \cite[Thm. 1.5, p.
86]{Nakajima} can be formulated as follows:

\begin{theorem}
[Nakajima's Classification of Toric L.C.I.'s]\label{Nak-thm}Let $N$ be a free
$\mathbb{Z}$-module of rank $r$, and $\sigma\subset N_{\mathbb{R}}$ \ a s.c.p.
cone of dimension $d\leq r$. Moreover, let $U_{\sigma}$ denote the affine
toric variety associated to $\sigma$, and $U_{\sigma^{\prime}}$ as in
\emph{\S3} $\mathsf{(a)}$. Then $U_{\sigma}$ is local complete
intersection\emph{\ }\textbf{if and only if }\ there exists an admissible
sequence $\mathbf{m}$ of free parameters of length $d-1$ $\emph{(}%
$w.r.t.$\emph{\ }\mathbb{Z}^{d}\emph{)}$, such that\smallskip\ for any
standard representative $U_{\tau_{P}}\cong U_{\sigma^{\prime}}$ of $U_{\sigma
}$ we have $P\sim P_{\mathbf{m}}^{\left(  d\right)  }$, i.e., $P$ is a
Nakajima $\left(  d-1\right)  $-dimensional polytope \emph{(}w.r.t.
$\mathbb{R}^{d}$\emph{).}
\end{theorem}

\begin{remark}
\label{REMNAK}\emph{(i) Theorem \ref{Nak-thm} was first proved in dimension
}$3$\ \emph{by Ishida \cite[Thm. 8.1, p. 136]{Ishida}. Previous
classification results, due to Watanabe \cite{Watanabe}, cover essentially
only the class of the }$\mathbb{Q}$\emph{-factorial toric l.c.i.'s in all
dimensions. (The term ``Watanabe simplex'' introduced in~\cite[5.13]{DHZ1}, 
can be used, up to lattice equivalence, as a synonym for a Nakajima
polytope which is simultaneously a simplex.)
\medskip\newline (ii) Obviously, }$U_{\sigma}$ \emph{is a
l.c.i.} $\Longleftrightarrow$ $U_{\sigma^{\prime}}\cong U_{\tau_{P}}$ \emph{is
a g.c.i. (Since in the setting of \ \cite{DHZ1}, it was always assumed that
}$d=r$\emph{, the abelian quotient ``g.c.i.''-spaces were abbreviated therein
simply as ``c.i.'s'').\medskip}\newline \emph{(iii) For }$P$\emph{\ a
non-basic Nakajima polytope,} $(U_{\tau_{P}},\emph{orb}\left(  \tau
_{P}\right)  )$ \emph{is a toric g.c.i.-singularity.\medskip\newline (iv) For
}$P$\emph{\ }$\ $\emph{a Nakajima }$\left(  d-1\right)  $\emph{-polytope and
}$\tau_{P}$ \emph{non-basic w.r.t. }$\mathbb{Z}^{d}$\emph{, orb}$\left(
\tau_{P}\right)  \in U_{\tau_{P}}$ \emph{has splitting codimension }%
$\varkappa$\emph{, with } $2\leq\varkappa\leq d-1$\emph{, iff }$P$ \emph{is
lattice-equivalent to the join }$\check{P}\star\mathbf{s}$ \emph{of } \emph{a
}$\left(  \varkappa-1\right)  $\emph{-dimensional (non-basic) Nakajima
polytope }$\check{P}$ \emph{with a basic }$\left(  d-\varkappa-1\right)
$\emph{-simplex }$\mathbf{s}$, \emph{which lie in adjacent lattice
hyperplanes, and }$\varkappa$\emph{ is minimal w.r.t. this
property}. \medskip\emph{\newline (v) It is easy for every
}$P\subset\mathbf{\bar{H}}^{\left(  d\right)  }$\emph{, with }$P\sim
P_{\mathbf{m}}^{\left(  d\right)  }$\emph{,} \emph{to verify that}
\begin{equation}
d\leq\#(\text{\emph{vert}}(P))\leq2^{d-1}\ \ \ \text{\emph{and\ \ \ }}%
d\leq\#(\{\text{\emph{facets of }}P\})\leq2\left(  d-1\right)  \label{BOUNDS}%
\end{equation}
\end{remark}

\section{Proof of Main Theorem and of Koszul-property\label{M-KOS}}

\noindent\textsf{(a) }By Prop. \ref{CREPANT} (ii), (iii), and Thm.
\ref{Nak-thm}, our Main Theorem \ref{MAIN} is equivalent to the following:

\begin{theorem}
\label{MAIN2}All Nakajima polytopes admit b.c.-triangulations in all dimensions.
\end{theorem}

\noindent Our proof of \ref{MAIN2} relies on the construction of the desired
lattice triangulations via the classical ``pulling operation'' of vertices of
a point configuration and the ``Key-Lemma'' \ref{KEY}.

\begin{definition}
[Pulling vertices]\emph{Consider a finite set of points} $\mathcal{V}=\left\{
\mathbf{v}_{1},\ldots,\mathbf{v}_{k}\right\}  \subset\mathbb{R}^{d}$ \emph{and
let }$\mathcal{S}=\{P_{1},P_{2},...,P_{\nu}\}$ \emph{denote} \emph{a}
\emph{polytopal subdivision} \emph{of conv}$\left(  \mathcal{V}\right)
$\emph{\ with vert}$\left(  \mathcal{S}\right)  \subseteq$ $\mathcal{V}$.
\emph{\ For any }$i\in\left\{  1,\ldots,k\right\}  $\emph{, we define a
refinement} $\mathsf{p}_{\mathbf{v}_{i}}\left(  \mathcal{S}\right)  $
\emph{of} $\mathcal{S}$ \ \emph{(called the }\textit{pulling} \emph{of}
$\mathbf{v}_{i}$\emph{) as follows:\medskip\newline (i)} $\mathsf{p}%
_{\mathbf{v}_{i}}\left(  \mathcal{S}\right)  $ \emph{contains all} $P_{j}%
$\emph{'s } \emph{for which} $\mathbf{v}_{i}\notin$ $P_{j}$\emph{,
and\medskip\newline (ii) if} $\mathbf{v}_{i}\in P_{j}$\emph{, then}
$\mathsf{p}_{\mathbf{v}_{i}}\left(  \mathcal{S}\right)  $ \emph{contains all
the polytopes having the form} \emph{conv}$\left(  F\cup\mathbf{v}_{i}\right)
$\emph{, with} $F$ \emph{a facet of } $P_{j}$ \emph{such that} $\mathbf{v}%
_{i}\notin F$.\smallskip
\end{definition}

\begin{lemma}
\label{LA}The refinements of $\mathcal{S}$ $\ $obtained by pulling
\textbf{all} the points of $\mathcal{V}$ \ \emph{(}in arbitrary order\emph{)
}are triangulations of \emph{conv}$\left(  \mathcal{V}\right)$ with
vertex set $\mathcal{V}.$
\end{lemma}

\noindent\textit{Proof.} This is an easy exercise\textit{\
  }(cf. \cite[\S2]{Lee1}).\textit{\ }$_{\Box}$

\begin{example}
[Realization of pullings by ``full flags'']\label{FLAGS}\emph{Suppose that
}$\mathcal{V}=\left\{  \mathbf{v}_{1},\ldots,\mathbf{v}_{k}\right\}
\subset\mathbb{R}^{d}$ \emph{is the set of vertices of a }$d$%
\emph{-dimensional polytope }$P$\emph{. In this special case, the
triangulation $\mathcal{T}$ \ obtained after performing the pulling operation
for all points of $\mathcal{V}$ has a nice geometric realization due to
Stanley (see \cite[\S 1]{Stanley2}). For every face }$F$\emph{\ of }$P$
\emph{define} $\mathbf{v}(F):=\mathbf{v}_{j}$\emph{, where }$j:=$
\emph{min}$\{i\,\left|  \,\mathbf{v}_{i}\right.  \in F\}$\emph{. A
}\textit{full flag }\emph{of }$P$ \emph{is a chain }$\mathcal{F}$ \emph{of
faces} $F_{0}\subset F_{1}\subset F_{2}\subset\cdots\subset F_{d}=P$\emph{,
such that dim}$\left(  F_{i}\right)  =i$\emph{, for all }$i$\emph{, }$0\leq
i\leq d$\emph{, and }$\mathbf{v}(F_{i-1})\neq\mathbf{v}(F_{i})$\emph{, for all
}$i$\emph{, }$1\leq i\leq d$\emph{. For any full flag }$\mathcal{F}$
\emph{define }$\mathbf{v}(\mathcal{F}):=\{\mathbf{v}\left(  F_{0}\right)
,..,\mathbf{v}\left(  F_{d}\right)  \}$\emph{. Then the simplices of the
triangulation}
\begin{equation}
\mathcal{T}=\mathsf{p}_{\mathbf{v}_{k}}(\mathsf{p}_{\mathbf{v}_{k-1}}%
\cdots\,\cdots(\mathsf{p}_{\mathbf{v}_{2}}(\mathsf{p}_{\mathbf{v}_{1}}\left(
\{P\}\right)  )))\label{PULLST}%
\end{equation}
\emph{constructed by pulling the points of }$\mathcal{V}$ \ \emph{in the order
}$\mathbf{v}_{1},\ldots,\mathbf{v}_{k}$ \emph{(by starting from the
``trivial''} \emph{subdivision }$\{P\}$\emph{) are exactly the elements of
the set} $\{$\emph{conv}$(\mathbf{v}(\mathcal{F}))\ \left|  \ \mathcal{F\,}%
\text{\emph{a full flag of \thinspace}}P\right.  \}.$\emph{\ }
\end{example}

\noindent The pulling of a vertex point of a polytope (and therefore
triangulations of the form (\ref{PULLST}) too) are known to be coherent (cf.
Lee \cite[p. 448]{Lee1}, and \cite[p. 275]{Lee2}). In fact, a slightly
stronger statement is also true:

\begin{lemma}
[Coherency preservation by pulling operation]\label{LB}Let $\mathcal{V}$ \ be
a set of finite points in $\mathbb{R}^{d}$\emph{\ }and $\mathcal{S}%
=\{P_{1},P_{2},...,P_{\nu}\}$ an arbitrary coherent polytopal subdivision of
\emph{conv}$\left(  \mathcal{V}\right)  $\emph{\ }with \emph{vert}%
$(\mathcal{S})\subseteq\mathcal{V}$. Then the refinement $\mathsf{p}%
_{\mathbf{v}_{0}}\left(  \mathcal{S}\right)  $ of $\mathcal{S}$, \ for a
$\mathbf{v}_{0}\in$ $\mathcal{V}$, forms a coherent polytopal subdivision of
\ \emph{conv}$\left(  \mathcal{V}\right)  $.
\end{lemma}

\noindent\textit{Proof. }Let $\omega: \mathcal{V}\rightarrow
\mathbb{R}$ be a height function on $\mathcal{V}$, for which
$\mathcal{S}=\mathcal{S}_{\omega}$ (in the notation of Lemma \ref{HEIGHTS}).
The (maximal dimensional) polytopes $P_{1},...,P_{\nu}$ of $\mathcal{S}$ are
images of the lower facets of the polytope $Q_{\omega}:=$ conv$(\{(\mathbf{v}%
,\omega(\mathbf{v}))\in\mathbb{R}^{d+1}\ \left|  \ \mathbf{v}\in
\text{vert}(\mathcal{S})\right.  \})$ under the projection $\pi:\mathbb{R}%
^{d+1}\rightarrow\mathbb{R}^{d}$ w.r.t. the last coordinate. Let $\mathfrak{l}%
_{1},\ldots,\mathfrak{l}_{\nu}\in(\mathbb{R}^{d})^{\vee}$ denote functionals for
which
\[
\mathfrak{l}_{i}\left(  \mathbf{v}\right)  -\omega\left(  \mathbf{v}\right)  \leq
c_{i},\ \ \ \text{for all }\mathbf{v}\in\text{vert}(\mathcal{S}),\ \ 1\leq
i\leq\nu,
\]
and appropriate $c_{i}$'s $\in\mathbb{R}$, so that the equality is valid only
for $\mathbf{v}$ $\in P_{i}$, i.e., so that $\mathfrak{l}_{i}$ ``determines''
$P_{i}$. We define
\[
t_{0}:=\text{max}\{\mathfrak{l}_{i}\left(  \mathbf{v}_{0}\right)  -c_{i}\ \left|
\ 1\leq i\leq\nu\right.  \}=\text{ min}\{t\in\mathbb{R}:(\mathbf{v}_{0},t)\in
Q_{\omega}\}\ .
\]
Obviously, $(\mathbf{v}_{0},t_{0})$ belongs to the boundary of $Q_{\omega}$
(see Fig. $\mathbf{6}$). Without loss of generality, we may assume that
$\nu\geq2$ and that the maximum of the differences $\mathfrak{l}_{i}\left(
\mathbf{v}_{0}\right)  -c_{i}$ is achieved for $1\leq i\leq j$, but not for
$j+1\leq i\leq\nu$, for some index $j\in\{1,2,..,\nu-1\}$. (This means that
$\mathbf{v}_{0}\in P_{i}$ for $1\leq i\leq j$, but $\mathbf{v}_{0}\notin
P_{i}$ for $j+1\leq i\leq\nu$). We define
\[
\omega^{\prime}:\text{vert}\left(  \mathsf{p}_{\mathbf{v}_{0}}\left(
\mathcal{S}\right)  \right)  \longrightarrow\mathbb{R},\text{ \ \ with
\ \ }\omega^{\prime}\left(  \mathbf{v}\right)  :=\left\{
\begin{array}
[c]{ll}%
\omega\left(  \mathbf{v}\right)  , & \text{if\ \ }\mathbf{v\medskip}%
\neq\mathbf{v}_{0}\\
t_{0}-\varepsilon, & \text{if\ \ }\mathbf{v}=\mathbf{v}_{0}%
\end{array}
\right.
\]
where $\varepsilon>0$ is chosen to be small enough for ensuring that
$\mathfrak{l}_{i}\left(  \mathbf{v}_{0}\right)  -\omega^{\prime}\left(
\mathbf{v}_{0}\right)  <c_{i}$, for all $i$, $j+1\leq i\leq\nu$, i.e., for
setting $(\mathbf{v}_{0},\omega^{\prime}\left(  \mathbf{v}_{0}\right)  )$ into
a ``general position'' w.r.t. the lower envelope of $Q_{\omega}$. If we define
$Q_{\omega^{\prime}}:=$ conv$(\{(\mathbf{v},\omega^{\prime}(\mathbf{v}%
))\in\mathbb{R}^{d+1}\ \left|  \ \mathbf{v}\in\text{vert}(\mathsf{p}%
_{\mathbf{v}_{0}}\left(  \mathcal{S}\right)  )\right.  \})$, then the faces of
$Q_{\omega^{\prime}}$ are the faces of $Q_{\omega}$ which do not contain
$(\mathbf{v}_{0},\omega^{\prime}\left(  \mathbf{v}_{0}\right)  )$, together
with the faces of type conv$(\{(\mathbf{v}_{0},\omega^{\prime}\left(
\mathbf{v}_{0}\right))  \}\cup F)$, where $F$ is a face of some facet of
$Q_{\omega}$ containing $(\mathbf{v}_{0},t_{0})$. Thus, the projection $\pi
($conv$(\{(\mathbf{v}_{0},\omega^{\prime}\left(  \mathbf{v}_{0}\right)
)\}\cup F))$ onto $\mathbb{R}^{d}$ is a subset of the projection of this
facet, and $\mathsf{p}_{\mathbf{v}_{0}}\left(  \mathcal{S}\right)  $ is
exactly the polytopal subdivision of conv$\left(  \mathcal{V}\right)  $
induced by the above defined height function $\omega^{\prime}$ (again in the
sense of Lemma \ref{HEIGHTS}).\textit{\ }$_{\Box}$

\begin{center}
\begin{picture}(0,0)%
\epsfig{file=pull.pstex}%
\end{picture}%
\setlength{\unitlength}{1816sp}%
\begingroup\makeatletter\ifx\SetFigFont\undefined%
\gdef\SetFigFont#1#2#3#4#5{%
  \reset@font\fontsize{#1}{#2pt}%
  \fontfamily{#3}\fontseries{#4}\fontshape{#5}%
  \selectfont}%
\fi\endgroup%
\begin{picture}(7224,6140)(1189,-9489)
\put(4276,-6436){\makebox(0,0)[lb]{\smash{\SetFigFont{11}{13.2}{\rmdefault}{\mddefault}{\updefault}$(\mathbf{v}_0,t_0)$}}}
\put(4351,-7861){\makebox(0,0)[lb]{\smash{\SetFigFont{11}{13.2}{\rmdefault}{\mddefault}{\updefault}$(\mathbf{v}_0,\omega'(\mathbf{v}_0))$}}}
\end{picture}

\\[.5\baselineskip]
\textbf{Figure 6}
\end{center}

\begin{theorem}
\label{COHTH}Let $\mathcal{V}$ be a set of finite points in $\mathbb{R}^{d}%
$\emph{\ }and $\mathcal{S}$ a coherent polytopal subdivision of \emph{conv}%
$\left(  \mathcal{V}\right)  $ with\emph{\ vert}$\left(  \mathcal{S}\right)
\subseteq$ $\mathcal{V}$. Then $\mathcal{S}$ can be always refined to a
coherent triangulation\emph{\ }$\mathcal{T}$ \ of \emph{conv}$\left(
\mathcal{V}\right)  $, such that $\emph{vert}\left(  \mathcal{T}\right)
=\mathcal{V}.$
\end{theorem}

\noindent\textit{Proof. }By sequentially pulling \textit{all} the points of
$\mathcal{V}$ (in arbitrary order), and by using Lemmas \ref{LA} and \ref{LB},
we may always construct such a coherent triangulation $\mathcal{T}$ of
conv$\left(  \mathcal{V}\right)  $. $_{\Box}$

\begin{lemma}
[Key-Lemma]\label{KEY}Let $\mathbf{s}\subset\mathbb{R}^{d}$ be a $\left(
d-2\right)  $-dimensional simplex with
\[
\emph{vert}\left(  \mathbf{s}\right)  =\left\{  \mathbf{v}_{1},\mathbf{v}%
_{2},\ldots,\mathbf{v}_{d-1}\right\}  \subset\mathbf{\bar{H}}^{\left(
d-1\right)  }\cap\mathbb{Z}^{d}\hookrightarrow\mathbf{\bar{H}}^{\left(
d\right)  }\cap\mathbb{Z}^{d},\ \ \ d\geq2,
\]
and let $\mathbf{s}^{\prime}\subset\mathbb{R}^{d}$ denote a $\left(
d-1\right)  $-dimensional simplex with
\[
\emph{vert}\left(  \mathbf{s}^{\prime}\right)  =\left\{  \mathbf{v}%
_{1}^{\prime},\mathbf{v}_{2}^{\prime},\ldots,\mathbf{v}_{d-1}^{\prime
},\mathbf{v}_{d}^{\prime}\right\}  \subset(\mathbf{Pr}_{\text{\emph{hl}}%
}^{\left(  d\right)  }\left(  \mathbf{s}\right)  )\cap\mathbb{Z}%
^{d}\,\hookrightarrow\mathbf{\bar{H}}^{\left(  d\right)  }\cap\mathbb{Z}^{d}.
\]
If $\mathbf{s}$ is a basic simplex and $\mathbf{s}^{\prime}$ an elementary
simplex \emph{(}w.r.t. $\mathbb{Z}^{d}$\emph{)}, then $\mathbf{s}^{\prime}$
has to be basic too.
\end{lemma}

\noindent\textit{Proof. }The property for a lattice simplex to be elementary
or basic remains invariant among all the members of its lattice equivalence
class. Since $\mathbf{s}$ is embedded into $\mathbf{\bar{H}}^{\left(
d\right)  }\hookrightarrow\mathbb{R}^{d}$ and is assumed to be basic, there
exists an affine integral transformation $\Phi:\mathbb{R}^{d}\rightarrow
\mathbb{R}^{d}$, such that $\Phi\left(  \mathbf{v}_{1}\right)  =e_{1}$,
$\Phi\left(  \mathbf{v}_{i}\right)  =e_{1}+e_{i}$, for all $i$, $2\leq i\leq
d-1$, and $\Phi(\mathbf{\bar{H}}^{\left(  d\right)  }\mathbf{)}=\mathbf{\bar
{H}}^{\left(  d\right)  }$, where $\left\{  e_{1},\ldots,e_{d-1}%
,e_{d}\right\}  $ is the standard basis of unit vectors of $\mathbb{R}^{d}$
(!). This induces the lattice equivalence:
\[
\mathbf{s\sim\,}\widetilde{\mathbf{s}},\ \ \ \text{with \ \ \ }\widetilde
{\mathbf{s}}:=\text{ conv}\left(  \left\{  e_{1},e_{1}+e_{2},e_{1}%
+e_{3},\ldots,e_{1}+e_{d-1}\right\}  \right)  \ .
\]
We define $\widetilde{\mathbf{s}}^{\prime}:=\Phi\left(  \mathbf{s}^{\prime
}\right)  $. Since $\widetilde{\mathbf{s}}^{\prime}$ is $\left(  d-1\right)
$-dimensional and
\[
\text{vert}(\widetilde{\mathbf{s}}^{\prime})=\left\{  \Phi\left(
\mathbf{v}_{1}^{\prime}\right)  ,\Phi\left(  \mathbf{v}_{2}^{\prime}\right)
,\ldots,\Phi\left(  \mathbf{v}_{d-1}^{\prime}\right)  ,\Phi\left(
\mathbf{v}_{d}^{\prime}\right)  \right\}  \subset(\mathbf{Pr}_{\text{hl}%
}^{\left(  d\right)  }\left(  \widetilde{\mathbf{s}}\right)  )\cap
\mathbb{Z}^{d},
\]
$d-1$ among the $d$ vertices of $\widetilde{\mathbf{s}}^{\prime}$, e.g., up to
enumeration of indices, say the first $d-1$ ones, must be of the form
\[
\Phi\left(  \mathbf{v}_{1}^{\prime}\right)  =e_{1}+\gamma_{1}\cdot
e_{d},\ \ \ \ \Phi\left(  \mathbf{v}_{i}^{\prime}\right)  =e_{1}+e_{i}%
+\gamma_{i}\cdot e_{d},\ \ \forall i,\ \ 2\leq i\leq d-1,
\]
for a $\left(  d-1\right)  $-tuple $\left(  \gamma_{1},\gamma_{2}%
,\ldots,\gamma_{d-1}\right)  \in\left(  \mathbb{Z}_{\geq0}\right)
^{d-1}$. Moreover,
\[
\Phi\left(  \mathbf{v}_{d}^{\prime}\right)  \in\{\mathbb{R}_{\geq0}\left(
e_{1}+\gamma_{1}\cdot e_{d}\right)  \cap\mathbb{Z}^{d}\}\cup\ (%
{\textstyle\bigcup\limits_{i=2}^{d-1}}
\ \{\mathbb{R}_{\geq0}\left(  e_{1}+e_{i}+\gamma_{i}\cdot e_{d}\right)
\cap\mathbb{Z}^{d}\})\ .
\]
On the other hand, $\widetilde{\mathbf{s}}^{\prime}\cap\mathbb{Z}^{d}=$
vert$(\widetilde{\mathbf{s}}^{\prime})$ means that
\[
\Phi\left(  \mathbf{v}_{d}^{\prime}\right)  \in\{e_{1}+\left(  \gamma_{1}%
\pm1\right)  \cdot e_{d}\}\cup\ (%
{\textstyle\bigcup\limits_{i=2}^{d-1}}
\ \{e_{1}+e_{i}+\left(  \gamma_{i}\pm1\right)  \cdot e_{d}\})\ ,
\]
because otherwise $\widetilde{\mathbf{s}}^{\prime}$ could not be elementary
(cf. Figure $\mathbf{7}$). Now since the matrices
\[
\left(
\begin{array}
[c]{ccccccc}%
1 & 1 & 1 & 1 & \cdots & 1 & 1\\
0 & 1 & 0 & 0 & \cdots & 0 & 0\\
0 & 0 & 1 & 0 & \cdots & 0 & 0\\
\vdots & \vdots & \vdots & \ddots &  & \vdots & \vdots\\
0 & 0 & 0 & 0 & \ddots & 0 & 0\\
0 & 0 & 0 & 0 & \cdots & 1 & 0\\
\gamma_{1} & \gamma_{2} & \gamma_{3} & \gamma_{4} & \cdots & \gamma_{d-1} &
\gamma_{1}\pm1
\end{array}
\right)  \text{ \ ,\ }\left(
\begin{array}
[c]{lllllll}%
1 & 1 & 1 & 1 & \cdots & 1 & 1\\
0 & 1 & 0 & 0 & \cdots & 0 & 0\\
\vdots & \vdots & \vdots & \vdots & \cdots & \vdots & \vdots\\
\vdots & \vdots & \vdots & \ddots &  & \vdots & 1\ \text{{\scriptsize (i-th
row)}}\\
\vdots & \vdots & \vdots & \vdots & \ddots & \vdots & \vdots\\
0 & 0 & 0 & 0 & \cdots & 1 & 0\\
\gamma_{1} & \gamma_{2} & \gamma_{3} & \gamma_{4} & \cdots & \gamma_{d-1} &
\gamma_{i}\pm1
\end{array}
\!\!\right)
\]
for all $i$, $2\leq i\leq d-1$, have always determinants equal to $\pm1$, we
obtain mult$(\widetilde{\mathbf{s}}^{\prime};\mathbb{Z}^{d})=1$, and
consequently both $\widetilde{\mathbf{s}}^{\prime}$ and $\mathbf{s}^{\prime}$
have to be basic simplices. $_{\Box}$\newline 

\begin{center}
\epsfig{file=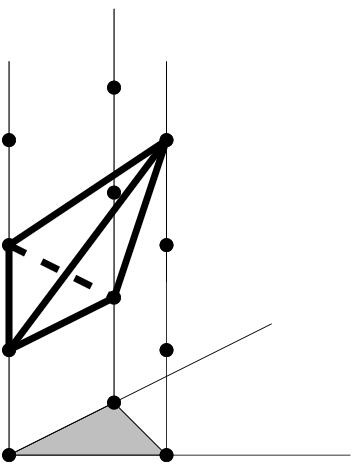}
\\[\baselineskip]
\textbf{Figure 7}
\end{center}

\noindent\textbf{Proof of Theorem \ref{MAIN2}: }By Thm. \ref{Nak-thm} it
suffices to show that all $\left(  d-1\right)  $-dimensional Nakajima
polytopes $P\subset\mathbf{\bar{H}}^{\left(  d\right)  }\subset\mathbb{R}^{d}
$ (with $P\sim P_{\mathbf{m}}^{\left(  d\right)  }$, where $\mathbf{m}$ is an
admissible sequence of length $d-1$) admit b.c.-triangulations. We shall use
induction on the dimension $d$ of the ambient space. For $d\leq3$ this is
obviously trivial. The proof will take place for any fixed $d\geq4$ by
assuming that the assertion is true for $d-1$\textbf{. }By Reduction Lemma
\ref{REDUCTION}, we may write $P$ as the intersection
\[
P=(\mathbf{Pr}_{\text{hl}}^{\left(  d\right)  }\left(  Q\right)  )\cap\left\{
\mathbf{x}=(1,x_{2},..,x_{d})^{\intercal}\in\mathbf{\bar{H}}^{\left(
d\right)  }\ \left|  \ x_{d}\leq\,\right.  \sum_{j=1}^{d-1}\lambda_{j}%
x_{j}\right\}
\]
of a half-line prism over Nakajima $\left(  d-2\right)  $-polytope
$Q\subset\mathbf{\bar{H}}^{\left(  d-1\right)  }\hookrightarrow\mathbf{\bar
{H}}^{\left(  d\right)  }\subset\mathbb{R}^{d}$ with a half-space in
$\mathbf{\bar{H}}^{\left(  d\right)  }$ determined by a non-trivial,
non-negatively-valued functional $(\lambda_{1},\ldots,\lambda_{i-1}%
,-1,0,\ldots,0)\in(\mathbb{Z}^{d})^{\vee}$ on $Q\times\left\{
  0\right\}  $. By
induction hypothesis, $Q$ possesses a coherent triangulation, say $\mathfrak{T}$,
with vert$\left(  \mathfrak{T}\right)  =Q\cap\mathbb{Z}^{d}$, into \textit{basic}
simplices\textbf{. }This means that
\[
\mathcal{S}_{\mathfrak{T}}:=\bigcup_{\text{all simplices }\mathbf{s}\in\mathfrak{T}%
}P_{\mathbf{s}},\text{ where \ }P_{\mathbf{s}}:=(\mathbf{Pr}_{\text{hl}%
}^{\left(  d\right)  }\left(  \mathbf{s}\right)  )\cap\left\{  \mathbf{x}%
=(1,x_{2},..,x_{d})^{\intercal}\in\mathbf{\bar{H}}^{\left(  d\right)
}\ \left|  \ x_{d}\leq\,\right.  \sum_{j=1}^{d-1}\lambda_{j}x_{j}\right\}  ,
\]
forms a polytopal lattice subdivision of $P$ into polytopes constructed by the
half-line prisms over all the simplices of $\mathfrak{T}$.\medskip\ \newline
$\bullet$ The polytopal subdivision $\mathcal{S}_{\mathfrak{T}}$ itself is
coherent. Indeed, if the coherent triangulation $\mathfrak{T}$ \ of $Q$ is induced
by a height function $\omega:$ vert$(\mathfrak{T})\rightarrow\mathbb{R} $ (as in
Lemma \ref{HEIGHTS}), then $\mathcal{S}_{\mathfrak{T}}$ will be induced by the
height function $\omega^{\prime}:$ vert$(\mathcal{S}_{\mathfrak{T}})\rightarrow
\mathbb{R}$ defined by
\[
\omega^{\prime}\left(  \mathbf{v}^{\prime}\right)  :=\omega\left(
\mathbf{v}\right)  ,\text{ \ for all \ }\mathbf{v}^{\prime}\in\text{vert}%
(\mathcal{S}_{\mathfrak{T}}),\ \mathbf{v}^{\prime}=\left(  \mathbf{v},t\right)
,\ \mathbf{v}\in\text{vert}(\mathfrak{T}),\ \ t\in\mathbb{Z}_{\geq0}\ .
\]
\newline $\bullet$ As we mentioned in Theorem \ref{COHTH}, pulling
sequentially all the points of $\mathcal{V}:=\left|  \mathcal{S}_{\mathfrak{T}%
}\right|  \cap$ $\mathbb{Z}^{d}$ (in arbitrary order), we arrive at a coherent
triangulation $\mathcal{T}$ of $P$, which is simultaneously a maximal lattice
triangulation.\medskip\ \newline $\bullet$ To show that $\mathcal{T}$ is a
b.c.-triangulation w.r.t. $\mathbb{Z}^{d}$, it is therefore enough to verify
its ``basicness''. Since $\mathcal{T}$ is by construction a refinement of
$\mathcal{S}_{\mathfrak{T}}$ (see \ref{PSUB} (iv)), all subtriangulations
$\left\{  \mathcal{T}\left|  _{P_{\mathbf{s}}}\right.  :\ \mathbf{s}\text{
simplices of }\mathfrak{T}\right\}  $ obtained by the restrictions of
$\mathcal{T}$ \ onto $P_{\mathbf{s}}$'s have to be maximal lattice
triangulations too. As all the simplices of them are elementary with vertices
belonging to the set of lattice points of half-line prisms over basic
simplices, we prove that all these subtriangulations have to be basic by
applying Lemma \ref{KEY}. Since these subtriangulations fit together to give
$\mathcal{T}$, $\mathcal{T}$ has to be basic as well. This completes the proof
of Theorems \ref{MAIN2} and \ref{MAIN}. $_{\Box}$

\begin{example}
\emph{Fixing b.c.-triangulations }$\mathfrak{T}$ \emph{of the 
``bases'' of the Nakajima polytopes which were
  shown in Figures \textbf{4} and \textbf{5}, we construct in Figures
  \textbf{8} and \textbf{9}, respectively, the subdivisions
  }$\mathcal{S}_{\mathfrak{T}}$ \emph{and afterwards} \emph{two}
\emph{b.c.-triangulations} $\mathcal{T}$ \emph{\ by pulling
  vertices. More precisely, in Figure \textbf{8} we pull the
  available points in the order }%
$\emph{(1,0,2)}^{\intercal}$ \emph{, }
$\emph{(1,1,1)}^{\intercal}$ \emph{, }
$\emph{(1,1,2)}^{\intercal}$ \emph{, }
$\emph{(1,1,0)}^{\intercal}$ \emph{, }
$\emph{(1,2,1)}^{\intercal}$ \emph{, }
$\emph{(1,2,2)}^{\intercal}$ \emph{, }
$\emph{(1,2,3)}^{\intercal}$ \emph{(and the remaining ones in
  arbitrary order).}
\\[.5\baselineskip]
\emph{In Figure \textbf{9 }we pull the
points in the order }
$\emph{(1,0,0,1)}^{\intercal}$ \emph{and }
$\emph{(1,1,0,1)}^{\intercal}$ \emph{(and the remaining ones in
  arbitrary order). The obtained subdivision is again a b.c.-triangulation.}
\end{example}

\begin{center}
\begin{minipage}{25mm}
{\centering
\epsfig{file=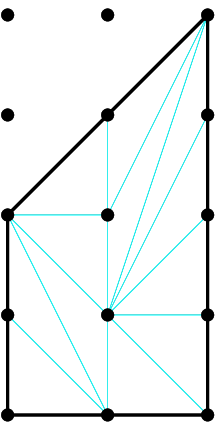}
\\[.5\baselineskip]
{\centering \textbf{Figure 8}}}
\end{minipage}
\hspace{40mm}
\begin{minipage}{28mm}
{\centering
\epsfig{file=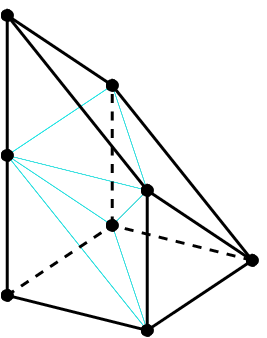}
\\[.5\baselineskip]
{\centering \textbf{Figure 9}}}
\end{minipage}
\end{center}

\begin{remark}
\emph{Theorem \ref{MAIN} has various applications to global geometrical
constructions. For instance, the Calabi-Yau varieties which arise from
(compactified, non-degenerate) hypersurfaces or ideal-theoretic complete
intersections of hypersurfaces embedded into compact toric Fano varieties, and
have at most l.c.i.-singularities, admit crepant, full, global
desingularizations in all dimensions (cf. \cite{Bat1,BB,BD}).\smallskip}
\end{remark}

\noindent\textsf{(b) }Another application of the proof of our Main
Theorem \ref{MAIN} is of 
purely algebraic nature and is related to the so-called\textit{\ Koszul
property }of graded algebras. (We restrict ourselves to graded algebras
defined over the field $\mathbb{C}$ of complex numbers).

\begin{definition}
[Koszul $\mathbb{C}$-algebras]\emph{A graded }$\mathbb{C}$\emph{-algebra }$R$
\emph{is called a }\textit{Koszul algebra} \emph{if }$\mathbb{C}$
\emph{(regarded as the }$R$\emph{-module }$R/\mathfrak{m}$ \emph{for }$\mathfrak{m}$
\emph{a maximal homogeneous ideal) has a linear free
resolution (in the sense of homological algebra), i.e., if there exists an
exact sequence }
\[
\cdots\longrightarrow\mathfrak{R}_{i+1}\overset{\varphi_{i+1}}{\longrightarrow
}\mathfrak{R}_{i}\overset{\varphi_{i}}{\longrightarrow}\cdots\overset{\varphi_{2}%
}{\longrightarrow}\mathfrak{R}_{1}\overset{\varphi_{1}}{\longrightarrow}%
\mathfrak{R}_{0}\longrightarrow R/\mathfrak{m}\longrightarrow0
\]
\emph{of graded free }$R$\emph{-modules all of whose matrices (determined by
the }$\varphi_{i}$\emph{'s) have entries which are linear forms (i.e., forms
of degree }$1$\emph{). Every Koszul algebra is generated by its component of
degree }$1$ \emph{and is defined by relations of degree }$2$.
\end{definition}

\begin{definition}
[``Non-faces'']\emph{Let }$\mathcal{V}$ \emph{be a finite set of points in}
$\mathbb{R}^{d}$ \emph{and }$\mathcal{T}$ \emph{a triangulation of
conv}$(\mathcal{V})$. \emph{A simplex whose vertices belong to }$\mathcal{V}$
$\ $\emph{but itself does not belong to $\mathcal{T}$ is defined to be a}
\textit{non-face} \emph{of} $\mathcal{T}$. \emph{A }\textit{minimal non-face}
\emph{of} $\mathcal{T}$ \emph{is a non-face of} $\mathcal{T}$ \emph{which is
minimal with respect to the} \emph{inclusion.}
\end{definition}

\begin{proposition}
[Koszulness and b.c.-triangulations]\label{KOS-T}If a $(d-1)$-dimensional
lattice polytope $P\subset\mathbf{\bar{H}}^{\left(  d\right)  }\hookrightarrow
\mathbb{R}^{d}$ \emph{(}w.r.t. $\mathbb{Z}^{d}$\emph{) }admits a
b.c.-triangulation whose minimal non-faces are $1$-dimensional, then
$R_{P}=\mathbb{C}\left[  \tau_{P}\cap\mathbb{Z}^{d}\right]  $ is a Koszul algebra.
\end{proposition}

\noindent\textit{Proof.} See Bruns, Gubeladze \& Trung \cite[2.1.3., p.
142]{BGT}. $_{\Box}$

\begin{proposition}
[From Nakajima polytopes to Koszulness]The coordinate rings \emph{(}%
$\mathbb{C}$-algebras\emph{)} $R_{P}=\mathbb{C}\left[  \tau_{P}\cap
\mathbb{Z}^{d}\right]  $ of the affine toric varieties $U_{\tau_{P}^{\vee}}$
being associated to the duals of the cones $\tau_{P}$ are Koszul for all
$(d-1)$-dimensional Nakajima polytopes $P$.
\end{proposition}

\noindent\textit{Proof. }By proposition \ref{KOS-T} it suffices to prove that
the b.c.-triangulations of such a $P$ which were constructed in the proof of
Theorem \ref{MAIN2} have exclusively $1$-dimensional minimal non-faces. We
shall again use induction on $d$. Assume that the assertion is true for $d-1$.
Let $\pi:\mathbf{\bar{H}}^{\left(  d\right)  }\rightarrow\mathbf{\bar{H}%
}^{\left(  d-1\right)  }$ denote the projection w.r.t. the last coordinate,
and $\mathcal{T}$ a b.c.-triangulation of a $(d-1)$-dimensional Nakajima
polytope $P\subset\mathbf{\bar{H}}^{\left(  d\right)  }$ induced by extending
a b.c.-triangulation $\mathfrak{T}$ of a $(d-2)$-dimensional Nakajima polytope
$Q\subset\mathbf{\bar{H}}^{\left(  d-1\right)  }$ as in the proof of
\ref{MAIN2}.\medskip\ \newline $\bullet$ \textit{Fact.} By construction, each
face of $\mathcal{T}$ is mapped by $\pi$ onto a face of $\mathfrak{T}$%
.\medskip\ \newline Now choose an arbitrary \textit{non-face} $\mathbf{s}$ of
$\mathcal{T}$ of dimension $\geq2$. It is enough to show that $\mathbf{s}$
contains an $1$-dimensional non-face of $\mathcal{T}$. \ We examine the two
possible cases separately:\medskip\newline (i) If the projection $\pi\left(
\mathbf{s}\right)  $ of $\mathbf{s}$ is a \textit{face} of the
b.c.-triangulation $\mathfrak{T}$ of $Q$, we consider the simplex $\mathbf{s}%
^{\prime}$ of $\mathcal{T}$ which contains the barycenter $\mathbf{bar}\left(
\mathbf{s}\right)  $ of $\mathbf{s}$ in its relative interior. (Such a simplex
$\mathbf{s}^{\prime}$ always exists, though it might be of dimension strictly
smaller than that of $\mathbf{s}$). Since both $\pi\left(  \mathbf{s}^{\prime
}\right)  $ and $\pi\left(  \mathbf{s}\right)  $ are faces of $\mathfrak{T}$, and
$\pi\left(  \mathbf{bar}\left(  \mathbf{s}\right)  \right)  $ belongs to the
intersection of their relative interiors, we have $\pi\left(  \mathbf{s}%
^{\prime}\right)  =\pi\left(  \mathbf{s}\right)  $. (Any point of $\left|
\mathfrak{T}\right|  $ belongs to the relative interior of exactly one simplex of
$\mathfrak{T}$). For each vertex $\mathbf{u}\in\pi\left(  \mathbf{s}^{\prime
}\right)  =\pi\left(  \mathbf{s}\right)  $, we define:
\[
t_{\mathbf{s}}^{\,\text{max}}\left(  \mathbf{u}\right)  :=\text{max}\left\{
t\in\mathbb{Z}_{\geq0}\ \left|  \ \left(  \mathbf{u},t\right)  \in
\text{vert}\left(  \mathbf{s}\right)  \right.  \right\}  ,\ \ t_{\mathbf{s}%
^{\prime}}^{\,\text{max}}\left(  \mathbf{u}\right)  :=\text{max}\left\{
t\in\mathbb{Z}_{\geq0}\ \left|  \ \left(  \mathbf{u},t\right)  \in
\text{vert}\left(  \mathbf{s}^{\prime}\right)  \right.  \right\}
\]
and
\[
t_{\mathbf{s}}^{\,\text{min}}\left(  \mathbf{u}\right)  :=\text{min}\left\{
t\in\mathbb{Z}_{\geq0}\ \left|  \ \left(  \mathbf{u},t\right)  \in
\text{vert}\left(  \mathbf{s}\right)  \right.  \right\}  ,\ \ t_{\mathbf{s}%
^{\prime}}^{\,\text{min}}\left(  \mathbf{u}\right)  :=\text{min}\left\{
t\in\mathbb{Z}_{\geq0}\ \left|  \ \left(  \mathbf{u},t\right)  \in
\text{vert}\left(  \mathbf{s}^{\prime}\right)  \right.  \right\}  ,
\]
respectively. Since $\mathbf{s}$ is a non-face and $\mathbf{s}^{\prime}$ a
face of $\mathcal{T}$, $\mathbf{s}$ cannot be contained in $\mathbf{s}%
^{\prime}$; so there must be at least one vertex $\mathbf{v}_{0}%
=(\mathbf{u}_{0},t_{0})\in\mathbf{s}\mathbb{r}\mathbf{s}^{\prime}$ of
$\mathcal{T}$, $\mathbf{u}_{0}\in$ vert$(\pi(\mathbf{s}))$, for which
\[
\text{\textit{either\ \ \ \ }}t_{0}>t_{\mathbf{s}^{\prime}}^{\,\text{max}%
}\left(  \mathbf{u}_{0}\right)  \ \ \ (\ast)\ \ \ \ \text{\textit{or}%
}\mathit{\ \ \ \ \ }t_{0}<t_{\mathbf{s}^{\prime}}^{\,\text{min}}\left(
\mathbf{u}_{0}\right)  \ \ \ (\ast\ast)
\]
\newline $\bullet$ \textit{Claim }A\textit{.} In case $(\ast)$ (resp. in case
$(\ast\ast)$) there is at least one vertex $\mathbf{u}_{\blacktriangledown}$
of $\pi\left(  \mathbf{s}^{\prime}\right)  =\pi\left(  \mathbf{s}\right)  $,
such that $t_{\mathbf{s}}^{\,\text{min}}\left(  \mathbf{u}_{\blacktriangledown
}\right)  <t_{\mathbf{s}^{\prime}}^{\,\text{max}}\left(  \mathbf{u}%
_{\blacktriangledown}\right)  $ (resp. $t_{\mathbf{s}}^{\,\text{max}}\left(
\mathbf{u}_{\blacktriangledown}\right)  >t_{\mathbf{s}^{\prime}}%
^{\,\text{min}}\left(  \mathbf{u}_{\blacktriangledown}\right)  $%
).\medskip\ \newline $\bullet$ \ \textit{Proof of Claim }A\textit{. }The proof
will be done only for the case $(\ast)$ because case $(\ast\ast)$ can be
treated similarly. Suppose that $t_{\mathbf{s}}^{\,\text{min}}\left(
\mathbf{u}\right)  \geq t_{\mathbf{s}^{\prime}}^{\,\text{max}}\left(
\mathbf{u}\right)  $, for all vertices $\mathbf{u}\in\pi\left(  \mathbf{s}%
^{\prime}\right)  =\pi\left(  \mathbf{s}\right)  $. Let
\[
\text{vert}\left(  \mathbf{s}\right)  =\left\{  \mathbf{v}_{i}=\left(
1,v_{i,1},v_{i,2},\ldots,v_{i,d-1}\right)  \in\mathbf{\bar{H}}^{\left(
d\right)  }\cap\mathbb{Z}^{d}\ \left|  \ 1\leq i\leq\text{ dim}\left(
\mathbf{s}\right)  +1\right.  \right\}
\]
be an enumeration of the vertex set of $\mathbf{s}$, and $\mathbf{bar}\left(
\mathbf{s}\right)  =(1,b_{1},b_{2},\ldots,b_{d-1})$ the coordinates of the
barycenter of $\mathbf{s}$ with
\[
\mathbf{bar}\left(  \mathbf{s}\right)  =\frac{1}{\text{dim}\left(
\mathbf{s}\right)  +1}\sum_{i=1}^{\text{dim}\left(  \mathbf{s}\right)
+1}\,\mathbf{v}_{i}\text{,\ \ i.e., \ \ \ }b_{j}=\frac{1}{\text{dim}\left(
\mathbf{s}\right)  +1}\sum_{i=1}^{\text{dim}\left(  \mathbf{s}\right)
+1}\,v_{i,j},\ \forall j,\ 1\leq j\leq d-1.
\]
The projection of $\mathbf{bar}\left(  \mathbf{s}\right)  $ equals
\[
\pi\left(  \mathbf{bar}\left(  \mathbf{s}\right)  \right)  =(1,b_{1}%
,b_{2},\ldots,b_{d-2})=\frac{1}{\text{dim}\left(  \mathbf{s}\right)  +1}%
\sum_{i=1}^{\text{dim}\left(  \mathbf{s}\right)  +1}\,\pi\left(
\mathbf{v}_{i}\right)  =\sum_{\mathbf{u}\in\text{vert}\left(  \pi\left(
\mathbf{s}\right)  \right)  }\mathfrak{r}\left(  \mathbf{u}\right)  \ \mathbf{u},
\]
where
\[
\mathfrak{r}\left(  \mathbf{u}\right)  :=\tfrac{1}{\text{dim}\left(
\mathbf{s}\right)  +1}\ \#\left\{  \pi^{-1}\left(  \mathbf{u}\right)
\cap\text{ vert}\left(  \mathbf{s}\right)  \right\}  =\tfrac{1}{\text{dim}%
\left(  \mathbf{s}\right)  +1}\ \#\left\{  \text{all }\mathbf{v}_{i}\text{'s
mapped onto }\mathbf{u}\text{ under }\pi\right\}  \ .
\]
Now let $\mathbf{v}_{\bullet}=(\mathbf{u}_{\bullet},t_{\bullet})$ denote an
arbitrary point of $\mathbf{s}^{\prime}$ with $\mathbf{u}_{\bullet}=\pi\left(
\mathbf{bar}\left(  \mathbf{s}\right)  \right)  $. Obviously,
\[
t_{\bullet}\leq\sum_{\mathbf{u}\in\text{vert}\left(  \pi\left(  \mathbf{s}%
\right)  \right)  }\mathfrak{r}\left(  \mathbf{u}\right)  \,t_{\mathbf{s}^{\prime
}}^{\,\text{max}}\left(  \mathbf{u}\right)  \leq\sum_{\mathbf{u}\in
\text{vert}\left(  \pi\left(  \mathbf{s}\right)  \right)  }\mathfrak{r}\left(
\mathbf{u}\right)  \,t_{\mathbf{s}}^{\,\text{min}}\left(  \mathbf{u}\right)
\leq\tfrac{1}{\text{dim}\left(  \mathbf{s}\right)  +1}\sum_{i=1}%
^{\text{dim}\left(  \mathbf{s}\right)  +1}\,v_{i,d-1}=b_{d-1}\ .
\]
If the last inequality were not strict, we would conclude that
\[
v_{i,d-1}=t_{\mathbf{s}}^{\,\text{min}}\left(
\mathbf{u}\right)  ,\ \forall i,\ 1\leq i\leq\text{ dim}\left(  \mathbf{s}%
\right)  +1,\text{ \ and \ \ }\forall\mathbf{u},\ \ \ \mathbf{u\in}\text{
vert}\left(  \pi\left(  \mathbf{s}\right)  \right)  .
\]
Since $\mathbf{u}_{0}\mathbf{\in}$ vert$\left(  \pi\left(  \mathbf{s}\right)
\right)  $ and $t_{0}=t_{\mathbf{s}}^{\,\text{min}}\left(  \mathbf{u}%
_{0}\right)  >t_{\mathbf{s}^{\prime}}^{\,\text{max}}\left(  \mathbf{u}%
_{0}\right)  $, this would mean that the second inequality is necessarily
strict. Hence, in each case, either the second or the third inequality has to
be \textit{strict}. This implies that the last coordinate of \textit{all}
points of $\mathbf{s}^{\prime}$ having the point $\pi\left(  \mathbf{bar}%
\left(  \mathbf{s}\right)  \right)  $ as their projection under $\pi$ is
$<b_{d-1}$, and therefore $\mathbf{bar}\left(  \mathbf{s}\right)
\notin\mathbf{s}^{\prime}$, which contradicts our initial assumption.\medskip
\ \newline $\bullet$ \textit{Claim }B\textit{.} conv$(\{\mathbf{v}%
_{0},\mathbf{(\mathbf{u}_{\blacktriangledown}},t_{\mathbf{s}}^{\,\text{min}%
}\mathbf{\left(  \mathbf{u}_{\blacktriangledown}\right)  )}\})$ in case
$(\ast)$ (resp. conv$(\{\mathbf{v}_{0},\mathbf{(\mathbf{u}_{\blacktriangledown
}},t_{\mathbf{s}}^{\,\text{max}}\mathbf{\left(  \mathbf{u}_{\blacktriangledown
}\right)  )}\})$ in case $(\ast\ast)$) is indeed an $1$-dimensional non-face
of $\mathcal{T}$.\medskip\newline $\bullet$ \ \textit{Proof of Claim
}B\textit{. }If it were a face of $\mathcal{T}$, then it would obviously
possess non-empty intersection with the face conv$(\{(\mathbf{u}%
_{0},t_{\mathbf{s}^{\prime}}^{\,\text{max}}\left(  \mathbf{u}_{0}\right)
),\mathbf{(\mathbf{u}_{\blacktriangledown}},t_{\mathbf{s}^{\prime}%
}^{\,\text{max}}\mathbf{\left(  \mathbf{u}_{\blacktriangledown}\right)  )}\})$
(resp. with the face conv$(\{(\mathbf{u}_{0},t_{\mathbf{s}^{\prime}%
}^{\,\text{min}}\left(  \mathbf{u}_{0}\right)  ),\mathbf{(\mathbf{u}%
_{\blacktriangledown}},t_{\mathbf{s}^{\prime}}^{\,\text{min}}\mathbf{\left(
\mathbf{u}_{\blacktriangledown}\right)  )}\})$). But this would mean that
$\mathcal{T}$ cannot be a triangulation.\medskip\newline (ii) Suppose now that
$\pi\left(  \mathbf{s}\right)  $ is a \textit{non-face} of the
b.c.-triangulation $\mathfrak{T}$ of $Q$. In this case, by induction hypothesis,
$\pi\left(  \mathbf{s}\right)  $ contains an $1$-dimensional minimal non-face
of $\mathfrak{T}$ , say conv$(\{\mathbf{u},\mathbf{u}^{\prime}\})$. Then both
$\pi^{-1}(\mathbf{u})$ and $\pi^{-1}(\mathbf{u}^{\prime})$ have to be faces of
$\mathcal{T}$ (cf. \cite[7.10]{Ziegler}), and for any two vertices
$\mathbf{v},\mathbf{v}^{\prime}$ of $\mathbf{s}$, with $\mathbf{v}\in\pi
^{-1}(\mathbf{u})\subset\mathbf{s}$ and $\mathbf{v}^{\prime}\in\pi
^{-1}(\mathbf{u}^{\prime})\subset\mathbf{s}$, conv$(\{\mathbf{v}%
,\mathbf{v}^{\prime}\})$ constitutes necessarily an $1$-dimensional non-face
of $\mathcal{T}$\textbf{\ }(by the above mentioned fact). This completes the
proof. $_{\Box}$

\begin{remark}
\emph{In fact we have shown the stronger statement that }$R_{P}\cong
\mathbb{C}\left[  \mathsf{T}_{1},\mathsf{T}_{2},\ldots,\mathsf{T}%
_{\#(P\cap\mathbb{Z}^{d})}\right]  \,/\,I_{P}$\emph{, where the binomial
ideal }$I_{P}$\emph{\ has a Gr\"{o}bner basis of degree }$2$\emph{\ (cf.
\cite{BGT}).}
\end{remark}

\section{On the computation of cohomology group dimensions\label{COHOMOLOGY}}

\noindent To compute the non-trivial (even) cohomology group dimensions of the
overlying spaces of crepant, full resolutions of toric l.c.i.-singularities we
need some basic concepts from enumerative combinatorics
(cf. \cite[\S4.6]{Stanley3}).\medskip

\noindent Let $N$ be a free $\mathbb{Z}$-module, $P\subset N_{\mathbb{R}}$ \ a
lattice polytope of dimension $k$ w.r.t. $N$, and $\nu$ a positive integer.
Let $\mathbf{Ehr}\left(  P,\nu\right)  :=\#\,\left(  \nu\,P\cap N_{P}\right)
=\sum_{j=0}^{k}\mathbf{a}_{j}\left(  P\right)  \nu^{j}\in\mathbb{Q}\left[
\nu\right]  $ denote the \textit{Ehrhart polynomial }of $P$ with $N_{P}$ the
affine sublattice aff$\left(  P\right)  \cap N$ of $N$, and
\[
\mathfrak{Ehr}\left(  P;q\right)  :=1+\sum_{\nu=1}^{\infty}\ \mathbf{Ehr}\left(
P,\nu\right)  \ q^{\nu}\in\mathbb{Q}_{\,}\left[  \!\left[  q\right]
\!\right]
\]
\textit{\ }the corresponding \textit{Ehrhart series}. Writing $\mathfrak{Ehr}%
\left(  P;q\right)  $ as
\[
\mathfrak{Ehr}\left(  P;q\right)  =\frac{\ \mathbf{\deltafett}_{0}\left(  P\right)
+\ \deltafett_{1}\left(  P\right)  \ q+\cdots+\ \deltafett%
_{k-1}\left(  P\right)  \ q^{k-1}+\mathbf{\ \deltafett}_{k}\left(  P\right)
\ q^{k}}{\left(  1-q\right)  ^{k+1}}%
\]
we get the so-called $\deltafett$-\textit{vector} $\deltafett\left(
P\right)  =\left(  \deltafett_{0}\left(  P\right)  ,\deltafett%
_{1}\left(  P\right)  ,\ldots,\deltafett_{k-1}\left(  P\right)
,\deltafett_{k}\left(  P\right)  \right)  $ of $P.$ (We should mention
that both $\mathbf{a}_{j}\left(  P\right)  $'s and $\deltafett_{j}\left(
P\right)  $'s are invariant under lattice equivalence).

\begin{lemma}
For all $j$, $0\leq j\leq k$, the $j$-th coordinate of the\textbf{\ }%
$\deltafett$-\textit{vector} of $P$ is given by the formula\emph{:}
\begin{equation}
\deltafett_{j}\left(  P\right)  =%
{\displaystyle\sum\limits_{i=0}^{k}}
\mathbf{\ }\left(
{\displaystyle\sum\limits_{\xi=0}^{j}}
\ \left(  -1\right)  ^{\xi}\ \binom{k+1}{\xi}\ \left(  j-\xi\right)
^{i}\right)  \ \mathbf{a}_{i}\left(  P\right)  \label{DELTA-A}%
\end{equation}
\end{lemma}

\noindent\textit{Proof. }Consider the sum
\[
\sum_{\nu=0}^{\infty}\ \left(  \sum_{i=0}^{k}\mathbf{a}_{i}\left(  P\right)
\nu^{i}\right)  \ \sum_{\mu=0}^{k+1}\ \left(  -1\right)  ^{\mu}\ \binom
{k+1}{\mu}\ q^{\,\mu+\nu}%
\]
and compute the coefficient of $q^{j}$ in its development. $_{\Box}$

\begin{theorem}
[Cohomology Group Dimensions]Let $X(\mathbb{Z}^{d},\widehat{\Delta}%
_{P})\longrightarrow U_{\tau_{P}}$ be any torus-equivariant crepant full
resolution of a $d$-dimensional standard singular representative of a
\emph{(}singular\emph{) }Gorenstein toric affine variety $U_{\sigma}$
\emph{(}as in \emph{\ref{SINR}} and in \emph{\S 3} $\mathsf{(c)}$\emph{)}.
Then the odd cohomology groups of its overlying space are trivial and the
dimension of the even ones equals\emph{:}
\begin{equation}
\text{\emph{dim}}_{\mathbb{Q}}H^{2j}\left(  X\left(  \mathbb{Z}^{d}%
,\widehat{\Delta}_{P}\right)  ;\mathbb{Q}\right)  =\deltafett_{j}\left(
P\right)  ,\ \ \forall j,\ \ 0\leq j\leq d-1\label{COH-DIM}%
\end{equation}
and is therefore \textbf{independent} of the particular choice of a basic
triangulation $\mathcal{T}$ of $P$ by means of which one constructs the fan
$\widehat{\Delta}_{P}$ $(=\widehat{\Delta}_{P}\left(  \mathcal{T}\right)  )$.
\end{theorem}

\noindent\textit{Proof. }See Batyrev-Dais \cite[Thm. 4.4., p. 909]{BD}.
$_{\Box}$

\begin{proposition}
For torus-equivariant crepant full resolutions $X(\mathbb{Z}^{d}%
,\widehat{\Delta}_{P})\rightarrow U_{\tau_{P}}$ of a $d$-dimensional standard
singular representative of a toric affine variety $U_{\sigma}$ which
is a
\emph{(}singular\emph{)} l.c.i. with $P\sim P_{\mathbf{m}}^{\left(
d\right)  }$ \emph{(}as in Thm. \emph{\ref{Nak-thm}), }the non-trivial
cohomology group dimensions of $X(\mathbb{Z}^{d},\widehat{\Delta}_{P})$ are
computable by means of the formulae \emph{(\ref{DELTA-A}), (\ref{COH-DIM}),
}and the coefficients of the Ehrhart polynomial
\begin{equation}
\mathbf{Ehr}\left(  P,\nu\right)  =\sum_{\mu_{1}=0}^{m_{1,1}\nu}\ \sum
_{\mu_{2}=0}^{m_{2,1\,}\nu+m_{2,2\,}\mu_{1}}\cdots\,\cdots\sum_{\mu_{d-1}%
=0}^{m_{d-1,1}\,\nu+m_{d-1,2}\,\mu_{1}+m_{d-1,3}\,\mu_{2}+\cdots
+m_{d-1,d-1}\,\mu_{d-2}}\ \mathbf{1}\label{EHR-POL}%
\end{equation}
which depend exlusively on the corresponding admissible
free-parameter-sequence $\mathbf{m}$\textbf{\ }defining $P_{\mathbf{m}%
}^{\left(  d\right)  }$.\newline \emph{(}Notice that $d\geq2$ and\emph{\ }by
convention\emph{: }$\mu_{0}=\nu$ for $d=2.$\emph{)}\textbf{\ }
\end{proposition}

\noindent\textit{Proof. \ }Since $P_{\mathbf{m}}^{\left(  d\right)  }$
contains always the ``origin'' of aff$(P_{\mathbf{m}}^{\left(  d\right)  })$
as one of its vertices, the $\nu$ times dilated $P_{\mathbf{m}}^{\left(
d\right)  }$ (with respect to aff$(P_{\mathbf{m}}^{\left(  d\right)  })$!) can
be described by Lemma \ref{INEQ} as\newline
\[
\nu\,P_{\mathbf{m}}^{\left(  d\right)  }=\left\{  \mathbf{x}=(x_{1}%
,..,x_{d})^{\intercal}\in\mathbb{R}^{d}\ \left|  \ x_{1}%
=1\ \text{and\emph{\ \ }}0\leq x_{j+1}\leq m_{j,1}\,\nu+\sum_{2\leq\kappa\leq
j}\right.  m_{j,\kappa}x_{\kappa},\,\forall\text{\emph{\ }}j,\,1\leq j\leq
d-1\right\}  \,.
\]
Thus, formula (\ref{EHR-POL}) expresses its canonical lattice point
enumerator. $_{\Box}$

\begin{example}
\emph{For }$d=2,3,4$\emph{, the Ehrhart polynomial of }$P_{\mathbf{m}%
}^{\left(  d\right)  }$ \emph{equals }$\mathbf{Ehr}(P_{\mathbf{m}}^{\left(
2\right)  },\nu)=\allowbreak m_{1,1}\nu+1$\emph{,}
\[
\mathbf{Ehr}(P_{\mathbf{m}}^{\left(  3\right)  },\nu)=(\frac{1}{2}%
m_{2,2}m_{1,1}^{2}+m_{2,1}m_{1,1})\,\nu^{2}\allowbreak+(m_{1,1}+\frac{1}%
{2}m_{2,2}m_{1,1}+m_{2,1})\nu+1
\]
$\emph{and\ }$ $\mathbf{Ehr}(P_{\mathbf{m}}^{\left(  4\right)  },\nu
)=(m_{3,1}m_{2,1}m_{1,1}+\frac{1}{2}m_{3,2}m_{2,1}m_{1,1}^{2}+\frac{1}%
{2}m_{3,3}\allowbreak m_{2,1}^{2}m_{1,1}+\smallskip$\newline $+\frac{1}%
{6}m_{3,3}m_{2,2}^{2}m_{1,1}^{3}+\frac{1}{2}m_{3,3}\allowbreak m_{2,2}%
m_{2,1}m_{1,1}^{2}+\frac{1}{2}m_{3,1}m_{2,2}m_{1,1}^{2}+\frac{1}{3}%
m_{3,2}\allowbreak m_{2,2}m_{1,1}^{3})\ \nu^{3}+\smallskip$\newline
$(m_{2,1}m_{1,1}+\frac{1}{2}m_{3,3}m_{2,1}^{2}+\frac{1}{2}m_{3,1}%
m_{2,2}m_{1,1}+\allowbreak\frac{1}{4}m_{3,3}m_{2,2}^{2}m_{1,1}^{2}%
+m_{3,1}m_{2,1}+\smallskip\allowbreak$\newline $+\frac{1}{2}m_{2,2}m_{1,1}%
^{2}+\frac{1}{2}m_{3,2}m_{1,1}^{2}+\frac{1}{2}m_{3,2}m_{2,2}m_{1,1}^{2}%
+\frac{1}{4}m_{3,3}\allowbreak m_{2,2}m_{1,1}^{2}+\frac{1}{2}m_{3,3}%
m_{2,1}m_{1,1}+\smallskip$\newline $m_{3,1}m_{1,1}+\allowbreak\frac{1}%
{2}m_{3,2}m_{2,1}m_{1,1}+\frac{1}{2}m_{3,3}m_{2,2}m_{2,1}m_{1,1}%
\allowbreak)\ \nu^{2}+\smallskip$\newline $+(\frac{1}{2}m_{3,2}m_{1,1}%
+m_{2,1}+m_{1,1}+\frac{1}{2}m_{3,3}m_{2,1}+\allowbreak\frac{1}{2}%
m_{2,2}m_{1,1}+\smallskip$\newline $+m_{3,1}+\frac{1}{12}m_{3,3}m_{2,2}%
^{2}m_{1,1}+\allowbreak\frac{1}{6}m_{3,2}m_{2,2}m_{1,1}+\frac{1}{4}%
m_{3,3}m_{2,2}m_{1,1})\ \nu+1\allowbreak$\emph{, respectively.}
\end{example}

\section{Extreme classes: $\left(  d,k\right)  $-hypersurfaces and
\textbf{RP}-singularities\label{EXTREME}}

\noindent Two-dimensional toric singularities are always msc-singularities.
Moreover, the underlying spaces of the Gorenstein ones (more precisely, the
standard singular representatives of them) are of the form
\begin{equation}
U_{\tau}=\text{Max-Spec}\left(  \mathbb{C}\left[  t,u,w\right]
\,/\,\left\langle t^{k}-u\ w\right\rangle \right)  \label{A-SING}%
\end{equation}
i.e., hypersurfaces depending on a free parameter $k\in\mathbb{Z}_{\geq2}$.
(These are nothing but the classically called $A_{k-1}$\textit{-singularities}%
.) Obviously, $U_{\tau}=U_{\tau_{\text{conv}\left(  \left\{  e_{1}%
,e_{1}+k\cdot e_{2}\right\}  \right)  }}$ with conv$\left(  \left\{
e_{1},e_{1}+k\cdot e_{2}\right\}  \right)  $ a lattice segment constructed by
a dilation of a unit interval by the scalar $k$. In this section, we
apply our results for two
classes (\ref{DK-SING}) and (\ref{RP-SING}) of toric msc-g.c.i.-singularities
which are direct generalizations of (\ref{A-SING}) and which are, in addition,
``extreme'', in the sense, that their corresponding Nakajima polytopes achieve
exactly the lowest and the highest bound, respectively, for the number
(\ref{BOUNDS}) of vertices/facets. Moreover, these Nakajima polytopes for both
classes are simultaneously examples for $\mathbb{H}_{d}$-compatible
polytopes.\bigskip\newline \textsf{(a) }For $k\in\mathbb{N}$, $d\in
\mathbb{Z}_{\geq2}$, let $\mathbf{s}_{k}^{\left(  d\right)  }\subset$
$\mathbf{\bar{H}}^{\left(  d\right)  }\hookrightarrow\mathbb{R}^{d}$ denote
the $\left(  d-1\right)  $-simplex
\[
\mathbf{s}_{k}^{\left(  d\right)  }:=\text{conv}\left(  \left\{  e_{1}%
,e_{1}+k\,e_{2},e_{1}+k(\,e_{2}+e_{3}),\ldots,e_{1}+k(\,e_{2}+e_{3}%
+\cdots+e_{d-1}+e_{d})\right\}  \right)
\]
being constructed by the $k$-th dilation of a basic $\left(  d-1\right)  $-simplex.

\begin{proposition}
[On $\left(  d;k\right)  $-hypersurfaces]\label{HYPERS}\emph{(i) }%
$\mathbf{s}_{k}^{\left(  d\right)  }$ is a Nakajima polytope \emph{(}w.r.t.
$\mathbb{R}^{d}$\emph{).\smallskip}\newline \emph{(ii) }For the corresponding
affine toric g.c.i.-variety we have\emph{:}
\begin{equation}
U_{\tau_{\mathbf{s}_{k}^{\left(  d\right)  }}}\cong\text{\emph{Max-Spec}%
}\left(  \mathbb{C}\left[  t,u_{1},u_{2},u_{3},\ldots,u_{d}\right]
\,/\,\left\langle t^{k}-%
{\textstyle\prod\nolimits_{j=1}^{d}}
u_{i}\right\rangle \right)  \label{DK-SING}%
\end{equation}
\emph{(}This is called, in particular,\emph{\ }$\left(  d;k\right)
$-\emph{hypersurface).\medskip\ }\newline \emph{(iii) }$(U_{\tau
_{\mathbf{s}_{k}^{\left(  d\right)  }}},\emph{orb}(\tau_{\mathbf{s}%
_{k}^{\left(  d\right)  }}))$ is a singularity \emph{(}in fact, an
msc-singularity\emph{)} if and only if $k\geq2.\smallskip\medskip$%
\newline \emph{(iv) }$\mathbf{s}_{k}^{\left(  d\right)  }$ is a $\mathbb{H}%
_{d}$-compatible polytope.\smallskip\newline \emph{(v) }\ For all
torus-equivariant crepant full desingularizations $X(\mathbb{Z}^{d}%
,\widehat{\Delta}_{\mathbf{s}_{k}^{\left(  d\right)  }})\longrightarrow
U_{\tau_{\mathbf{s}_{k}^{\left(  d\right)  }}}$ we obtain\emph{:}
\[
\text{\emph{dim}}_{\mathbb{Q}}H^{2j}\left(  X(\mathbb{Z}^{d},\widehat{\Delta
}_{\mathbf{s}_{k}^{\left(  d\right)  }});\mathbb{Q}\right)  =%
{\displaystyle\sum\limits_{i=0}^{j}}
\,\left(  -1\right)  ^{i}\,\dbinom{d}{i}\,\dbinom{k\left(  j-i\right)
+d-1}{d-1},\ \forall j,\ \ 0\leq j\leq d-1.
\]
\end{proposition}

\noindent\textit{Proof. }(i) Obviously, $\mathbf{s}_{k}^{\left(  d\right)
}=P_{\mathbf{m}}^{\left(  d\right)  }$ with $\mathbf{m}$ denoting the $\left(
\left(  d-1\right)  \times d\right)  $-matrix having entries $m_{1,1}=k$,
$m_{i,i}=1$ in its diagonal, $\forall i$, $2\leq i\leq d-1$, and zero entries
otherwise. For (ii), (iii), (iv) and (v) see Dais-Henk-Ziegler
\cite[Prop. 5.10, 6.1, and Cor. 7.4]{DHZ1}. (Notice that actually
$U_{\tau_{\mathbf{s}%
_{k}^{\left(  d\right)  }}}\cong\mathbb{C}^{d}/G\left(  d;k\right)  $ is an
abelian quotient space with $G\left(  d;k\right)  \cong\left(  \mathbb{Z\,}%
/\,k\,\mathbb{Z}\right)  ^{d-1}$). $_{\Box}\bigskip$

\noindent\textsf{(b) }Let $k_{1},k_{2},\ldots,k_{d-1}$ be a $\left(
d-1\right)  $-tuple of positive integers ($d\geq2$), and let
\[%
\begin{array}
[c]{ll}%
\mathbf{RP\medskip}\left(  k_{1},k_{2},\ldots,k_{d-1}\right)  & =\left\{
\left(  x_{1},..,x_{d}\right)  ^{\intercal}\in\mathbb{R}^{d}\ \left|
\ x_{1}=1\medskip,\ 0\leq x_{j+1}\leq k_{j},\ \forall j,\ 1\leq j\leq
d-1\right.  \right\} \\
& \ \\
\, & =\left\{  1\right\}  \times\left[  0,k_{1}\right] \times\left[  0,k_{2}\right] \times\cdots\times\left[  0,k_{d-1}\right]
\end{array}
\]
denote the $\left(  d-1\right)  $-dimensional \textit{rectangular
parallelepiped} in $\mathbf{\bar{H}}^{\left(  d\right)  }\hookrightarrow
\mathbb{R}^{d}$ having them as lengths of its edges.

\begin{proposition}
[On \textbf{RP}-singularities]\label{RPS}\emph{(i) }$\mathbf{RP}\left(
k_{1},k_{2},\ldots,k_{d-1}\right)  $ is a Nakajima polytope \emph{(}w.r.t.
$\mathbb{R}^{d}$\emph{).\medskip}\newline \emph{(ii) }For the corresponding
affine toric g.c.i.-variety we have\emph{:}
\begin{equation}
U_{\tau_{\mathbf{RP}\left(  k_{1},..,k_{d-1}\right)  }}\cong
\text{\emph{Max-Spec}}(\mathbb{C}\left[  t,u_{1},u_{2},..,u_{d-1},w_{1}%
,w_{2},..,w_{d-1}\right]  \,/\,\left\langle \{t^{k_{i}}-u_{i}w_{i}\left|
_{1\leq i\leq d-1}\right.  \}\right\rangle )\label{RP-SING}%
\end{equation}
\emph{(iii) }$(U_{\tau_{\mathbf{RP}\left(  k_{1},..,k_{d-1}\right)  }%
},\emph{orb}(\tau_{\mathbf{RP}\left(  k_{1},..,k_{d-1}\right)  }))$ is an
msc-singularity $($unless $d=2$ and $k_1=1)${}$.\medskip$\newline \emph{(iv) }$\mathbf{RP}\left(  k_{1}%
,k_{2},\ldots,k_{d-1}\right)  $ is a $\mathbb{H}_{d}$-compatible
polytope.\medskip\newline \emph{(v) }For all torus-equivariant crepant full
desingularizations
\[
X(\mathbb{Z}^{d},\widehat{\Delta}_{\mathbf{RP}\left(  k_{1},k_{2}%
,\ldots,k_{d-1}\right)  })\longrightarrow U_{\tau_{\mathbf{RP}\left(
k_{1},k_{2},\ldots,k_{d-1}\right)  }}%
\]
we obtain for all $j$, $0\leq j\leq d-1$\emph{:}
\[
\text{\emph{dim}}_{\mathbb{Q}}H^{2j}\left(  X(\mathbb{Z}^{d},\widehat{\Delta
}_{\mathbf{RP}\left(  k_{1},..,k_{d-1}\right)  });\mathbb{Q}\right)
=\sum_{i=0}^{d-1}\mathbf{\ }\left[  \sum_{\xi=0}^{j}\ \left(  -1\right)
^{\xi}\ \binom{d}{\xi}\ \left(  j-\xi\right)  ^{i}\right]  \ \mathfrak{s}%
_{i}\left(  k_{1},k_{2},\ldots,k_{d-1}\right)  ,
\]
where
\[
\mathfrak{s}_{0}\left(  k_{1},..,k_{d-1}\right)  =1,\ \ \mathfrak{s}_{i}\left(
k_{1},..,k_{d-1}\right)  =\sum_{1\leq\mu_{1}<\mu_{2}<\cdots<\mu_{i}\leq
d-1}\ k_{\mu_{1}}\cdot k_{\mu_{2}}\cdot\cdots\cdot k_{\mu_{i}},\ \forall
i,\ \ 1\leq i\leq d-1,
\]
are the elementary symmetric polynomials w.r.t. the variables $k_{1}%
,k_{2},\ldots,k_{d-1}$.
\end{proposition}

\noindent \textit{Proof. }(i) $\mathbf{RP}\left( k_{1},k_{2},\ldots
,k_{d-1}\right) $ equals $P_{\mathbf{m}}^{\left( d\right) }$ with $\mathbf{m}
$ denoting the $\left( \left( d-1\right) \times d\right) $-matrix with
entries $m_{i,1}=k_{i}$ in its first column, for all $i$, $1\leq i\leq d-1$,
and zero entries otherwise. Moreover, it has $2^{d-1}$ vertices; namely
\begin{equation*}
\text{vert}(\mathbf{RP}\left( k_{1},..,k_{d-1}\right) )=
\left\{e_1 + \varepsilon_1 \cdot k_{1}\cdot e_{2} + \varepsilon_2 \cdot
k_{2}\cdot e_{3} + \cdots + \varepsilon_{d-1} \cdot k_{d-1} \cdot
e_{d} \medspace | \medspace \varepsilon_1,\ldots,\varepsilon_{d-1} \in
\{0,1\} \right\}.
\end{equation*}
\\
(ii) As it was pointed out by Nakajima \cite[p. 92]{Nakajima}, the set
\begin{equation*}
\left\{ e_{1}^{\vee },e_{2}^{\vee },\ldots ,e_{d-1}^{\vee },e_{d}^{\vee
},k_{1}\cdot e_{1}^{\vee }-e_{2}^{\vee },k_{2}\cdot e_{1}^{\vee
}-e_{3}^{\vee },\ldots ,k_{d-1}\cdot e_{1}^{\vee }-e_{d}^{\vee }\right\}
\end{equation*}
forms a system of generators for the monoid $\tau _{\mathbf{RP}\left(
k_{1},..,k_{d-1}\right) }^{\vee }\cap (\mathbb{Z}^{d})^{\vee }$. Using (\ref
{Hilbbasis}) it is easy to see that the above set is exactly the Hilbert
basis of $\tau _{\mathbf{RP}\left( k_{1},..,k_{d-1}\right) }^{\vee }$ w.r.t.
$(\mathbb{Z}^{d})^{\vee }$. Hence,
\begin{equation*}
\left\{ \mathbf{e}\left( e_{1}^{\vee }\right) ,\mathbf{e}\left( e_{2}^{\vee
}\right) ,\ldots ,\mathbf{e}\left( e_{d-1}^{\vee }\right) ,\mathbf{e}\left(
e_{d}^{\vee }\right) ,\mathbf{e}\left( k_{1}\cdot e_{1}^{\vee }-e_{2}^{\vee
}\right) ,\mathbf{e}\left( k_{2}\cdot e_{1}^{\vee }-e_{3}^{\vee }\right)
,\ldots ,\mathbf{e}\left( k_{d-1}\cdot e_{1}^{\vee }-e_{d}^{\vee }\right)
\right\}
\end{equation*}
generates $\mathbb{C}[\tau _{\mathbf{RP}\left( k_{1},..,k_{d-1}\right)
}^{\vee }\cap (\mathbb{Z}^{d})^{\vee }]$, and the affine toric variety $%
U_{\tau _{\mathbf{RP}\left( k_{1},..,k_{d-1}\right) }}$ has embedding
dimension $2d-1$ (and is, in particular, a g.c.i. of $d-1$ \textit{binomials}
by (i), Thm. \ref{Nak-thm}, Rem. \ref{REMNAK}(ii), and Thm. \ref{EMB}). The
map
\begin{equation*}
\theta :\mathbb{C}\left[ t,u_{1},u_{2},..,u_{d-1},w_{1},w_{2},..,w_{d-1}%
\right] \longrightarrow \mathbb{C}[\tau _{\mathbf{RP}\left(
k_{1},..,k_{d-1}\right) }^{\vee }\cap (\mathbb{Z}^{d})^{\vee }]
\end{equation*}
defined by $\theta \left( t\right) :=\mathbf{e}\left( e_{1}^{\vee }\right) $%
, $\theta \left( u_{i}\right) :=\mathbf{e}\left( e_{i+1}^{\vee }\right) $, $%
\theta \left( w_{i}\right) :=\mathbf{e}\left( k_{i}\cdot e_{1}^{\vee
}-e_{i+1}^{\vee }\right) $, $\forall i$, $1\leq i\leq d-1$, is a $\mathbb{C}$%
-algebra epimorphism. It suffices to show that Ker$\left( \theta \right) =I$
with $I:=\left\langle \{t^{k_{i}}-u_{i}w_{i}\ \left| \ 1\leq i\leq
d-1\right. \}\right\rangle $. For $d=2$ this is obvious. We shall hereafter
assume that $d\geq 3$. \medskip \newline
$\bullet $ \textit{Claim }A. This ideal is contained in the kernel of $%
\theta $, i.e., $I\subseteq $ Ker$\left( \theta \right) .\medskip $\newline
$\bullet $ \textit{Proof of Claim }A. Consider the lattice $\Lambda $ of the
defining binomial equations of $U_{\tau _{\mathbf{RP}\left(
k_{1},..,k_{d-1}\right) }}\hookrightarrow \mathbb{C}^{2d-1}$,
\begin{equation*}
\Lambda =\left\{ \left( a_{1},a_{2},\ldots ,a_{2d-1}\right) \in \mathbb{Z}%
^{2d-1}\ \left| \ \sum_{i=1}^{d}\,a_{i}\,e_{i}^{\vee
}+\sum_{i=d+1}^{2d-1}\,a_{i}\,\left( k_{i-d}\cdot e_{1}^{\vee
}-e_{i-d+1}^{\vee }\right) \right. =0\right\} \ .
\end{equation*}
Since the extra relations can be written in the form
\begin{equation*}
a_{1}+\sum_{i=d+1}^{2d-1}\,a_{i}\,k_{i-d}=0,\ \ \ \ \ \ \
a_{j}-a_{d+j-1}=0,\ \ \ \forall j,\ \ \ 2\leq j\leq d,
\end{equation*}
setting $\xi _{j-1}:=-a_{j}=-a_{d+j-1}$, for all $j$, $2\leq j\leq d$, as
auxiliary parameters, we may express every point of $\Lambda $ as follows
\begin{align*}
\left( a_{1},a_{2},\ldots ,a_{2d-1}\right) &
=(\sum\nolimits_{i=1}^{d-1}\,\xi _{i}\,k_{i},-\xi _{1},-\xi _{2},\ldots
,-\xi _{d-1},-\xi _{1},-\xi _{2},\ldots ,-\xi _{d-1})=\smallskip  \\
& =\sum\nolimits_{i=1}^{d-1}\,\xi _{i}\,\,(k_{i}\cdot e_{1}^{\vee
}-e_{i+1}^{\vee }-e_{d+i}^{\vee }).
\end{align*}
Now these $d-1$ vectors are also $\mathbb{Z}$-linearly independent. So they
constitute a $\mathbb{Z}$-basis of $\Lambda $, and
\begin{equation*}
k_{i}\cdot e_{1}^{\vee }-e_{i+1}^{\vee }-e_{d+i}^{\vee }=k_{i}\cdot
e_{1}^{\vee }-(e_{i+1}^{\vee }+e_{d+i}^{\vee })
\end{equation*}
is the difference of two vectors with non-negative coordinates having
disjoint support for all indices $i$, $1\leq i\leq d-1$. Hence, $I\subseteq $
Ker$\left( \theta \right) $ by \cite[4.3-4.4, p. 32]{Sturmfels}.\medskip\
\newline
$\bullet $ \textit{Claim }B. The opposite inclusion $I\supseteq $ Ker$\left(
\theta \right) $ is true too$.\medskip $\newline
$\bullet $ \textit{Proof of Claim }B. For every $\kappa \in \mathbb{N}$, $%
1\leq \kappa \leq d-1$, and every subset of indices $2\leq
i_{1}<i_{2}<\cdots <i_{\kappa }\leq d$ of length $\kappa $, we define the
cone
\begin{equation*}
\begin{array}{ll}
C_{i_{1},i_{2},\ldots ,i_{\kappa }} & :=\left\{ \left( \mathbf{y}_{1},%
\mathbf{y}_{2},\ldots ,\mathbf{y}_{d}\right) \in \tau _{\mathbf{RP}\left(
k_{1},..,k_{d-1}\right) }^{\vee }\,\left|
\begin{array}{l}
\mathbf{y}_{1}\geq 0,\ \mathbf{y}_{i_{1}}\geq 0,\ \mathbf{y}_{i_{2}}\geq 0,\
\ldots ,\ \mathbf{y}_{i_{\kappa }}\geq 0\medskip  \\
\text{and \ }\mathbf{y}_{j}\leq 0,\ \forall j,\ j\in \{2,\ldots ,d\}\mathbb{r%
}\{i_{1},\ldots ,i_{\kappa }\}
\end{array}
\right. \right\}  \\
\  & \  \\
\  & =\text{ pos}\left( \left\{ e_{1}^{\vee }\right\} \cup \left\{
e_{i_{1}}^{\vee },\ldots ,e_{i_{\kappa }}^{\vee }\right\} \cup \left\{
k_{j-1}\cdot e_{1}^{\vee }-e_{j}^{\vee }\ \left| \ j\right. \in \{2,\ldots
,d\}\mathbb{r}\{i_{1},\ldots ,i_{\kappa }\}\right\} \right) \,.
\end{array}
\end{equation*}
Then
\begin{equation*}
\tau _{\mathbf{RP}\left( k_{1},..,k_{d-1}\right) }^{\vee }=\bigcup_{\kappa
=1}^{d-1}\ \bigcup_{2\leq i_{1}<i_{2}<\cdots <i_{\kappa }\leq d}\
C_{i_{1},i_{2},\ldots ,i_{\kappa }}
\end{equation*}
is a fan-subdivision of $\tau _{\mathbf{RP}\left( k_{1},..,k_{d-1}\right)
}^{\vee }$ into $2^{d-1}$  s.c.p. (and actually basic) cones. Now since
\begin{equation*}
\mathbf{Hilb}_{(\mathbb{Z}^{d})^{\vee }}(\tau _{\mathbf{RP}\left(
k_{1},..,k_{d-1}\right) }^{\vee })=\left\{ e_{1}^{\vee },e_{2}^{\vee
},\ldots ,e_{d}^{\vee },k_{1}\cdot e_{1}^{\vee }-e_{2}^{\vee },k_{2}\cdot
e_{1}^{\vee }-e_{3}^{\vee },\ldots ,k_{d-1}\cdot e_{1}^{\vee }-e_{d}^{\vee
}\right\} ,
\end{equation*}
each lattice point $m\in \tau _{\mathbf{RP}\left( k_{1},..,k_{d-1}\right)
}^{\vee }\cap (\mathbb{Z}^{d})^{\vee }$ can be written as a linear
combination
\begin{equation}
m=\lambda \cdot e_{1}^{\vee }+\sum_{i=1}^{d-1}\,\mu _{i}\cdot e_{i+1}^{\vee
}+\sum_{j=1}^{d-1}\,\nu _{j}\cdot (k_{j}\cdot e_{1}^{\vee }-e_{j+1}^{\vee })
\label{LIN-COMB}
\end{equation}
for \textit{uniquely} determined coefficients $(\lambda ,\mu _{1},\ldots
,\mu _{d-1},\nu _{1},\ldots ,\nu _{d-1})\in (\mathbb{Z}_{\geq 0})^{2d-1}$
which have to satisfy the extra conditions
\begin{equation}
\mu _{1}\cdot \nu _{1}=\mu _{2}\cdot \nu _{2}=\cdots =\mu _{d-1}\cdot \nu
_{d-1}=0  \label{EXTRA-COND}
\end{equation}
because $m$ necessarily belongs to a cone of the form $C_{i_{1},i_{2},\ldots
,i_{\kappa }}$. For the rest of the proof it is enough to make use of an
elegant trick due to Ishida (see \cite[p. 143]{Ishida}). We define a
homomorphism
\begin{equation*}
g:\mathbb{C}[\tau _{\mathbf{RP}\left( k_{1},..,k_{d-1}\right) }^{\vee }\cap (%
\mathbb{Z}^{d})^{\vee }]\longrightarrow \mathbb{C}\left[
t,u_{1},u_{2},..,u_{d-1},w_{1},w_{2},..,w_{d-1}\right]
\end{equation*}
of $\mathbb{C}$-vector spaces by mapping the character of any $m$ as in (\ref
{LIN-COMB}) onto
\begin{equation*}
\mathbb{C}[\tau _{\mathbf{RP}\left( k_{1},..,k_{d-1}\right) }^{\vee }\cap (%
\mathbb{Z}^{d})^{\vee }]\ni \mathbf{e}(m)\longmapsto g(\mathbf{e}%
(m)):=t^{\lambda }\,u_{1}^{\mu _{1}}\,\cdots \,u_{d-1}^{\mu
_{d-1}}\,w_{1}^{\nu _{1}}\,\cdots \,w_{d-1}^{\nu _{d-1}}\ .
\end{equation*}
$g$ is obviously a section of $\theta $ (i.e., $\theta \circ g=$ Id). For
proving $I\supseteq $ Ker$\left( \theta \right) $ it is therefore sufficient
to show that
\begin{equation}
\mathbb{C}\left[ t,u_{1},\ldots ,u_{d-1},w_{1},\ldots ,w_{d-1}\right] =g(%
\mathbb{C}[\tau _{\mathbf{RP}\left( k_{1},..,k_{d-1}\right) }^{\vee }\cap (%
\mathbb{Z}^{d})^{\vee }])+I\ .  \label{TO-PROVE}
\end{equation}
Let $\phi =t^{p}\,u_{1}^{q_{1}}\,\cdots
\,u_{d-1}^{q_{d-1}}\,w_{1}^{r_{1}}\,\cdots \,w_{d-1}^{r_{d-1}}$ denote an
arbitrary monomial in $\mathbb{C}\left[ t,u_{1},\ldots ,u_{d-1},w_{1},\ldots
,w_{d-1}\right] $. To conclude (\ref{TO-PROVE}) we shall prove that
\begin{equation}
\phi \in g(\mathbb{C}[\tau _{\mathbf{RP}\left( k_{1},..,k_{d-1}\right)
}^{\vee }\cap (\mathbb{Z}^{d})^{\vee }])+I\ .\   \label{PROPERTY}
\end{equation}
Define $\mathfrak{A}_{\phi }:=\{i\in \{1,\ldots ,d-1\}\ \left| \ q_{i}\cdot
r_{i}>0\right. \}$ and  $\beta _{\phi }:=\#(\mathfrak{A}_{\phi }).$ If $\beta
_{\phi }=0$, then $\phi =g(\mathbf{e}(m))$, for some $m\in \tau _{\mathbf{RP}%
\left( k_{1},..,k_{d-1}\right) }^{\vee }\cap (\mathbb{Z}^{d})^{\vee }$ (by (%
\ref{LIN-COMB}) and (\ref{EXTRA-COND})). Let now $\beta _{\phi }\in
\{1,..,d-1\}$ and consider an index $i_{0}\in \mathfrak{A}_{\phi }$. Suppose
that property (\ref{PROPERTY}) is false for $\phi .$ Without loss of
generality, we may further assume that the product $q_{i_{0}}\cdot r_{i_{0}}$
of the degrees of $\phi $ in the variables $u_{i_{0}}$ and $w_{i_{0}}$ is
chosen to be \textit{minimal} with respect to the violation of (\ref
{PROPERTY}).\medskip\ \newline
Then
\begin{equation*}
t^{p+k_{i_{0}}}\,\,u_{1}^{q_{1}}\,\cdots
\,u_{i_{0}-1}^{q_{i_{0}-1}}\,u_{i_{0}}^{q_{i_{0}}-1}%
\,u_{i_{0}+1}^{q_{i_{0}+1}}\,\cdots
\,\,u_{d-1}^{q_{d-1}}\,w_{1}^{r_{1}}\,\cdots
\,w_{i_{0}-1}^{r_{i_{0}-1}}\,w_{i_{0}}^{r_{i_{0}}-1}%
\,w_{i_{0}+1}^{r_{i_{0}+1}}\,\cdots \,\,w_{d-1}^{r_{d-1}\noindent }-\phi =
\end{equation*}
\begin{equation*}
=t^{p}\,u_{1}^{q_{1}}\,\cdots
\,u_{i_{0}-1}^{q_{i_{0}-1}}\,u_{i_{0}}^{q_{i_{0}}-1}%
\,u_{i_{0}+1}^{q_{i_{0}+1}}\,\cdots
\,\,u_{d-1}^{q_{d-1}}\,w_{1}^{r_{1}}\,\cdots
\,w_{i_{0}-1}^{r_{i_{0}-1}}\,w_{i_{0}}^{r_{i_{0}}-1}%
\,w_{i_{0}+1}^{r_{i_{0}+1}}\,\cdots \,\,w_{d-1}^{r_{d-1}\noindent }\,\left(
t^{k_{i_{0}}}-u_{i_{0}}w_{i_{0}}\right) \in I\,.
\end{equation*}
This implies that
\begin{equation*}
t^{p+k_{i_{0}}}\,\,u_{1}^{q_{1}}\,\cdots
\,u_{i_{0}-1}^{q_{i_{0}-1}}\,u_{i_{0}}^{q_{i_{0}}-1}%
\,u_{i_{0}+1}^{q_{i_{0}+1}}\,\cdots
\,\,u_{d-1}^{q_{d-1}}\,w_{1}^{r_{1}}\,\cdots
\,w_{i_{0}-1}^{r_{i_{0}-1}}\,w_{i_{0}}^{r_{i_{0}}-1}%
\,w_{i_{0}+1}^{r_{i_{0}+1}}\,\cdots \,\,w_{d-1}^{r_{d-1}\noindent }
\end{equation*}
does not belong to $g(\mathbb{C}[\tau _{\mathbf{RP}\left(
k_{1},..,k_{d-1}\right) }^{\vee }\cap (\mathbb{Z}^{d})^{\vee }])+I$,
contradicting the minimality assumption for $q_{i_{0}}\cdot r_{i_{0}}$.
Hence, (\ref{PROPERTY}) is always true.\medskip

\noindent (iii) and (iv) are obvious. Finally, (v) follows from the
determination of the coefficients of the Ehrhart polynomial
\begin{equation*}
\mathbf{Ehr}\left( \mathbf{RP}\left( k_{1},k_{2},\ldots ,k_{d-1}\right) ;\nu
\right) =\prod_{i=1}^{d-1}\ \left( k_{i}\ \nu +1\right) \ ,
\end{equation*}
combined with the formulae (\ref{DELTA-A}) and (\ref{COH-DIM}). $_{\Box }$

\begin{remark}
\emph{As both toric g.c.i.-varieties (\ref{DK-SING}) and (\ref{RP-SING}) are
constructible by means of \ }$\mathbb{H}_{d}$\emph{-compatible Nakajima
polytopes, the most natural choice of a crepant birational morphism} \emph{to
desingularize them is } $f_{\left.  \mathcal{T}_{\mathbb{H}_{d}}\right|
\text{\emph{restr.}}}.$ \emph{\ For this choice the precise nature of the
occuring exceptional prime divisors is known by Thm. \ref{ZONOS}. \smallskip}
\end{remark}

\noindent\textsf{(c) }At the end of the paper we devote a few words to
\textit{non-l.c.i.'s}: In complete analogy to the case of non-l.c.i.\
Gorenstein abelian
quotient spaces (cf. \cite{DHH,DH}), we expect that also the
underlying spaces of toric non-l.c.i.\
Gorenstein singularities will be \textit{only
rarely\/} overall resolvable by crepant birational morphisms. Let us
nevertheless give two examples of non-Nakajima polytopes admitting
b.c.-triangulations.

\begin{example} \label{DELPEZZO}
\emph{Let }$d$ \emph{be an odd integer }$\geq3$\emph{. Define the }$\left(
d-1\right)  $\emph{-dimensional lattice polytope}
\begin{align*}
Q &=
\text{\emph{conv}}\left(  \left\{  e_{1}\pm e_{2},e_{1}\pm e_{3}%
,\ldots,e_{1}\pm e_{d-1},e_{1}\pm e_{d},e_{1}\pm(%
{\textstyle\sum\nolimits_{j=2}^{d}}
e_{j})\right\}  \right)
\\
&= \left\{ \mathbf{x} =\left( x_{1},\ldots,x_{d}\right)^{\intercal}%
\in \mathbb{R}^{d}\ | \ x_1 = 1 , -1 \le x_1 + \sum_{i=2}^d
  \varepsilon_i x_i \le 1 , \forall \varepsilon_i \in \{\pm1\}: 
  \sum_{i=2}^d \varepsilon_i \in \{-1,0,1\} \right\}   .
\end{align*}
\emph{Since}
\[
\#(\left\{  \text{\emph{facets of \ }}Q\right\}  ) =
d \cdot \tbinom 
{d-1\smallskip}{\frac{1}{2}(d-1)} = \frac{d\,!}{\left[  (\frac
{1}{2}\left(  d-1\right)  )!\right]  ^{2}}>2\,\left(  d-1\right)  ,
\]
$Q \subset\mathbf{\bar{H}}^{\left(  d\right)  }\hookrightarrow\mathbb{R}%
^{d}$ \emph{cannot be lattice equivalent to a Nakajima polytope (by
(\ref{BOUNDS})). On the other hand, }$Q$ \emph{is a non-simplex, Fano
polytope, and the Gorenstein non-l.c.i., non-quotient, msc-singularity
}$(U_{\tau_{Q}},$\emph{orb}$\left(  \tau_{Q}\right)  )$ \emph{can therefore be
overall resolved by a crepant projective birational morphism (by Prop.
\ref{FANOS}).}
\end{example}

\begin{example}
\emph{The zonotope }$\mathcal{Z}^{(d)}$\emph{ defined earlier in the 
  proof of Theorem \ref{ZONOS} is
  }$\mathbb{H}_{d}$\emph{-compatible, but for }$d \ge 3$ \emph{, it has }$ 
d (d-1) > 2(d-1)$\emph{ facets.}
\end{example}

\noindent To create more examples of non-Nakajima polytopes having
b.c.-triangulations, one may start with $\mathcal{Z}^{(d)}$ or with a
$Q$ as above in \ref{DELPEZZO}, and consider
``joins'' of it with further (finitely many) Nakajima polytopes, or,
alternatively, combine or mix all those with suitably triangulated dilations
of basic simplices. (For the ``good'' behaviour of joins and dilations under
b.c.-triangulations, see \cite[\S3 and \S6]{DHZ1}).


\begin{thebibliography}{99}
\bibitem{Bat1}\textsc{Batyrev V.V.}: \textit{Dual polyhedra and mirror
symmetry for Calabi-Yau hypersurfaces in toric varieties}, Jour. of Algebraic
Geometry \textbf{3}, (1994), 493-535.

\bibitem {Bat2}\textsc{Batyrev V.V.}: \textit{Stringy Hodge numbers of
varieties with Gorenstein canonical singularities}, alg-geom/9711008; to
appear in the Proc. of Taniguchi Symposium 1997: ``Integrable Systems and
Algebraic Geometry'', Kobe/Kyoto.

\bibitem {Bat3}\textsc{Batyrev V.V.}: \textit{Non-Archimedian integrals and
stringy Euler numbers of log-terminal pairs}, alg-geom/9803071; to appear in
the Jour. of the European Math. Soc. \textbf{1}, (1999).

\bibitem {BB}\textsc{Batyrev V.V., Borisov L.A.}: \textit{Mirror duality and
string-theoretic Hodge numbers}, Inventiones Math. \textbf{126}, (1996), 183-203.

\bibitem {BD}\textsc{Batyrev V.V., Dais D.I.}: \textit{Strong McKay
correspondence, string-theoretic Hodge numbers and mirror symmetry}, Topology
\textbf{35}, (1996), 901-929.

\bibitem {BFS}\textsc{Billera L.J., Filliman P., Sturmfels B.}:
\textit{Construction and complexity of secondary polytopes,} Advances in Math.
\textbf{83}, (1990), 155-179.

\bibitem {Borisov1}\textsc{Borisov L.A.}: \textit{String cohomology of a
toroidal singularity}, preprint, alg-geom/9802052.

\bibitem {Borisov2}\textsc{Borisov L.A.}: \textit{Vertex algebras and mirror
symmetry}, preprint, alg-geom/9809094.

\bibitem {BGT}\textsc{Bruns W., Gubeladze J., Trung N. V.}: \textit{Normal
polytopes, triangulations and Koszul algebras}, Jour. f\"{u}r die reine und
ang. Math. \textbf{485}, (1997), 123-160.

\bibitem {DHH}\textsc{Dais D.I., Haus U.-U., Henk M.}: \textit{On crepant
resolutions of }$2$\textit{-parameter series of} \textit{Gorenstein cyclic
quotient singularities}, Results in Math. \textbf{33}, (1998), 208-265.

\bibitem {DH}\textsc{Dais D.I., Henk M.}: \textit{On a series of Gorenstein
cyclic quotient singularities admitting a unique projective crepant
resolution}, alg-geom/9803094; to appear in: ``Combinatorial Convex Geometry
and Toric Varieties'', (ed. by G.Ewald \& B.Teissier), Birkh\"{a}user.

\bibitem {DHZ1}\textsc{Dais D.I., Henk M., Ziegler G.M.}: \textit{All abelian
quotient c.i.-singularities admit projective crepant resolutions in all
dimensions}, Advances in Math. \textbf{139}, (1998), 192-239.

\bibitem {DHZ2}\textsc{Dais D.I., Henk M., Ziegler G.M.}: \textit{On the
existence of crepant resolutions of Gorenstein abelian quotient singularities
in dimensions }$\geq4$, in preparation.

\bibitem {DELOERA}\textsc{de Loera J.A.}: \textit{Triangulations of Polytopes
and Computational Algebra}, PhD Thesis, Cornell University, (1995). [The
author's computer algebra program \textsc{puntos} which determines coherent
triangulations of the convex hull of point configurations is available via
anonymous ftp at \texttt{ftp://geom.umn.edu}, directory \texttt{priv/deloera}.]

\bibitem {DL-H-S-S}\textsc{de Loera J.A., Ho}\c{s}\textsc{ten S., Santos F.,
Sturmfels B.}: \textit{The polytope of all triangulations of a point
configuration}, Documenta Math. J. DMV \textbf{1}, (1996), 103-119.

\bibitem {Ewald}\textsc{Ewald G.}: \textit{Combinatorial Convexity and
Algebraic Geometry}, Graduate Texts in Mathematics, Vol. \textbf{168},
Springer-Verlag, (1996).

\bibitem {Fulton}\textsc{Fulton W.}: \textit{Introduction to Toric Varieties},
Annals of Math. Studies, Vol. \textbf{131}, Princeton University Press, (1993).

\bibitem {GKZ}\textsc{Gelfand I.M., Kapranov M.M., Zelevinsky A.V.}:
\textit{Discriminants, Resultants and Multidimensional Determinants},
Birkh\"{a}user, (1994).

\bibitem {Gro1}\textsc{Grothendieck A.}: \textit{Techniques de construction en
g\'{e}om\'{e}trie analytique, }S\'{e}minaire H. Cartan, \'{E}.N.S., 1960/61,
Exp. 13.

\bibitem {Gro2}\textsc{Grothendieck A.}: \textit{G\'{e}om\'{e}trie
alg\'{e}brique et g\'{e}om\'{e}trie analytique}, (Notes by M.Raynaud). In :
\textit{Rev\^{e}tements \'{e}tales et Groupe Fondemental }(SGA 1), Exp. 12,
Lecture Notes in Math., Vol. \textbf{224}, Springer-Verlag, (1971), pp. 311-343.

\bibitem {Hi-O}\textsc{Hibi T., Ohsugi H.}: \textit{A normal 
$\left(0,1\right)  $-polytope none of whose regular triangulations is
unimodular}, Discrete \ \& Comp.\ Geometry \textbf{21} (1999), 201-204.

\bibitem {Hochster}\textsc{Hochster M.}: \textit{Rings of invariants of tori,
Cohen-Macaulay rings generated by monomials, and polytopes}, Annals of Math.
\textbf{96}, (1972), 318-337.

\bibitem {Ishida}\textsc{Ishida M.-N.}: \textit{Torus embeddings and dualizing
complexes}, T\^{o}hoku Math. Jour. \textbf{32}, (1980), 111-146.

\bibitem {Ito-Reid}\textsc{Ito Y., Reid M.}: \textit{The McKay correspondence
for finite subgroups of }SL$\left(  3,\mathbb{C}\right)  $. In : ``Higher
Dimensional Complex Varieties'', Proceedings of the International Conference
held in Trento, Italy, June 15-24, 1994, (M.Andreatta \& Th.Peternell, eds.);
Walter de Gruyter, (1996), 221-240.

\bibitem {KKMS}\textsc{Kempf G., Knudsen F., Mumford D., Saint-Donat D.}:
\textit{Toroidal Embeddings I}, Lecture Notes in Mathematics, Vol.
\textbf{339}, Springer-Verlag, (1973).

\bibitem {Kontsevich}\textsc{Kontsevich M.}: \textit{p-adic integrals, loop
spaces and generalized McKay's correspondence}, two talks given at E.N.S.
(Paris) and M.P.I. (Bonn), (February-March 1996).

\bibitem {Kunz}\textsc{Kunz E.}: \textit{Introduction to Commutative Algebra
and Algebraic Geometry}, Birkh\"{a}user, (1985).

\bibitem {Lee1}\textsc{Lee C.W.}: \textit{Regular triangulations of convex
polytopes}. In: ``Applied Geometry and Discrete Mathematics-The Victor Klee
Festschrift'', (P.Gritzmann \& B.Sturmfels, eds.), DIMACS Series in Discrete
Math. and Theoretical Comp. Science, Vol. \textbf{4}, A.M.S., (1991), 443-456.

\bibitem {Lee2}\textsc{Lee C.W.}: \textit{Subdivisions and triangulations of
polytopes}, in: ``Discrete and Computational Geometry'', (J.E.Goodman \&
J.O'Rourke, eds.), CRC-Press, New York, (1997), pp. 271-290.

\bibitem {Matsumura}\textsc{Matsumura H.}: \textit{Commutative Ring Theory},
Cambridge Studies in Adv. Math., Vol. \textbf{8}, Cambridge University Press, (1986).

\bibitem {Mo-Ste}\textsc{Morrison D.R., Stevens G.}: \textit{Terminal quotient
singularities in dimension three and four}, Proc. A.M.S.\textbf{\ 90}, (1984), 15-20.

\bibitem {Nakajima}\textsc{Nakajima H.}: \textit{Affine torus embeddings which
are complete intersections}, T\^{o}hoku Math. Jour. \textbf{38}, (1986), 85-98.

\bibitem {Oda}\textsc{Oda T.}:\textit{\ Convex Bodies and Algebraic Geometry.
An Introduction to the theory of toric varieties}, Ergebnisse der Mathematik
und ihrer Grenzgebiete, 3. Folge, Bd. \textbf{15}, Springer-Verlag, (1988).

\bibitem {OP}\textsc{Oda T., Park H.S.}: \textit{Linear Gale transforms and
GKZ-decompositions}, T\^{o}hoku Math. Jour. \textbf{43}, (1991), 375-399.

\bibitem {Reid1}\textsc{Reid M.}: \textit{Canonical threefolds}, Journ\'{e}e
de G\'{e}om\'{e}trie Alg\'{e}brique d'Angers, (A. Beauville, ed.), Sijthoff
and Noordhoff, Alphen aan den Rijn, (1980), 273-310.

\bibitem {Reid2}\textsc{Reid M.}: \textit{Decompositions of toric morphisms}.
In: ``Arithmetic and Geometry II'', (M.Artin \& J.Tate, eds.), Progress in
Math. \textbf{36}, Birkh\"{a}user, (1983), 395-418.

\bibitem {Reid3}\textsc{Reid M.}: \textit{\ McKay correspondence}, preprint, alg-geom/9702016.

\bibitem {Serre}\textsc{Serre J.}-\textsc{P.}: \textit{G\'{e}om\'{e}trie
alg\'{e}brique et g\'{e}om\'{e}trie analytique}, Ann. Inst. Fourier
\textbf{6}, (1956), 1-42.

\bibitem {Stanley1}\textsc{Stanley R.P.}: \textit{Hilbert functions of graded
algebras}, Advances in Math. \textbf{28}, (1978), 57-81.

\bibitem {Stanley2}\textsc{Stanley R.P.}: \textit{Decompositions of rational
convex polytopes}, Annals of Discrete Math. \textbf{6}, (1980), 333-342.

\bibitem {Stanley3}\textsc{Stanley R.P.}: \textit{Enumerative Combinatorics},
Vol. \textbf{I}, Wadsworth \& Brooks/Cole Math. Series, (1986); second
printing: Cambridge University Press, (1997).

\bibitem {Sturmfels}\textsc{Sturmfels B.}: \textit{Gr\"{o}bner Bases and
Convex Polytopes}, University Lecture Series, Vol. \textbf{8}, A.M.S., (1996).

\bibitem {Watanabe}\textsc{Watanabe K.}: \textit{Invariant subrings which are
complete intersections I, (Invariant subrings of finite Abelian groups)},
Nagoya Math. Jour. \textbf{77}, (1980), 89-98.

\bibitem {Ziegler}\textsc{Ziegler G.M.}: \textit{Lectures on Polytopes},
Graduate Texts in Mathematics, Vol. \textbf{152}, Springer-Verlag, (1995).
\end{thebibliography}
\end{document}